\input amstex
\documentstyle {amsppt}

\pagewidth{32pc} 
\pageheight{45pc} 
\mag=1200
\baselineskip=15 pt

\hfuzz=5pt
\topmatter
\NoRunningHeads 
\title discrete localities I 
\endtitle
\author Andrew Chermak and Alex Gonzales 
\endauthor 

\address Mathematics Department, Kansas State University, Manhattan, KS 66506 (USA) 
\endaddress 

\email chermak\@math.ksu.edu 
\endemail 

\address Institut Arraona, 
43 Praga Street, 
08207 Sabadell (Spain)
\endaddress 

\email agondem\@gmail.com
\endemail 

\date
March 2022 
\enddate 

\endtopmatter

\redefine\pce{\preccurlyeq} 

\define\w{\widetilde}
\define\wh{\widehat}
\redefine\norm{\trianglelefteq}

\redefine\bar{\overline}

\redefine\maps{\mapsto}
\redefine\i{^{-1}}

\redefine\l{\lambda}
\redefine\s{\sigma}
\redefine\a{\alpha}
\redefine\b{\beta}
\redefine\d{\delta}
\redefine\g{\gamma}

\redefine\r{\rho}

\redefine\G{\Gamma}

\redefine\L{\Lambda}

\redefine\<{\langle}
\redefine\>{\rangle}

\redefine\ca{\Cal}

\redefine\D{\Delta}

\redefine\sub{\subseteq}

\redefine\nset{\emptyset} 

\redefine\1{\bold 1} 

\redefine\up{\uparrow} 
\redefine\down{\downarrow}

\redefine\bull{\bullet}

\redefine\pce{\preccurlyeq}

\document

\vskip .2in 
\noindent 
{\bf Introduction} 
\vskip .1in 
 Finite groups, Lie groups, and linear algebraic groups have 
various structural properties in common, but these properties tend to be deduced in different 
ways. For example; in finite group theory one has the notion of the generalized Fitting subgroup 
$F^*(G)=E(G)F(G)$ of a group $G$, where $E(G)$ plays much the same role for finite groups as reductive 
Lie groups or reductive algebraic groups play in their respective categories. 
The structure of centralizers of elements in simple groups belonging to any of these three 
classes reflects a similar convergence after the fact. 

There remains some mystery as to why this should be so, even though most of the finite simple groups are - 
as we know -  derivable from algebraic groups or from Lie groups by the methods pioneered by 
Chevalley. The alternating groups can be construed as being degenerate cases of algebraic or Lie groups 
via the theory of buildings, so one is left with only a small number (twenty-six) of exceptional 
cases. There are various strategies for explaining away or dismissing these  
exceptions, and thus the entire enterprise of classifying the finite simple groups may be viewed as 
a vindication of the outlook which places Lie groups and algebraic groups first and foremost. 
Complacency in the face of these phenomena becomes a bit more difficult when one 
considers not only the sporadic groups but also the ``exotic fusion systems" which have been 
proliferating over the last twenty years or so - and for which there are as yet no general 
organizing principles. 

The aim of this series of papers is to begin at the other end of the thread, and 
to present a unified approach, from the view-point of a single prime $p$, to the 
three classes of groups already mentioned, to exotic fusion systems on finite $p$-groups, to 
the $p$-compact ``groups" of Dwyer and Wilkerson, and to other structures 
whose investigation it will be an aim of these papers to initiate. 
The main objects of study will not be groups at all, but rather ``partial groups" which may be thought of 
as ``locally grouped spaces" in somewhat the same way that, in algebraic geometry, schemes are  
locally ringed spaces. The ``affine spaces" in this analogy 
will be countable, locally finite groups $G$ having the property that a certain lattice $\Omega_S(G)$ 
of subgroups of a maximal $p$-subgroup $S$ of $G$ is ``finite-dimensional". 

For example, let $F$ be the algebraic closure of a finite field and let $G$ be a group having a faithful, 
finite-dimensional representation over $F$ (or let $G$ be a homomorphic image of such a group). Then 
$G$ is countable and locally finite, and from this it follows that for any given prime $p$ there exists 
a maximal $p$-subgroup $S$ of $G$. Let  $\Omega_S(G)$ be the set, partially ordered by inclusion, whose 
members are those 
subgroups $X$ of $S$ such that $X$ is the intersection of some set of $G$-conjugates of $S$. We show 
in Appendix A to this Part I that there exists an upper bound to the lengths of monotone chains in 
$\Omega_S(G)$; and it is in this  sense that we say that $\Omega_S(G)$ is finite-dimensional. This type 
of consideration will form the basis for all that will be done here, once the basic definitions and 
the basic properties concerning partial groups, objective partial groups, and localities, have been 
laid down. 
\vskip .1in 

What follows is a brief synopsis of what will be covered in Part I. 

\vskip .1in 
Let $G$ be a group and let $\bold W(G)$ be the free monoid on $G$. Thus, $\bold W(G)$ is the set of all 
words in the alphabet $G$, with the binary operation given by concatenation of words. The product 
$G\times G\to G$ extends, by generalized associativity, to a ``product" $\Pi:\bold W(G)\to G$, whereby a 
word $w=(g_1,\cdots,g_n)\in\bold W(G)$ is mapped to 
$g_1\cdots g_n$. The inversion map on $G$ induces an ``inversion" on $\bold W(G)$, sending  
$w$ to $(g_n\i,\cdots,g_1\i)$. In fact, one may easily replace the standard definition of ``group"
by a definition given in terms of $\Pi$ and the inversion on $\bold W(G)$. One obtains the notion of 
{\it partial group} by restricting the domain of $\Pi$ to a subset $\bold D$ of $\bold W(G)$, where 
$\bold D$, the product, and the inversion, are required to satisfy conditions (see definition 1.1) that 
preserve the outlines of the 
strictly group-theoretic setup. When one looks at things in this way, a group is simply a partial group 
$G$ having the property that $\bold D=\bold W(G)$. 

The notions of partial subgroup and homomorphism of partial groups will immediately suggest themselves,  
and a partial subgroup of a partial group $\ca L$ may in fact be a group. We say that the partial group 
$\ca L$ is ``objective" (see definition 2.1) provided that the domain $\bold D$ of the product is 
determined in a certain way by a collection $\D$ of subgroups of $\ca L$ (the set of ``objects"), and 
provided that $\D$ has a certain closure property. If moreover there exists an object $S\in\D$ such that 
$\D$ is a collection of subgroups of $S$ then $(\ca L,\D,S)$ is a {\it pre-locality}. 

Let $(\ca L,\D,S)$ be a pre-locality, let $w=(g_1,\cdots,g_n)$ be a non-empty word in $\bold W(\ca L)$, 
and let $x\in S$. Suppose that $x$ has the property that $((g_1)\i,x,g_1)$ is in the domain $\bold D(\ca L)$ 
of the product, and that the element $x_1=\Pi((g_1)\i,x,g_1)$ lies in $S$. Now suppose that $x_1$ has the 
property that $((g_2)\i,x_1,g_2)$ is in $\bold D(\ca L)$, and that the element $x_2=\Pi((g_2)\i,x_1,g_2)$ 
lies in $S$. If this procedure can be continued all the way to $((g_n)\i,x_{n-1},g_n)\in\bold D(\ca L)$ 
with $\Pi((g_n)\i,x_{n-1},g_n)\in S$, then we say that ``$x$ is conjugated sequentially into $S$ by 
the entries of $w$". The set of all such $x\in S$ is denoted $S_w$, and it is a consequence of the 
axioms (see 2.10) that $S_w$ is a subgroup of $S$. Indeed, if $\ca L$ were a bona fide group then 
$S_w$ would simply be the intersection of a collection of $\ca L$-conjugates of $S$. One has the set 
$$
 \Omega=\Omega_S(\ca L)\overset\text{def}\to = \{S_w\mid w\in\bold W(\ca L)\} 
$$ 
of all such $S_w$; with $S_w$ defined to be $S$ if $w$ is the empty word. We regard $\Omega$ as a poset 
via inclusion, and we say that $\Omega$ is {\it finite-dimensional} if there is an upper bound to the 
lengths of monotone chains in $\Omega$. 

The pre-locality $(\ca L,\D,S)$ is a {\it locality} if the following three conditions hold.  
\roster 

\item "{(L1)}" All subgroups of $\ca L$ are locally finite and countable.  

\item "{(L2)}" $S$ is a $p$-group, and is maximal (with respect to inclusion) in the set of all $p$-subgroups 
of $\ca L$. 

\item "{(L3)}" $\Omega_S(\ca L)$ is finite-dimensional, 

\endroster  

The basic properties of partial groups, objective partial groups, and localities, are derived in 
sections 1 through 3. We then begin in section 4 to consider partial normal subgroups of localities in  
detail.  One of the key results in section 4 is Stellmacher's splitting lemma (4.12), which leads to 
the partition of $\ca L$ into a collection of ``maximal cosets" of $\ca N$, and to a partial 
group structure on the set $\ca L/\ca N$ of 
maximal cosets. In section 5 it is shown that $\ca L/\ca N$ is in fact a locality, and we obtain 
versions of the first N\" other isomorphism theorem and of its familiar consequences. 
Section 6 concerns products of partial normal subgroups, and the main results here (Theorems 6.7 and 
6.8) are based entirely on the treatment by Ellen Henke [He] of the case of finite localities.  
\vskip .1in

Finite localities were introduced by the first named author in [Ch1], in order to give a positive solution 
to a basic existence/uniqueness question concerning fusion systems over finite $p$-groups. 
The solution that was given in [Ch1] was thus tied to a narrow goal, and did not allow for 
a complete development of ideas. The aim here is to attempt such a development, to thereby obtain  
the prerequisites for a version of finite group theory itself from a strictly ``$p$-local" point of view, 
and to extend that point of view beyond the finite context.  Concommitantly, our aim is to enrich the theory 
of fusion systems and to broaden its scope, so as to provide a  
 deeper connection with the homotopy theory in which - through Bob Oliver's proof 
[O1 and O2] of the Martino-Priddy conjecture - the entire enterprise has its roots. 

A locality is something which need not be a group, but which can have plenty of subgroups. All 
of the groups to be considered here will be subgroups of localities, and a basic assumption 
throughout is that any subgroup of a locality should be countable and locally finite. Thus, when 
we speak of a $p$-group $P$, we mean a countable, locally finite 
group all of whose elements have order a power of $p$. Such a group $P$ need not have the property 
(important when working with finite $p$-groups) that the normalizer in $P$ of a proper subgroup 
$Q$ of $P$ necessarily contains $Q$ properly. It will become necessary to impose this ``normalizer-increasing 
property" in Parts II and III, but in this Part I we can proceed without it. A weak version of this property 
(see 3.2 and 3.3 below) is already implied by the hypothesis of finite-dimensionality. 

The division into Parts closely parallels the extent to which fusion systems on the one hand, and primes 
other than $p$ on the other, are drawn into the developing picture. This Part I can be characterized by its 
having no direct involvement with fusion systems at all, and by there being no attention paid to $p'$-groups. 
The main results here, after the foundations for partial groups and localities have been laid down in the 
first three sections, concern partial normal subgroups of localities, the 
corresponding quotient localities, and finally the result (based in part on a result 
of Henke [He]) that products of partial normal normal subgroups are again partial normal subgroups. 

Readers who wish to consult a version of these papers in which all the localities under consideration are 
finite are referred to [Ch2] and [Ch3]. 

\vskip .1in 
Concerning notation: We adopt the group-theorist's use of ``$X^g$" to denote conjugation of a subset or 
element $X$ of a group $G$ by an element $g\in G$. Consistent with this practice, mappings will 
(almost always) be written to the right of their arguments.

\vskip .2in 
\noindent 
{\bf Section 1: Partial groups} 
\vskip .1in 

The reader is asked to forget what a group is, and to trust that what was forgotten will soon be 
recovered. 

\vskip .1in 
For any set $X$ write $\bold W(X)$ for the free monoid on $X$. Thus, an element 
of $\bold W(X)$ is a finite sequence of (or {\it word in}) the elements of $X$. The multiplication in 
$\bold W(X)$ consists of concatenation of words, to be denoted $u\circ v$. The {\it length} of the word 
$(x_1,\cdots,x_n)$ is $n$. The {\it empty word} is the word $(\emptyset)$ of length 0. We make no 
distinction between $X$ and the set of words of length $1$. 

\definition {Definition 1.1} Let $\ca L$ be a non-empty set, let $\bold W=\bold W(\ca L)$ be the free monoid 
on $\ca L$, and let $\bold D$ be a subset of $\bold W$ such that: 
\roster 

\item "{(1)}" $\ca L\sub\bold D$ (i.e. $\bold D$ contains all words of length 1), and 
$$ 
u\circ v\in\bold D\implies u,v\in\bold D. 
$$
\endroster 
Notice that since $\ca L$ is non-empty, (1) implies that also the empty word is in $\bold D$. 

A mapping $\Pi:\bold D\to\ca L$ is a {\it product} if: 
\vskip .1in
\roster 

\item "{(2)}" $\Pi$ restricts to the identity map on $\ca L$, and 
\vskip .1in 

\item "{(3)}" $u\circ v\circ w\in\bold D\implies u\circ(\Pi(v))\circ w\in\bold D$, and 
$\Pi(u\circ v\circ w)=\Pi(u\circ(\Pi(v))\circ w)$. 
\vskip .1in 
\endroster 

An {\it inversion} on $\ca L$ consists of an involutory bijection $x\maps x\i$ on $\ca L$, 
together with the mapping $w\maps w\i$ on $\bold W$ given by 
$$ 
(x_1,\cdots,x_n)\maps(x_n\i,\cdots x_1\i). 
$$ 
We say that $\ca L$, with the product $\Pi:\bold D\to\ca L$ and inversion $(-)\i$, is a 
{\it partial group} if: 
\roster 

\item "{(4)}" $w\in\bold D\implies w\i\circ w\in\bold D$ and $\Pi(w\i\circ w)=\1$, 

\endroster 
where $\1$ denotes the image of the empty word under $\Pi$. (Notice that (1) and (4) yield 
$w\i\in\bold D$ if $w\in\bold D$. As $(w\i)\i=w$, condition (4) is then symmetric.) 
\enddefinition 

\definition {Example 1.2} Let $\ca L$ be the $3$-element set $\{\1,a,b\}$ and let $\bold D$ be the subset of 
$\bold W(\ca L)$ consisting of all words $w$ such that the word obtained from $w$ by deleting all entries 
equal to $\1$ is an alternating string of $a$'s and $b$'s (of odd or even length, beginning with $a$ or 
beginning with $b$). Define $\Pi:\bold D\to\ca L$ by the formula: $\Pi(w)=\1$ if the number 
of $a$-entries in $w$ is equal to the number of $b$'s; $\Pi(w)=a$ if the number of $a$'s  
exceeds the number of $b$'s (necessarily by 1); and $\Pi(w)=b$ if the number of $b$'s 
exceeds the number of $a$'s. Define inversion on $\ca L$ by $\1\i=\1$, $a\i=b$, and $b\i=a$. 
It is then easy to check that $\ca L$, with these structures, is a partial group. In fact, $\ca L$ is 
the ``free partial group on one generator", as will be made clear in 1.12. 
\enddefinition 
 
It will be convenient to make the definition: a {\it group} is a partial group $\ca L$ in which 
$\bold W(\ca L)=\bold D$. In order to distinguish between this definition and the usual one, we shall 
use the expression {\it binary group} for a set $G$ with an associative binary operation, identity 
element, and inverses in the usual sense. The following lemma shows that the distinction is subtle. 

\proclaim {Lemma 1.3} 
\roster  

\item "{(a)}" Let $G$ be a binary group and let $\Pi:\bold W(G)\to G$ be the ``multivariable 
 product" on $G$ given by $(g_1,\cdots,g_n)\maps g_1\cdots g_n$. Then $G$, together with $\Pi$ and 
the inversion in $G$, is a partial group, with $\bold D=\bold W(G)$. 

\item "{(b)}" Let $\ca L$ be a group; i.e. a partial group for which $\bold W=\bold D$. Then $\ca L$ is a 
binary group with respect to the operation given by restricting $\Pi$ to words of length 2, and 
with  respect to the inversion in $\ca L$.  Moreover, $\Pi$ is then the multivariable product on $\ca L$ 
defined as in (a). 

\endroster 
\endproclaim 

\demo {Proof} Point (a) is given by generalized associativity in the binary group $G$. Point 
(b) is a straightforward exercise, and is left to the reader. 
\qed 
\enddemo 

We now list some elementary consequences of definition 1.1. 

\proclaim {Lemma 1.4} Let $\ca L$ $($with $\bold D$, $\Pi$, and the inversion$)$ be a partial group. 
\roster 

\item "{(a)}" $\Pi$ is {\it $\bold D$-multiplicative}. That is, if $u\circ v$ is in $\bold D$ then 
the word $(\Pi(u),\Pi(v))$ of length 2 is in $\bold D$, and 
$$
\Pi(u\circ v)=\Pi(u)\Pi(v), 
$$
where ``$\Pi(u)\Pi(v)$" is an abbreviation for $\Pi((\Pi(u),\Pi(v))$.

\item "{(b)}" $\Pi$ is {\it $\bold D$-associative}. That is: 
$$
u\circ v\circ w\in\bold D\implies\Pi(u\circ v)\Pi(w)=\Pi(u)\Pi(v\circ w). 
$$ 

\item "{(c)}" If $u\circ v\in\bold D$ then $u\circ(\1)\circ v\in\bold D$ 
and $\Pi(u\circ(\1)\circ v)=\Pi(u\circ v)$. 

\item "{(d)}" If $u\circ v\in\bold D$ then both $u\i\circ u\circ v$ and 
$u\circ v\circ v\i$ are in $\bold D$, 
$\Pi(u\i\circ u\circ v)=\Pi(v)$, and $\Pi(u\circ v\circ v\i)=\Pi(u)$. 

\item "{(e)}" The cancellation rule: If $u\circ v$, $u\circ w\in\bold D$, and 
$\Pi(u\circ v)=\Pi(u\circ w)$, then $\Pi(v)=\Pi(w)$ $($and similarly for 
right cancellation$)$. 

\item "{(f)}" If $u\in\bold D$ then $u\i\in\bold D$, and $\Pi(u\i)=\Pi(u)\i$. 
In particular, $\1\i=\1$. 

\item "{(g)}" The uncancellation rule: Suppose that both 
$u\circ v$ and $u\circ w$ are in $\bold D$ and that $\Pi(v)=\Pi(w)$. Then 
$\Pi(u\circ v)=\Pi(u\circ w)$. $($Similarly for right uncancellation.$)$ 

\endroster 
\endproclaim 

\demo {Proof} Let $u\circ v\in\bold D$. Then 1.1(3) applies to 
$(\nset)\circ u\circ v$ and yields 
$(\Pi(u))\circ v\in\bold D$ with $\Pi(u\circ v)=\Pi((\Pi(u))\circ v)$. Now 
apply 1.1(3) to $(\Pi(u))\circ v\circ(\nset)$, to obtain (a). 

Let $u\circ v\circ w\in\bold D$. Then $u\circ v$ and 
$w$ are in $\bold D$ by 1.1(1), and $\bold D$-multiplicativity yields 
$\Pi(u\circ v\circ w)=\Pi(u\circ v)\Pi(w)$. Similarly, 
$\Pi(u\circ v\circ w)=\Pi(u)\Pi(v\circ w)$, and (b) holds. 

Since $\1=\Pi(\nset)$, point (c) is immediate from 1.1(3). 

Let $u\circ v\in\bold D$. Then 
$v\i\circ u\i\circ u\circ v\in\bold D$ by 1.1(4), and then 
$u\i\circ u\circ v\in\bold D$ by 1.1(1). Multiplicativity then yields 
$$ 
\Pi(u\i\circ u\circ v)=\Pi(u\i\circ u)\Pi(v)=\1\Pi(v)= 
\Pi(\emptyset)\Pi(v)=\Pi(\emptyset\circ v)=\Pi(v). 
$$
As $(w\i)\i=w$ for any $w\in\bold W$, one obtains $w\circ w\i\in\bold D$ 
for any $w\in\bold D$, and $\Pi(w\circ w\i)=\1$. From this one easily completes 
the proof of (d). 

Now let $u\circ v$ and $u\circ w$ be in $\bold D$, with 
$\Pi(u\circ v)=\Pi(u\circ w)$. Then (d) (together with multiplicativity and 
associativity, which will not be explicitly mentioned hereafter) yield 
$$ 
\Pi(v)=\Pi(u\i\circ u\circ v)=\Pi(u\i)\Pi(u)\Pi(v)=\Pi(u\i)\Pi(u)\Pi(w)= 
\Pi(u\i\circ u\circ w)=\Pi(w), 
$$ 
and (e) holds. 

Let $u\in\bold D$. Then $u\circ u\i\in\bold D$, and then $\Pi(u)\Pi(u\i)=\1$. 
But also $(\Pi(u),\Pi(u)\i)\in\bold D$, and $\Pi(u)\Pi(u)\i=\1$. Now (f) 
follows by 1.1(2) and cancellation. 

Let $u,v,w$ be as in (g). Then $u\i\circ u\circ v$ and $u\i\circ u\circ w$ 
are in $\bold D$ by (d). By two applications of (d), 
$\Pi(u\i\circ u\circ v)=\Pi(v)=\Pi(w)=\Pi(u\i\circ u\circ w)$, so 
$\Pi(u\circ v)=\Pi(u\circ w)$ by (e), and (g) holds. 
\qed 
\enddemo

It will often be convenient to eliminate the symbol ``$\Pi$" and to speak of ``the product $g_1\cdots g_n$" 
instead of $\Pi(g_1,\cdots,g_n)$. More generally, if $\{X_i\}_{1\leq i\leq n}$ is a collection of subsets of 
$\ca L$ then the ``product set $X_1\cdots X_n$" is by definition the image under $\Pi$ of the set of words 
$(g_1,\cdots,g_n)\in\bold D$ such that $g_i\in X_i$ for all $i$. If $X_i=\{g_i\}$ is a singleton then we may 
write $g_i$ in place of $X_i$ in such a product. Thus, for example, the product 
$g\i Xg$ stands for the set of all $\Pi(g\i,x,g)$ with $(g\i,x,g)\in\bold D$, and with $x\in X$.

A Word of Urgent Warning: In writing products in the above way one may be drawn into imagining that 
associativity holds in a stronger sense than that which is given by 1.4(b). This is an error that is to 
be avoided. For example one should not suppose, if $(f,g,h)\in\bold W$, and both $(f,g)$ and $(fg,h)$ are in 
$\bold D$, that $(f,g,h)$ is in $\bold D$. That is, it may be that ``the 
product $fgh$" is undefined, even though the product $(fg)h$ is defined. Of 
course, one is tempted to simply extend the domain $\bold D$ to include such 
triples $(f,g,h)$, and to ``define" the product $fgh$ to be $(fg)h$. The 
trouble is that it may also be the case that $gh$ and $f(gh)$ are defined, but that $(fg)h\neq f(gh)$.

\vskip .1in 
For $\ca L$ a partial group and $g\in\ca L$, write $\bold D(g)$ for the 
set of all $x\in\ca L$ such that the product $g\i xg$ is defined. There is then a mapping 
$$
c_g:\bold D(g)\to\ca L 
$$
given by $x\maps g\i xg$ (and called {\it conjugation by $g$}). Our 
preference is for right-hand notation for mappings, so we write 
$$ 
x\maps(x)c_g\quad\text{or}\quad x\maps x^g 
$$ 
for conjugation by $g$.

\vskip .1in 
The following result provides an illustration of the preceding notational conventions, 
and introduces a theme which will be developed further as we pass from partial groups 
to objective partial groups, localities, and (in Part III) regular localities.

\proclaim {Lemma 1.5} Let $\ca L$ be a partial group, and let $f,g\in\ca L$. 
\roster 

\item "{(a)}" Suppose that the products $fg$ and $gf$ are defined and that 
$fg=gf$. Suppose further that $f\in\bold D(g)$. Then $f^g=f$. 

\item "{(b)}" Suppose that $f\in\bold D(g)$, $g\in\bold D(f)$, and $f^g=f$. Then 
$fg=gf$ and $g^f=g$. 

\endroster 
\endproclaim 

\demo {Proof} (a): We are given $(f,g)\in\bold D$, so $(f\i,f,g)\in\bold D$ and $\Pi(f\i,f,g)=g$, 
by 1.4(d) and $\bold D$-associativity. We are given also $f\in\bold D(g)$ and $fg=gf$, so  
$$ 
f^g=\Pi(g\i,f,g)=\Pi((g\i,fg)=\Pi(g\i,gf)=\Pi(g\i,g,f)=f. 
$$
\vskip .1in 
\noindent 
(b): As $(g\i,f,g)\in\bold D$ we obtain $(f,g)\in\bold D$ from 1.1(1). As $(f\i,g,f)\in\bold D$, 
we get $(g,g\i,f,g)\in\bold D$ by 1.4(d). Then $\bold D$-associativity yields 
$fg=\Pi(g,g\i,f,g)=gf^g$. As $f^g=f$ by hypothesis, we obtain $fg=gf$. Finally, since 
$(f\i,f,g)$ and $(f\i,g,f)$ are in $\bold D$ the uncancellation rule yields 
$f\i fg=f\i gf$, and so $g^f=g$. 
\qed 
\enddemo

\definition {Notation} From now on, in any given partial group $\ca L$, usage of the symbol ``$x^g$" shall 
be taken to imply $x\in\bold D(g)$. More generally, for $X$ a subset of $\ca L$ and $g\in\ca L$, usage of 
``$X^g$" shall be taken to mean that $X\sub\bold D(g)$; whereupon $X^g$ is by definition the set of all 
$x^g$ with $x\in X$. 
\enddefinition

At this early point, and in the context of arbitrary partial groups, one can say very little about the 
maps $c_g$. The cancellation rule 1.4(e) implies 
that each $c_g$ is injective, but beyond that, the following lemma may be the best that can be obtained.

\proclaim {Lemma 1.6} Let $\ca L$ be a partial group and let $g\in\ca L$. Then the following hold. 
\roster

\item "{(a)}" $\1\in\bold D(g)$ and ${\1}^g=\1$. 

\item "{(b)}" $\bold D(g)$ is closed under inversion, and $(x\i)^g=(x^g)\i$ for all $x\in\bold D(g)$. 

\item "{(c)}" $c_g$ is a bijection $\bold D(g)\to\bold D(g\i)$, and $c_{g\i}=(c_g)\i$. 

\item "{(d)}" $\ca L=\bold D(\1)$, and $x^{\1}=x$ for each $x\in\ca L$. 

\endroster 
\endproclaim 

\demo {Proof} By 1.1(4), $g\circ\nset\circ g\i=g\circ g\i\in\bold D$, so $\1\in\bold D(g)$ and then 
$\1^g=\1$ by 1.4(c). Thus (a) holds. Now let $x\in\bold D(g)$ and set $w=(g\i,x,g)$. Then $w\in\bold D$, 
and $w\i=(g\i,x\i,g)$ by definition in 1.1. Then 1.1(4) yields $w\i\circ w\in\bold D$, and so 
$w\i\in\bold D$ by 1.1(1). This shows that $\bold D(g)$ is closed under inversion. Also, 1.1(4) yields 
$\1=\Pi(w\i\circ w)=(x\i)^g x^g$, and then $(x\i)^g=(x^g)^{\i}$ by 1.4(f). This completes the proof of (b). 

As $w\in\bold D$, 1.4(d) implies that $g\circ w$ and then $g\circ w\circ g\i$ are in $\bold D$. Now 1.1(3) 
and two applications of 1.4(d) yield 
$$ 
gx^gg\i=\Pi(g,g\i,x,g,g\i)=\Pi((g,g\i,x)\circ g\circ g\i)=\Pi(g,g\i,x)=x. 
$$ 
Thus $x^g\in\bold D(g\i)$ with $(x^g)^{g\i}=x$, and thus (c) holds. 

Finally, $\1=\1\i$ by 1.4(f), and $\nset\circ x\circ\nset=x\in\bold D$ for any $x\in\ca L$, proving (d). 
\qed 
\enddemo 

\definition {Definition 1.7} Let $\ca L$ be a partial group and let $\ca H$ be a non-empty subset of 
$\ca L$. Then $\ca H$ is a {\it partial subgroup} of $\ca L$ (denoted $\ca H\leq\ca L$) if $\ca H$ is 
closed under inversion ($g\in\ca H$ implies $g\i\in\ca H$) and with respect to products. The latter condition 
means that $\Pi(w)\in\ca H$ whenever $w\in\bold W(\ca H)\cap\bold D$. A partial subgroup 
$\ca N$ of $\ca L$ is {\it normal} in $\ca L$ (denoted $\ca N\norm\ca L$) if $x^g\in\ca N$ for all 
$x\in\ca N$ and all $g\in\ca L$ for which $x\in\bold D(g)$. 
We say that $\ca H$ is a {\it subgroup} of $\ca L$ if $\ca H\leq\ca L$ and $\bold W(\ca H)\sub\bold D$. 

An equivalent way to state the condition for normality is to say that the 
partial subgroup $\ca N$ of $\ca L$ is normal in $\ca L$ if $g\i\ca Ng\sub\ca N$ for all $g\in\ca L$. 
(This formulation relies on a notational convention introduced above for interpreting 
product sets $XYZ$.) 
\enddefinition

We leave it to the reader to check that if $\ca H\leq\ca L$ then $\ca H$ is a partial group, with 
$\bold D(\ca H)=\bold W(\ca H)\cap\bold D(\ca L)$.

\proclaim {Lemma 1.8} Let $\ca H$ and $\ca K$ be partial subgroups of a partial group $\ca L$, and let 
$\{\ca H_i\}_{i\in I}$ be a set of partial subgroups of $\ca L$. 
\roster 

\item "{(a)}" Each partial subgroup of $\ca H$ is a partial subgroup of $\ca L$. 

\item "{(b)}" Each partial subgroup of $\ca L$ which is contained in $\ca H$ is a partial subgroup of $\ca H$. 

\item "{(c)}" If $\ca H$ is a subgroup of $\ca L$ then $\ca H\cap\ca K$ is a subgroup of both $\ca H$ and 
$\ca K$.  

\item "{(d)}" Suppose $\ca K\norm\ca L$. Then $\ca H\cap\ca K\norm\ca H$. 
Moreover, $\ca H\cap\ca K$ is a normal subgroup of $\ca H$ if $\ca H$ is a subgroup of $\ca L$. 

\item "{(e)}" $\bigcap\{\ca H_i\mid i\in I\}$ is a partial subgroup of $\ca L$, and is 
a partial normal subgroup of $\ca L$ if $\ca H_i\norm\ca L$ for all $i$.  

\endroster  
\endproclaim 

\demo {Proof} One observes that in all of the points (a) through (e) the requisite closure with respect to 
inversion obtains. Thus, we need only be concerned with products. 
\vskip .1in 
\noindent 
(a) Let $\ca E\leq\ca H$ be a partial subgroup of $\ca H$. Then 
$$ 
\bold D(\ca E)=\bold W(E)\cap\bold D(\ca H)=\bold W(E)\cap(\bold W(\ca H)\cap\bold D(\ca L)=
\bold W(E)\cap\bold D(\ca L), 
$$ 
and (a) follows. 
\vskip .1in 
\noindent 
(b) Suppose $\ca K\sub\ca H$ and let $w\in\bold W(\ca K)\cap\bold D(\ca H)$. As 
$\bold D(\ca H)\leq\bold D(\ca L)$, and since $\ca K\leq\ca L$ by hypothesis, we obtain $\Pi(w)\in\ca K$. 
\vskip .1in 
\noindent 
(c) Assuming now that $\ca H$ is a subgroup of $\ca L$, we have $\bold W(\ca H)\sub\bold D(\ca L)$, and 
then $\bold D(\ca H\cap\ca K)\sub\bold D(\ca H)\cap\bold D(\ca K)$, so that $\ca H\cap\ca K$ is a subgroup of 
both $\ca H$ and $\ca K$. 
\vskip .1in 
\noindent 
(d) Let $\ca K\norm\ca L$ and let $x\in\ca H\cap\ca K$ and $h\in\ca H$ with $(h\i,x,h)\in\bold D(\ca H)$. 
Then $(h\i,x,h)\in\bold D(\ca L)$, and $x^h\in\ca K$. As $\ca H\leq\ca L$ we have also $x^h\in\ca H$, and 
so $\ca H\cap\ca K\norm\ca H$. Now suppose further that $\ca H$ is a subgroup of $\ca L$. That is, assume 
that $\bold W(\ca H)\sub\bold D(\ca L)$. Then $\bold W(\ca H\cap\ca K)\sub\bold D(\ca L)$, hence 
$\ca H\cap\ca K$ is a subgroup of $\ca H$, and evidently a normal subgroup.   
\vskip .1in 
\noindent 
(e) Set $\ca X=\bigcap\{\ca H_i\}_{i\in I}$. Then $\Pi(w)\in\ca X$ for all 
$w\in\bold W(\ca X)\cap\bold D(\ca L)$, 
and so $\ca X\leq\ca L$. The last part of (e) may be left to the reader. 
\qed 
\enddemo

For any subset $X$ of a partial group $\ca L$ define the partial subgroup $\<X\>$ of $\ca L$ 
{\it generated by $X$} to be the 
intersection of the set of all partial subgroups of $\ca L$ containing $X$. 
Then $\<X\>$ is itself a partial subgroup of $\ca L$ by 1.8(e).

\proclaim {Lemma 1.9} Let $X$ be a subset of $\ca L$ such that $X$ is closed under inversion. Set 
$X_0=X$ and recursively define $X_n$ for $n>0$ by 
$$ 
X_n=\{\Pi(w)\mid w\in\bold W(X_{n-1})\cap\bold D\}. 
$$ 
Then $\<X\>=\bigcup\{X_n\}_{n\geq 0}$. 
\endproclaim 

\demo {Proof} Let $Y$ be the union of the sets $X_i$. Each $X_i$ is closed under inversion by 1.4(f), and 
$Y\neq\nset$ since $\1=\Pi(\nset)$. Since $Y$ is closed under products, by construction, we get $Y\leq\<X\>$, 
and then $Y=\<X\>$ by the definition of $\<X\>$. 
\qed 
\enddemo

\proclaim {Lemma 1.10 (Dedekind Lemma)} Let $\ca H$, $\ca K$, and $\ca A$ be partial subgroups of a 
partial group $\ca L$, and assume that $\ca H\ca K$ is a partial subgroup of $\ca L$.  
\roster 

\item "{(a)}" If $\ca K\leq\ca A$ then $\ca A=(\ca A\cap\ca H)\ca K$. 

\item "{(b)}" If $\ca H\leq\ca A$ then $\ca A=\ca H(\ca A\cap\ca K)$.

\endroster  
\endproclaim 

\demo {Proof} The proof is identical to the proof for binary groups, and is left to the reader. 
\qed 
\enddemo 

\definition {Definition 1.11} Let $\ca L$ and $\ca L'$ be partial groups, let $\b:\ca L\to\ca L'$ be a 
mapping, and let $\b^*:\bold W\to\bold W'$ be the induced mapping of free monoids. Then $\b$ is a 
{\it homomorphism $($of partial groups$)$} if: 
\roster 

\item "{(H1)}" $\bold D\b^*\sub\bold D'$, and 

\item "{(H2)}" $(\Pi(w))\b=\Pi'(w\b^*)$ for all $w\in\bold D$. 

\endroster 
The {\it kernel} of $\b$ is the set $Ker(\b)$ of all $g\in\ca L$ such that $g\b=\1'$. We say that $\b$ is an 
{\it isomorphism} if there exists a homomorphism $\b':\ca L'\to\ca L$ such that $\b\circ\b'$ and 
$\b'\circ\b$ are identity mappings. (Equivalently, $\b$ is an isomorphism if $\b$ is bijective and 
$\bold D\b=\bold D'$.) 
\enddefinition 

\definition {Example 1.12} Let $\ca L=\{\1,a,b\}$ be the partial group from example 1.2, let $\ca L'$ be 
any partial group, and let $x\in\ca L'$. Then the mapping $\b:\ca L\to\ca L'$ given by 
$$ 
\1\maps\1', \quad a\maps x, \quad  b\maps x\i. 
$$ 
is a homomorphism. In fact, $\b$ is the unique homomorphism $\ca L\to\ca L'$ which maps $a$ to $x$, 
by the following lemma. Thus, $\ca L$ is the (unique up to a unique invertible homomorphism) 
{\it free partial group} on one generator. Free partial groups in general can be obtained as ``free 
products" of copies of $\ca L$ (see 1.17 below). 
\enddefinition

\proclaim {Lemma 1.13} Let $\b:\ca L\to\ca L'$ be a homomorphism of partial 
groups. Then $\1\b=\1'$, and $(g\i)\b=(g\b)\i$ for all $g\in\ca L$. 
\endproclaim 

\demo {Proof} Since $\1 \1=\1$, (H1) and (H2) yield $\1\b=(\1\1)\b=(\1\b)(\1\b)$, and then $\1\b=\1'$ 
by left or right cancellation. Since $(g,g\i)\in\bold D$ for any $g\in\ca L$ by 1.4(d),
(H1) yields $(g\b,(g\i)\b)\in\bold D'$, and then
$\1\b=(gg\i)\b=(g\b)((g\i)\b)$ by (H2). As $\1\b=\1'=(g\b)(g\b)\i$,
left cancellation yields $(g\i)\b=(gb)\i$.
\qed
\enddemo

\proclaim {Lemma 1.14} Let $\b:\ca L\to\ca L'$ be a homomorphism of partial groups, and set $\ca N=Ker(\b)$. 
Then $\ca N$ is a partial normal subgroup of $\ca L$.  
\endproclaim 

\demo {Proof} By 1.13 $\ca N$ is closed under inversion. For $w$ in $\bold W(\ca N)\cap\bold D$ the map 
$\b^*:\bold W\to\bold W'$ sends $w$ to a word of the form $(\1',\cdots,\1')$. Then $\Pi'(w\b^*)=\1'$, 
and thus $\Pi(w)\in\ca N$ and $\ca N$ is a partial subgroup of $\ca L$. Now let 
$f\in\ca L$ and let $g\in\ca N\cap\bold D(f)$. Then 
$$
(f\i,g,f)\b^*=((f\b)\i,{\1}',f\b) \quad\text{(by 1.13)}, 
$$
so that 
$$
(g^f)\b=\Pi'((f\i,g,f)\b^*)=\Pi'(f\b)\i,\1',f\b)=\1'. 
$$ 
Thus $\ca N\norm\ca L$. 
\qed 
\enddemo

It will be shown later (cf. 4.6) that partial normal subgroups of ``localities" are always 
kernels of homomorphisms.

\proclaim {Lemma 1.15} Let $\b:\ca L\to\ca L'$ be a homomorphism of partial groups and let $M$ be a 
subgroup of $\ca L$. Then $M\b$ is a subgroup of $\ca L'$. 
\qed 
\endproclaim 

\demo {Proof} We are given $\bold W(M)\sub\bold D(\ca L)$, so $\b^*$ maps $\bold W(M)$ into 
$\bold D(\ca L')$. (Note, however, example 1.12.) 
\qed 
\enddemo

\proclaim {Lemma 1.16} Let $G$ and $G'$ be groups (and hence also binary groups in the sense of 1.3). 
A map $\a:G\to G'$ is a homomorphism of partial groups if and only if $\a$ is a homomorphism of 
binary groups. 
\endproclaim  

\demo {Proof} We leave to the reader the proof that if $\a$ is a homomorphism of partial groups then 
$\a$ is a homomorphism of binary groups. Now suppose that $\a$ is a homomorphism of binary groups. 
As $\bold W(G)=\bold D(G)$ (and similarly for $G'$, it is immediate that $\a^*$ maps $\bold D(G)$ into 
$\bold D(G)$. Assume that $\a$ is not a 
homomorphism of partial groups and let $w\in\bold D(G)$ be of minimal length subject to 
$\Pi'(w\a^*)\neq(\Pi(w))\a$. Then $n>1$ and we can write $w=u\circ v$ with both $u$ and $v$ non-empty. Then 
$$ 
\Pi'(w\a^*)=\Pi'(u\a^*\circ v\a^*)=\Pi'(u\a^*)\Pi'(v\a^*)=((\Pi(u))\a)((\Pi(v))\a)=(\Pi(u)\Pi(v))\a, 
$$ 
as $\a$ is a homomorphism of binary groups. Since $(\Pi(u)\Pi(v))\a=(\Pi(w))\a$, the proof is complete. 
\qed 
\enddemo

\vskip .2in 
\noindent 
{\bf Section 2: Objective partial groups and pre-localities} 
\vskip .1in

Recall the convention: if $X$ is a subset of the partial group $\ca L$, and $g\in\ca L$, then any statement 
involving the expression ``$X^g$" is to be understood as carrying the assumption that $X\sub\bold D(g)$. 
Thus, the statement ``$X^g=Y$" means: $(g\i,x,g)\in\bold D$ for all $x\in X$, and $Y$ is the set of 
products $g\i xg$ with $x\in X$.

\definition {Definition 2.1} Let $\ca L$ be a partial group. For any collection $\D$ of subgroups of 
$\ca L$ define $\bold D_\D$ to be the set of all $w=(g_1,\cdots,g_n)\in\bold W(\ca L)$ such that:
\roster

\item "{(*)}" there exists a sequence $(X_0,\cdots,X_n)$ of elements of $\D$ such that 
$(X_{i-1})^{g_i}=X_i$ for all $i$ ($1\leq i\leq n$).

\endroster 
Then $\ca L$ is {\it objective} if there exists a set $\D$ of subgroups of $\ca L$ 
such that the following two conditions hold.
\roster

\item "{(O1)}" $\bold D=\bold D_{\D}$.

\item "{(O2)}" Whenever $X$ and $Y$ are in $\D$ and $g\in\ca L$ such that $X^g$ is a 
subgroup of $Y$, then every subgroup of $Y$ containing $X^g$ is in $\D$. 

\endroster 
We say also that $\D$ is a set of {\it objects}) for $\ca L$, if (O1) and (O2) hold.  
\enddefinition

It will often be convenient to somewhat over-emphasize the role of $\D$ in the above definition by saying 
that ``$(\ca L,\D)$ is an objective partial group". What is meant by this is that $\ca L$ is an objective 
partial group and that $\D$ is a set (there will often be more than one) of objects for $\ca L$.

We emphasize that the condition (O2) requires more than that $X^g$ be a sub{\it set} of $Y$, in order to 
conclude that overgroups of $X$ in $Y$ are objects. This is a non-vacuous distinction, since the 
conjugation map $c_g:X\to X^g$ need not send $X$ to a subgroup (or even a partial subgroup) 
of $\ca L$, in a general partial group. 

\definition {Example 2.2} Let $G$ be a group, let $S$ be a subgroup of $G$, and let $\D$ be 
a collection of subgroups of $S$ such that $S\in\D$. Assume that $\D$ satisfies (O2). That is, assume that 
$Y\in\D$ for every 
subgroup $Y$ of $S$ such that $X^g\leq Y$ for some $X\in\D$ and some $g\in G$. Let $\ca L$ 
be the set of all $g\in G$ such that $S\cap S^g\in\D$, and let $\bold D$ be the subset 
$\bold D_\D$ of $\bold W(\ca L)$. Then $\ca L$ is a partial group (via the multivariable 
product in $G$ and the inversion in $G$), and $(\ca L,\D)$ is an objective partial group. Specifically: 
\roster 

\item "{(a)}" If $\D=\{S\}$ then $\ca L=N_G(S)$ (and so $\ca L$ is a group in this case). 

\item "{(b)}" Take $G=O_4^+(2)$. Thus, $G$ is a semidirect product $V\rtimes S$ where $V$ is elementary 
abelian of order 9 and $S$ is a dihedral group of order 8 acting faithfully on $V$. Let $\D$ be the set of 
all non-identity subgroups of $S$. One may check that $S\cap S^g\in\D$ for all $g\in G$, and hence 
$\ca L=G$ (as {\it sets}). But $\ca L$ is not a group, as $\bold D_\D\neq\bold W(G)$. 

\item "{(c)}" Take $G=GL_3(2)$, $S\in Syl_2(G)$, and let $M_1$ and $M_2$ be the two maximal subgroups of $G$ 
containing $S$. Set $P_i=O_2(M_i)$ and set $\D=\{S,P_1,P_2\}$. Then $\ca L=M_1\cup M_2$ (in fact the ``free 
amalgamated product" of $M_1$ and $M_2$ over $S$ in the category of partial groups). On 
the other hand, if $\D$ is taken to be the set of all non-identity subgroups of $S$ then $\ca L$ is 
somewhat more complicated. Its underlying set is $M_1M_2\cup M_2M_1$. 

\endroster 
\enddefinition

\vskip .1in 
In an objective partial group $(\ca L,\D)$ we say that the word $w=(g_1,\cdots,g_n)$ is 
{\it in $\bold D$ via $(X_0,\cdots,X_n)$} if the condition (*) in 2.1 applies specifically 
to $w$ and $(X_0,\cdots,X_n)$. We may also say, more simply, that {\it $w$ is in $\bold D$ 
via $X_0$}, since the sequence $(X_0,\cdots,X_n)$ is determined by $w$ and $X_0$. 

\vskip .1in 
For any partial group $\ca L$ and subgroups $X,Y$ of $\ca L$, set   
$$ 
N_{\ca L}(X,Y)=\{g\in\ca L\mid X\sub\bold D(g),\ X^g\sub Y\},   
$$ 
and set 
$$ 
N_{\ca L}(X)=\{g\in\ca L\mid X^g=X\}. 
$$ 

\proclaim {Lemma 2.3} Let $(\ca L,\D)$ be an objective partial group. 
\roster 

\item "{(a)}" $N_{\ca L}(X)$ is a subgroup of $\ca L$ for each $X\in\D$. 

\item "{(b)}" Let $g\in\ca L$ and let $X\in\D$ with $Y:=X^g\in\D$. Then $N_{\ca L}(X)\sub\bold D(g)$, and 
$$
c_g:N_{\ca L}(X)\to N_{\ca L}(Y) 
$$
is an isomorphism of groups. More generally:

\item "{(c)}" Let $w=(g_1,\cdots,g_n)\in\bold D$ via $(X_0,\cdots,X_n)$. 
Then 
$$
c_{g_1}\circ\cdots\circ c_{g_n}=c_{\Pi(w)} 
$$
as isomorphisms from $N_G(X_0)$ to $N_G(X_n)$. 

\endroster 
\endproclaim 

\demo {Proof}(a) Let $X\in\D$ and let $u\in\bold W(N_{\ca L}(X))$. Then $u\in\bold D$ 
via $X$, $\1\in N_{\ca L}(X)$ (1.6(d)), and $N_{\ca L}(X)\i=N_{\ca L}(X)$ (1.6(c)). 
\vskip .1in 
\noindent 
(b) Let $x,y\in N_{\ca L}(X)$ and set $v=(g\i,x,g,g\i,y,g)$. Then $v\in\bold D$ via $Y$, 
and then $\Pi(v)=(xy)^g=x^gy^g$ (using points (a) and (b) of 1.4). Thus, the 
conjugation map $c_g:N_{\ca L}(X)\to N_{\ca L}(Y)$ is a homomorphism of binary groups  
(see 1.3), and hence a homomorphism of partial groups (1.16). Since $c_{g\i}=c_g\i$ by 1.6(c), 
$c_g$ is an isomorphism of groups. 
\vskip .1in 
\noindent 
(c) Let $x\in N_{\ca L}(X_0)$, set $u_x=w\i\circ(x)\circ w$, and observe that $u_x\in\bold D$ 
via $X_n$. Then $\Pi(u_x)$ can be written as $(\cdots(x)^{g_1}\cdots)^{g_n}$, and this yields (c). 
\qed 
\enddemo

The next lemma provides two basic computational tools.

\proclaim {Lemma 2.4} Let $(\ca L,\D)$ be an objective partial group. 
\roster 

\item "{(a)}"  Let $(a,b,c)\in\bold D$ and set $d=abc$. Then $bc=a\i d$ and $ab=dc\i$ 
(and all of these products are defined). 

\item "{(b)}" Let $(f,g)\in\bold D$ and let $X\in\D$. Suppose that both 
$X^f$ and $X^{fg}$ are in $\D$. Then $X^{fg}=(X^f)^g$. 

\endroster 
\endproclaim 

\demo {Proof} Point (a) is a fact concerning partial groups in general, and is immediate from 1.4(c). 
Now consider the setup in (b). As $(f,g)\in\D$ we have also $(f\i,f,g)\in\D$, and $g=\Pi(f\i,f,g)=f\i(fg)$. 
Now observe that $(f\i,fg)\in\bold D$ via $P^f$, and apply 2.3(c) to obtain 
$P^{fg}=((P^f)^{f\i})^{fg}=(P^f)^g$.  
\qed 
\enddemo 

The following result is a version of lemma 1.5(b) for objective partial groups. 
The hypothesis is weaker than that of 1.5(b), and the conclusion is stronger.

\proclaim {Lemma 2.5} Let $(\ca L,\D)$ be an objective partial group and let $f,g\in\ca L$. 
Suppose that $f^g=f$. Then $fg=gf$ and $g^f=g$. 
\endproclaim 

\demo {Proof} Suppose that $(g\i,f,g)\in\bold D$ via $(P_0,P_1,P_2,P_3)$. One then has the following 
commutative square of conjugation maps, in which the arrows are labeled by elements that perform the 
conjugation.
$$ 
\CD 
P_2@>g>>P_3 \\ 
@AfAA  @AA{f^g}A \\ 
P_1@>>g>P_0 
\endCD 
$$ 
Now assume that $f^g=f$. Since any of the arrows in the diagram may be reversed using 1.6(c), one reads 
off that $(f\i,g,f)\in\bold D$ via $P_2$. Then $g^f=g$ and $fg=gf$ by 1.5(b). 
\qed 
\enddemo 

\definition {Definition 2.6} A {\it pre-locality} is an objective partial group $\ca L$ having the property 
that there exists a collection $\D$ of objects for $\ca L$ such that $\D$ is a set of subgroups of some 
$S\in\D$. 
\enddefinition 

When we wish  to emphasize the role of $S$ and of $\D$, we may say that $(\ca L,\D,S)$ is a 
pre-locality. 
\vskip .1in 

\proclaim {Lemma 2.7} Let $(\ca L,\D,S)$ be a pre-locality. Then $\bold D$ is an $N_{\ca L}(S)$-biset. 
That is, if $w\in\bold D$  and $g,h\in N_{\ca L}(S))$, then $(g)\circ w\circ(h)\in\bold D$. 
In particular, $N_{\ca L}(S)$ acts on $\ca L$ by conjugation. 
\endproclaim 

\demo {Proof} If $w\in\bold D$ via $X\in\D$, then $(g)\circ v\circ(h)\in\bold D$ via $X^{g\i}$. 
The action of $N_{\ca L}(S)$ on $\ca L$ is then given by restricting to the case where $g=h\i$ and 
where $w$ is of length $1$. 
\qed 
\enddemo 

The following result and its corollary are fundamental to the entire enterprise. The proof is 
due to Bernd Stellmacher. 

\proclaim {Proposition 2.8} Let $(\ca L,\D,S)$ be a pre-locality. For each $g\in\ca L$ define 
$S_g$ to be the set of all $x\in\bold D(g)\cap S$ such that $x^g\in S$. Then: 
\roster 

\item "{(a)}" $S_g\in\D$. In particular, $S_g$ is a subgroup of $S$. 

\item "{(b)}" The conjugation map $c_g:S_g\to(S_g)^g$ is an isomorphism of groups, and 
$S_{g\i}=(S_g)^g$. 

\item "{(c)}" $P^g$ is defined and is a subgroup of $S$ for every subgroup $P$ of $S_g$. In 
particular, $P^g\in\D$ for any $P\in\D$ with $P\leq S_g$. 

\endroster 
\endproclaim 

\demo {Proof} Fix $g\in\ca L$. Then the word $(g)$ of length 1 is in $\bold D$ by 1.1(2), and since  
$\bold D=\bold D_\D$ by (O1) there exists $X\in\D$ such that $Y:=X^g\in\D$. Let $a\in S_g$ 
and set $b=a^g$. Then $X^a$ and $X^b$ are subgroups of $S$ (as $a,b\in S$), so $X^a$ and 
$Y^b$ are in $\D$ by (O2). Then $(a\i,g,b)\in\bold D$ via $X^a$, so also $(g,b)\in\bold D$. Also  
$(a,g)\in\bold D$ via $X^{a\i}$. Since $g\i ag=b$ we get $ag=gb$ by cancellation, and hence 
$$ 
a\i gb=a\i(gb)=a\i(ag)=(a\i a)g=g 
$$ 
by $\bold D$-associativity. Since $a\i gb$ conjugates $X^a$ to $Y^b$, we 
draw the following conclusion. 
\roster 

\item "{(1)}" $X^a\leq S_g$ and $(X^a)^g\in\D$ for all $a\in S_g$. 

\endroster 
Now let $c,d\in S_g$. Then (1) shows that both $X^c$ and $X^{cd}$ are members  
of $\D$ which are conjugated to members of $\D$ by $g$. Setting 
$w=(g\i,c,g,g\i,d,g)$, we conclude (by following $X^g$ along the chain of 
conjugations given by $w$) that $w\in\bold D$ via $X^g$. Then $\bold D$-associativity  
yields 
$$ 
\Pi(w)=(cd)^g=c^g d^g.  \tag2
$$ 
Since $c^g$ and $d^g$ are in $S$, we conclude that $cd\in S_g$. Since $S_g$ is closed under inversion 
by 1.6(b), $S_g$ is a subgroup of $S$. As $X\leq S_g\leq S$, where $X$ and $S$ are in $\D$, (02) now 
yields $S_g\in\D$. Thus (a) holds. 

Since $c_{g\i}=(c_g)\i$ by 1.6(c), it follows that $S_{g\i}=(S_g)^g$. Points (b) and (c) are then 
immediate from (a) and 2.3(b). 
\qed 
\enddemo

It will be convenient to extend the notation introduced preceding 1.5. Thus, for $\ca L$ a partial group 
and non-empty word $w=(g_1,\cdots,g_n)\in\bold W(\ca L)$, write 
$\bold D(w)$ for the set of all $x\in\ca L$ such that $x:=x_0\in\bold D(g_1)$, 
$x_1:=(x_0)^{g_1}\in\bold D(g_2)$, and so on until $y:=x_n=(x_{n-1})^{g_n}$. We then write 
$$
c_w:\bold D(w)\to\ca L 
$$
for the mapping $x\maps y$, and call this mapping {\it conjugation by $w$}. 
Thus, $\bold D(w)$ is the largest subset of $\ca L$ on which the composition 
$$ 
c_{g_1}\circ\cdots\circ c_{g_n} 
$$ 
is defined. If $w$ is the empty word define $\bold D(w)$ to be $\ca L$, and define $c_w$ to be the identity 
map on $\ca L$. For any subset or element $X$ of $\bold D(w)$, write $X^w$ for the image of $X$ under $c_w$. 

For each $w\in\bold W(\ca L)$ define $S_w$ to be the set of all $x\in\bold D(w)$ such that, in the preceding 
description of $c_w$, we have $x_i\in S$ for all $i$ ($0\leq i\leq n)$. In other words: 
$S_w$ is the set of all $x\in S$ 
such that $x$ is conjugated consecutively into $S$ by the sequence $(c_{g_1},\cdots,c_{g_n})$ of 
conjugation maps. Notice that $(S_w)^w=S_{w\i}$. The following lemma is a straightforward consequence 
of these definitions. 

\proclaim {Lemma 2.9} Let $(\ca L,\D,S)$ be a pre-locality, and let $u,v\in\bold W(\ca L)$. Then: 
\roster 

\item "{(a)}" $S_{u\circ v}=((S_u)^u\cap S_v)^{u\i}$, and 

\item "{(b)}" $S_u\cap S_v=S_{u\circ u\i\circ v}$. 

\endroster 
\qed 
\endproclaim 

\proclaim {Corollary 2.10} Let $(\ca L,\D,S)$ be a pre-locality, and let $w\in\bold W(\ca L)$. Then 
$S_w$ is a subgroup of $S$ for all $w\in\bold W(\ca L)$, and $S_w\in\D$ if and only if $w\in\bold D$. 
\endproclaim 

\demo {Proof} If the length of $w$ is at most $1$ then $S_w$ is a subgroup of $S$ by 2.8(a). This provides 
the basis for an induction which is completed by 2.9(a). 
If $S_w\in\D$ then $w\in\bold D$ by the condition (O1) in definition 2.1. Conversely, if 
$w\in\bold D$ then (O1) says that $S_w$ contains a member of $\D$, and hence $S_w\in\D$ by (O2). 
\qed 
\enddemo 

\definition {Definition 2.11} Let $(\ca L,\D,S)$ be a pre-locality. For each subgroup $X$ of $S$ set 
$$ 
X^\star=\bigcap\{S_w\mid w\in\bold W(\ca L),\ X\leq S_w\},  
$$ 
and set 
$$ 
\Omega:=\Omega_S(\ca L)=\{X^\star\mid X\leq S\}. 
$$ 
\enddefinition

From now on $(\ca L,\D,S)$ will be a pre-locality, and we shall write $\Omega$ for $\Omega_S(\ca L)$.

\proclaim {Lemma 2.12} Let $(\ca L,\D,S)$ be a pre-locality, and let $X,Y$ be subgroups of $S$. 
\roster 

\item "{(a)}" $X^\star=(X^\star)^\star$. 

\item "{(b)}" If $X\leq Y$ then $X^\star\leq Y^\star$. 

\item "{(c)}" $\<X^\star,Y^\star\>\leq \<X,Y\>^\star$. 

\endroster 
\endproclaim 

\demo {Proof} Points (a) and (b) are immediate from definition 2.11. Point (c) is an application 
of (b), with $\<X,Y\>$ in the role of $Y$. 
\qed 
\enddemo 

\proclaim {Lemma 2.13} Let $(\ca L,\D,S)$ be a pre-locality, let $u\in\bold W(\ca L)$, and let 
$X$ be a subgroup of $S_u$. Then $(X^\star)^u=(X^u)^\star$. 
\endproclaim 

\demo {Proof} Set $\bold W=\bold W(\ca L)$. Then 
$$  
X^\star=\bigcap\{S_w\mid w\in\bold W,\ X\leq S_w\}=\bigcap\{(S_u\cap S_w)\mid w\in\bold W,\ X\leq S_w\}
\leq S_u. 
$$ 
Then also 
$$
(X^\star)^u=\bigcap\{(S_u\cap S_w)^u\mid w\in\bold W,\ X\leq S_w\}. 
$$ 
Let $w\in\bold W(\ca L)$ with $X\leq S_w$. Then 2.9(a) yields  
$$ 
(S_u\cap S_w)^u=((S_{u\i})^{u\i}\cap S_w)^u=S_{u\i\circ w}, 
$$ 
and so 
$$ 
\align
(X^\star)^u&=\bigcap\{S_{u\i\circ w}\mid w\in\bold W,\ X\leq S_w\} \\ 
&\leq \bigcap\{S_{w'}\mid w'\in\bold W,\ X^u\leq S_{w'}\}\leq (X^u)^\star. 
\endalign 
$$ 
We may apply this result with $X^u$ in the role of $X$ and with $u\i$ in place of $u$, obtaining 
$(X^u)^\star\leq S_{u\i}$ and 
$$ 
((X^u)^\star)^{u\i}\leq ((X^u)^{u\i})^\star=X^\star. 
$$ 
Conjugation by $u$ then yields $(X^u)^\star\leq(X^\star)^u$, and we thus have the two inclusions 
required for equality. 
\qed 
\enddemo 

Regard $\Omega$ as a poset via the partial order given by inclusion, and define $dim(\Omega)$ to be 
the supremum (if such a supremum exists) of the numbers $d$ for which there exists a chain 
$$ 
X_0<\cdots<X_d \tag*
$$ 
of proper inclusions in $\Omega$. Similarly, for any $X\in\Omega$ 
write $dim_\Omega(X)$ for the supremum (if it exists) of the numbers $d$ for which there exists a chain 
(*) with $X_d=X$. The pre-locality $(\ca L,\D,S)$ will be said to be {\it finite-dimensional} if 
$dim(\Omega)<\infty$.

\proclaim {Proposition 2.14} Let $(\ca L,\D,S)$ be a finite-dimensional pre-locality, and set  
$\Omega=\Omega_S(\ca L)$. Then $\Omega=\{S_w\mid w\in\bold W(\ca L)\}$, and $\Omega$ is closed with 
respect to arbitrary intersections. 
\endproclaim 

\demo {Proof} By the definition of $X^\star$ for $X$ a subgroup of $S$, it suffices to show that 
$\Omega$ is closed with respect to arbitrary intersections. Such is the case, by an application of  
2.9(b), and by finite-dimensionality. 
\qed 
\enddemo

We now extend the notion of dimension to all subgroups of $S$, by the formula: 
$$ 
dim(X)\overset{def}\to= dim(X^\star). 
$$

\proclaim {Corollary 2.15} Let $(\ca L,\D,S)$ be a finite-dimensional pre-locality, and let 
$X$ and $Y$ be subgroups of $S$. Then the following hold. 
\roster 

\item "{(a)}" $dim(X)=dim(X^w)$ for all $w\in\bold W(\ca L)$ with $X\leq S_w$. 

\item "{(b)}" If $X\leq Y\leq S$ then $dim(X)\leq dim(Y)$. 

\item "{(c)}" If $X\leq Y\leq S$ with $X=X^\star$, then $dim(X)=dim(Y)$ if and only if $X=Y$. 

\endroster 
\endproclaim 

\demo {Proof} Let $w$ be as in (a), and let $\s=(X_0\leq\cdots\leq X_d=X^\star)$ be a chain of strict 
inclusions in $\Omega$, of length $d=dim(X)$. As $X^\star\leq S_w$ by 2.13, conjugation by $w$ sends 
$\s$ to a chain of strict inclusions terminating in $(X^\star)^w$. As $(X^\star)^w=(X^w)^\star$ (by 
2.13), we obtain (a). Point (b) is immediate from 2.12(b). Now suppose that $X=X^\star\leq Y$. 
Then $X\leq Y^\star$, and $dim(X)=dim(Y)$ if and only if $X=Y^\star$. As $X\leq Y\leq Y^\star$, (c) follows. 
\qed 
\enddemo 

\proclaim {Lemma 2.16} Let $(\ca L,\D,S)$ and $(\ca H,\G,R)$ be pre-localities, and let $\b:\ca H\to\ca L$ 
be an injective homomorphism of partial groups such that $R\b\leq S$. For any subgroup $X$ of $R$, write 
$X^\bullet=\bigcap\{R_w\mid w\in\bold W(\ca H),\ P\leq R_w\}$.  
\roster 

\item "{(a)}" There is an injective homomorphism $\Omega_R(\ca H)\to\Omega_S(\ca L)$ of posets, given 
by $X\maps (X\b)^\star$. 

\item "{(b)}" Assume that $\ca L$ is finite-dimensional. Then $dim(\Omega_R(\ca H))\leq dim(\Omega_S(\ca L))$,  
and if equality holds then $(\1^\bull)\b\leq\1^\star$ and $(R\b)^\star=S$. 

\endroster 
\endproclaim 

\demo {Proof} In proving (a) and (b) we may identify $\ca H$ with the image of $\b$, and so 
we may assume that $\b$ is an inclusion map. 

By 2.12(b), the map $X\maps X^\star$ is a homomorphism $\Omega_R(\ca H)\to\Omega_S(\ca L)$ of posets. 
In order to complete the proof of (a), the key observation is that $X^\star\cap R\leq X^\bullet$ for 
any subgroup $X$ of $R$. Thus $X^\star\cap R=X$ if $X\in\Omega_R(\ca H)$. If also $Y\in\Omega_R(\ca H)$ 
and $X^\star=Y^\star$ then $X=R\cap X^\star=R\cap Y^\star=Y$, and so (a) holds. 

Assume now that $\ca L$ is finite-dimensional, and let $\1^\bullet<X_1<\cdots<X_d$ be a chain of proper 
inclusions in $\Omega_R(\ca H)$. Since $\1\leq\1^\bullet\leq X$ for all $X\in\Omega_R(\ca H)$, we obtain 
$$ 
\1^\star\leq(\1^\bullet)^\star<(X_1)^\star<\cdots<(X_d)^\star\leq S  
$$
which is a chain of length greater than $d$ unless $\1^\star=(\1^\bullet)^\star$ and $(X_d)^\star=S$. 
This proves (b). 
\qed 
\enddemo

\proclaim {Proposition 2.17} Let $(\ca L,\D,S)$ be a finite-dimensional pre-locality and let $H$ be a 
subgroup of $\ca L$. Then there exists $P\in\D\cap\Omega$ such that $H\leq N_{\ca L}(P)$. Moreover, 
$P$ may be chosen so that $P=S_u$ for some $u\in\bold W(H)$, and then $P$ is the unique largest 
subgroup of $S$ normalized by $H$. 
\endproclaim 

\demo {Proof} Among all $u\in\bold W(H)$, choose $u$ so that $dim(S_u)$ is as small as possible, and 
set $P=S_u$. As $H$ is a subgroup of $\ca L$ we have $u\in\bold D$. Thus $P\in\D\cap\Omega$, and 
$P$ contains every subgroup $Q$ of $S$ such that $H\leq N_{\ca L}(Q)$.
Let $v$ be 
an arbitrary element of $\bold W(H)$, and set $w=u\circ u\i\circ v$. Then $S_w=S_u\cap S_v$ by 2.11(b), and 
$S_w\leq S_u$. Then $dim(S_w)=dim(S_u)$ by the minimality of $dim(S_u)$, and so $S_w=P$ by 2.15(c). Thus: 
$$ 
P=\bigcap\{S_v\mid v\cap\bold W(H)\}. 
$$ 
Now let $h\in H$. Then $P=S_{(h,h\i)\circ u}$ and $P=S_{(h\i,h)\circ u}$, whereas evidently 
$P^h=S_{(h\i,h)\circ u}$. Thus $P=P^h$, and $H\leq N_{\ca L}(P)$. 
by 2.12. Replacing $P$ with $P^\star$ yields $H\leq N_{\ca L}(P)$. If also 
$Q\in\D\cap\Omega$ with $H\leq N_{\ca L}(Q)$ then $Q\leq S_u=P$. 
\qed 
\enddemo

The proof of the following result is a straightforward exercise with definition 1.7, and it is left to 
the reader.

\proclaim {Lemma 2.18} Let $(\ca L,\D,S)$ be a pre-locality and let $X\leq S$ be a subgroup of $S$. Set 
$$ 
N_{\ca L}(X)=\{g\in\ca L\mid X\leq S_g\ \ \text{and}\ \ X^g=X\}, 
$$  
and 
$$ 
C_{\ca L}(X)=\{g\in L\mid X\leq S_g\ \ \text{and}\ \ x^g=a\ \ \text{for all}\ \ x\in X\}.  
$$ 
Then $N_{\ca L}(X)$ and $C_{\ca L}(X)$ are partial subgroups of $\ca L$, and 
$C_{\ca L}(X)\norm N_{\ca L}(X)$. 
\qed 
\endproclaim

\vskip .2in 
\noindent 
{\bf Section 3: Localities} 
\vskip .1in

We may now introduce the main object of study. The reader should recall the definition of 
pre-locality from 2.6, and the definition of the poset $\Omega_S(\ca L)$ from 2.11. Recall also that 
a poset $\Omega$ is defined to be finite-dimensional if there exists an upper bound on the lengths of 
strictly monotone chains in $\Omega$.

\definition {Definition 3.1} A {\it discrete locality} (or, for short, a {\it locality}) is a 
finite-dimensional pre-locality $(\ca L,\D,S)$ satisfying the following two conditions. 
\roster 

\item "{(L1)}" $S$ is a $p$-group for some prime $p$, maximal among the $p$-subgroups of $\ca L$.  

\item "{(L2)}" Each subgroup of $\ca L$ is locally finite and countable. 

\endroster 
In more detail, the partial group $\ca L$ is a locality if there exists a subgroup $S$ of $\ca L$ and a 
set $\D$ of subgroups of $S$ such that $(\ca L,\D)$ is objective, $(\ca L,\D,S)$ is a 
finite-dimensional pre-locality, and (L1) and (L2) hold. 
\enddefinition

\definition {Remarks} (1) By a $p$-group we mean a torsion group all of whose elements are of order a power 
of $p$. The $p$-subgroups of a partial group form a poset via inclusion, and it is with respect to this 
poset that we say that a given $p$-group is maximal. 
\vskip .1in 
\noindent 
(2) Notice that by 2.17, (L2) is equivalent to 
the requirement that each of the groups $N_{\ca L}(P)$, for $P\in\D$, be locally finite and countable. 
Another equivalent formulation of (L2) is that each subgroup of $\ca L$ is a nested union of 
a countable collection of finite groups. Since 
we shall make use of this last formulation from time to time, it will be convenient to establish  
the following terminology: a {\it framing} of a group $G$ (by finite subgroups) is a collection 
$\{G_n\}$ of finite subgroups of $G$, indexed by the non-negative integers, such that $G_n\leq G_{n+1}$ 
for all $n$, and such that $G=\bigcup\{G_n\}$. We shall write also $G=lim\{G_n\}$ if $\{G_n\}$ is 
a framing of $G$. 
\vskip .1in 
\noindent 
(3) The results 2.12 through 2.15 from the preceding section, which concern formal properties of 
$\Omega_S(\ca L)$ and of the operation $V\maps V^\star$ on subgroups of $S$ will be employed so often, and 
are (we believe) sufficiently natural, that explicit reference to them may usually be omitted.  
\vskip .1in 
\noindent  
(4) It will often be convenient, within discussions involving $p$-groups, to adopt the notation $P<Q$ to 
indicate that $P$ is a proper subgroup of $Q$. 
\enddefinition  

\proclaim {Lemma 3.2} Let $(\ca L,\D,S)$ be a locality, let $V$ be a subgroup of $S$, and let $X$ be a 
$p$-subgroup of $\ca L$ such that $V<X$. Then $V<N_X(V^\star)$, and if $V=V^\star\cap X$ then $V<N_X(V)$. 
In particular, we have $V<N_X(V)$ if $V\in\Omega$. 
\endproclaim 

\demo {Proof} Let $\{X_n\}$ be a framing of $X$ by finite subgroups, and set $V_n=V\cap X_n$. For $n$ 
sufficiently large  
we have $(V_n)^\star=V^\star$ (by finite-dimensionality) and $V_n<X_n$. As $X_n$ is a finite 
$p$-group there exists $x\in N_{X_n}(V_n)$ with $x\notin V_n$. Then $x\notin V$, while $x\in N_X(V^\star)$ 
by 2.13. Thus $V< N_X(V^\star)$, and if $V=V^\star\cap X$ then $V<N_X(V)$. 
\qed 
\enddemo

\proclaim {Corollary 3.3} Let $V\in\Omega$, and let $w\in\bold W(\ca L)$ with $V\leq S_w$. Set 
$P=N_{S_w}(V)$, set $Q=N_S(V)$, and assume that $P<Q$. Then $dim(P)<dim(N_Q(P))$. 
\endproclaim 

\demo {Proof} We have $P^\star\leq S_w$, so $P^\star\cap Q\leq N_{S_w}(V)=P$. As $P<Q$ we then have  
$P<N_Q(P)$ by 3.2. Suppose that $P^\star=N_Q(P)^\star$. Then 
$$ 
P<N_Q(P)\leq N_Q(P)^\star\cap Q= P^\star\cap Q=P,   
$$ 
for a contradiction. Thus $P^\star<N_Q(P)^\star$, and so $dim(P)<dim(N_Q(P))$. 
\qed 
\enddemo 

For each $w\in\bold W(\ca L)$ there is a conjugation map $c_w:X\to Y$, for any subgroup $X$ of 
$S_w$, and for any subgroup $Y$ of $S$ containing $X^w$ (cf. the discussion following 2.8).

\definition {Definition 3.4} Let $(\ca L,\D,S)$ be a locality. The {\it fusion system} $\ca F_S(\ca L)$ 
is the category whose objects are the subgroups of $S$, and whose morphisms are the conjugation maps 
$$ 
c_w:X\to Y \quad(\text{where $w\in\bold W(\ca L)$, $X\leq S_w$, and $X^w\leq Y$)}.
$$ 
\enddefinition 

The fusion system $\ca F:=\ca F_S(\ca L)$ will be part of the focus of Part II of this series, and it 
will play an important role in the remaining Parts. In this Part I, the only reason for 
introducing it is for the sake of some convenient terminology and notation, as follows. 
\vskip .1in 
\noindent 
$\cdot\ $ If $X\leq S_w$ then $X^w$ is an {\it $\ca F$-conjugate} of $X$, and the set of all 
$\ca F$-conjugates of $X$ is denoted $X^{\ca F}$. 

\vskip .1in 
\noindent 
$\cdot\ $ A subgroup $X$ of $S$ is {\it fully normalized} in $\ca F$ if $dim(N_S(X))\geq dim(N_S(X'))$ for 
every $\ca F$-conjugate $X'$ of $X$. A subgroup $Y$ of $S$ is {\it fully centralized} in $\ca F$ if 
$dim(C_S(Y)Y)\geq dim(C_S(Y')Y')$ for every $\ca F$-conjugate $Y'$ of $Y$. Notice that since $dim(\Omega)$ 
is finite, every subgroup of $S$ has a fully normalized $\ca F$-conjugate and a fully centralized 
$\ca F$-conjugate 
 
\vskip .1in 
\noindent 
$\cdot\ $ A collection $\G$ of subgroups of $S$ is {\it $\ca F$-invariant} if $\G$ is closed with respect to 
$\ca F$-isomorphisms; and $\G$ is {\it $\ca F$-closed} if $\G$ is non-empty and closed with respect to all 
$\ca F$-homomorphisms ($X\in\G$, $w\in\bold W(\ca L)$ with $X\leq S_w$, and $X^w\leq Y\leq S$ implies 
$Y\in\G$). For example, $\D$ is $\ca F$-closed.

\vskip .1in 
From this point forth, $(\ca L,\D,S)$ will be a fixed locality. We write $\ca F$ for $\ca F_S(\ca L)$, and 
we say also that $\ca L$ is a locality {\it on} $\ca F$. Write $\Omega$ for $\Omega_S(\ca L)$. For 
any subgroup $X$ of $S$ write $dim(X)$ for $dim_{\Omega}(X^\star)$. 
\vskip .1in 

The following result is immediate from 2.3(c). 

\proclaim {Lemma 3.5} For $P,Q\in\D$, $Hom_{\ca F}(P,Q)$ is the set of conjugation maps 
$c_g:P\to Q$ such that $g\in\ca L$, $P\leq S_g$, and $P^g\leq Q$. In particular, for each $P\in\D$ 
there exists $g\in\ca L$ such that $P^g$ is fully normalized in $\ca F$. 
\qed
\endproclaim

\definition {Definition 3.6} A $p$-subgroup $R$ of a partial group $\ca G$ is a {\it Sylow $p$-subgroup} if 
\roster 

\item "{(1)}" $R$ is a maximal $p$-subgroup of $\ca G$ and, 

\item "{(2)}" for each $p$-subgroup $P$ of $\ca G$ there exists $g\in\ca G$ such that $P^g$ is a subgroup 
of $R$. 

\endroster 
We write $Syl_p(\ca G)$ for the set (possibly empty) of all Sylow $p$-subgroups of $\ca G$. 
\enddefinition

\proclaim {Lemma 3.7} Set $O_p(\ca L)=\1^\star$. Then $O_p(\ca L)$ is a normal $p$-subgroup of $\ca L$, 
and contains every normal $p$-subgroup of $\ca L$. Moreover, if $\ca L$ is a group then 
$O_p(\ca L)=\bigcap Syl_p(\ca L)$. 
\endproclaim 

\demo {Proof} Let $X$ be a normal $p$-subgroup of $\ca L$. Then $SX$ is a $p$-group, and so $X\leq S$ by 
the maximality of $S$ among the $p$-subgroups of $\ca L$. Then $X\leq S_w$ for all $w\in\bold W(\ca L)$, and 
thus $X\leq\1^\star$. Since $\1^\star\norm\ca L$ by 2.13, we have the desired characterization of 
$O_p(\ca L)$. If $\ca L$ is a group then $Syl_p(\ca L)$ is a conjugacy class of subgroups of $\ca L$, 
and so $O_p(\ca L)=\bigcap Syl_p(\ca L)$ in that case. 
\qed 
\enddemo

\proclaim {Proposition 3.8} Let $(\ca L,\D,S)$ be a locality. Then $S\in Syl_p(\ca L)$. 
\endproclaim  

\demo {Proof} Set $\ca F=\ca F_S(\ca L)$ and $\Omega=\Omega_S(\ca L)$. We first show: 
\roster 

\item "{(1)}" Let $V\leq S$ with $V$ fully normalized in $\ca F$, and assume that for each $p$-subgroup $X$ 
of $N_{\ca L}(V)$ with $N_S(V)\norm X$ there exists $g\in\ca L$ such that $X^g\leq S$. Then $N_S(V)$ is 
a maximal $p$-subgroup of $N_{\ca L}(V)$. 

\endroster 
Indeed, set $R=N_S(V)$, let $X$ be a $p$-subgroup of $N_{\ca L}(V)$ containing $R$, and set $D=N_X(R)$. By 
assumption there exists $g\in\ca L$ with $D^g\leq S$, and so $R^g\leq D^g\leq N_S(V^g)$. As $V$ is fully 
normalized in $\ca F$ we have $dim(R^g)=dim(D^g)=dim(N_S(V^g))$. Then $D^g\leq(R^g)^\star=(R^\star)^g$. 
This yields $D\leq R^\star$, so $D\leq S$, and $D\leq S\cap N_{\ca L}(V)=R$. 
Thus $R=N_X(D)$, and so $R=X$ by 3.2. This completes the proof of (1). 

The proof of 3.8 will now proceed by contradiction. Among all counter-examples, choose $\ca L$ so 
that $dim(\Omega)$ is as small as possible. Here $S$ is a maximal $p$-subgroup of $\ca L$ by (L1), so 
the assumption that $S\notin Syl_p(\ca L)$ means that there exists a $p$-subgroup $X$ of $\ca L$ such 
that no $\ca L$-conjugate of $X$ is contained in $S$. By 2.17 there exists a 
unique largest $P_X\in\D\cap\Omega$ with $X\leq N_{\ca L}(P_X)$. Then $XP_X$ is a $p$-group, and no 
$\ca L$-conjugate of $XP_X$ is contained in $S$. Let $\bold X$ be the set of all $p$-subgroups 
$X$ of $\ca L$ such that $P_X\leq X$ and such that there exists no $g\in\ca L$ with $X^g\leq S$. Thus   
$\bold X$ is non-empty. Choose $X\in\bold X$ so that $dim(P_X)$ is as large as possible, 
and set $P=P_X$. 

Let $w\in\bold W(\ca L)$ such that $P\leq S_w$. Then $w\in\bold D$, and 2.3(c) shows that $X^w$ is defined 
and is a $p$-subgroup 
of $N_{\ca L}(P^w)$. If there exists $g\in\ca L$ such that $(X^w)^g\leq S$ then $X^f\leq S$, where 
$f=\Pi(w\circ(g))$ (and where $w\circ(g)\in\bold D$ via $P$). Thus, we are free to replace $X$ by any 
$X^g$ such that $P^g\leq S$, and we may therefore assume that $X$ has been chosen so that $P$ is fully 
normalized in $\ca F$. Set $R=N_S(P)$ and $H=N_{\ca L}(P)$. 

If $P=S$ then $X\leq S$ by (L1), contrary to 
$X\in\bold X$. Thus $P<S$, and since $P\in\Omega$ we obtain $P<R$ from 3.2. Then $dim(P)<dim(R)$ as $P\in\D$.  
Now let $X'$ be a $p$-subgroup of $H$ such that $R\norm X'$. Then $R\leq P_{X'}$, so the maximality 
of $dim(P)$ in the choice of $X$ yields $X'\notin\bold X$. That is, $X'$ is 
conjugate in $\ca L$ to a subgroup of $S$, and then (1) implies that $R$ is a maximal $p$-subgroup 
of $H$. Then $H$ may be viewed as a locality $(H,\G,R)$, in which $\G$ is the set of all subgroups of $R$. 
As $X\in\bold X$, $X$ is not conjugate in $H$ to a subgroup of $R$, and thus $(H,\G,R)$ is a 
counter-example to the proposition. As $dim(\Omega_R(H))\leq dim(\Omega)$ by 2.16, we may then 
assume that $H=\ca L$. Thus $R=S$, and $P=O_p(H)$ by 3.7.  

Let $\bold M$ be the set of all maximal $p$-subgroups of $H$, and let $\w{\bold M}$ be the set of all 
$S'\in\bold M$ such that there exists a framing $\{H_n\}$ of $H$ by finite subgroups with  
$S'\cap H_n\in Syl_p(H_n)$ for all $n$ (see the remark (2) following 3.1). 
For any framing $\{H_n\}$ of $H$ by finite subgroups we may choose $Q_n\in Syl_p(H_n)$ with 
$Q_n\leq Q_{n+1}$ for all $n$, and then $\cup\{Q_n\}\in\w{\bold M}$. Thus $\w{\bold M}$ is non-empty. 
Noting that $\w{\bold M}$ is invariant under $H$-conjugation, it follows that 
if $S\notin\w{\bold M}$ then no member of $\w{\bold M}$ is conjugate to $S$ in $H$. 
Thus we may choose $S'\in\bold M$ such that $S'$ is not conjugate to $S$, and we may assume that either 
$S$ or $S'$ is in $\w{\bold M}$. Among all such $S'$, choose $S'$ so that $dim(S\cap S')$ is as large as 
possible. Further, we may fix a framing $\{H_n\}$ of $H$ by finite subgroups such that, upon setting   
$S_n=S\cap H_n$ and $S'_n=S'\cap H_n$, either $S_n\in Syl_p(H_n)$ for all $n$ or 
$S'_n\in Syl_p(H_n)$ for all $n$. We are free to replace $S'$ by any $H$-conjugate of $S'$, so: 
\roster 

\item "{(2)}" For any given $n$ we may assume that either $S_n\leq S'_n$ or $S'_n\leq S_n$. 

\endroster 
If $S=P$ then $S\norm H$ then $S$ contains every $p$-subgroup of $\ca L$ by (L1), and then $H$ is not 
a counter-example to 3.8. Thus $P<S$. Thus there exists $n$ with $S_n\nleq P$. Also $S'\nleq P$ since 
$S'\nleq S$, and so for $n$ sufficiently large we have also $S'_n\nleq P$. Then (2) implies that we may 
take $P<S\cap S'$. As $P\in\Omega$ by 3.7, we then have $dim(S\cap S')>dim(P)$.  

Set $Y=(S\cap S')^\star$. Then there exists $x\in G$ with $Y^x$ fully normalized in $\ca F$. Upon replacing 
$S'$ with $(S')^x$, we may therefore assume that $Y$ is fully normalized in $\ca F$. Then (1) implies that  
$N_S(Y)$ is a maximal $p$-subgroup of $N_H(Y)$. Here $dim(\Omega_{N_S(Y)}(N_H(Y))<dim(\Omega)$ by 
2.16(b), so the locality $N_H(Y)$ is not a counterexample, and $N_S(Y)$ is a Sylow $p$-subgroup of 
$N_H(Y)$. Set $T=N_{S'}(S\cap S')$. 
Then $T\leq N_H(Y)$, so there exists $h\in N_H(Y)$ with $X^h\leq N_S(Y)$. The maximality condition on 
$dim(S\cap S')$ then yields $dim(S\cap S')=dim(S\cap(S')^h)$. Then, since 
$(S\cap S')^h\leq T^h\leq S\cap(S')^h$, we conclude that  
$$ 
Y^h=(T^h)^\star=(T^\star)^h\geq T^h. 
$$ 
Conjugation by $h\i$ now yields $Y\geq T$, and so $T=S\cap S'$. We now 
appeal to 3.2, and conclude that $S\cap S'=S'$. The maximality of $S'$ then yields 
$S'=S$, thereby providing a contradiction and completing the proof. 
\qed 
\enddemo

\proclaim {Lemma 3.9} Let $X\leq S$ with $X$ fully normalized in $\ca F$.  
Then $N_S(V)$ is a maximal $p$-subgroup of $N_{\ca L}(V)$. Similarly, if $Y\leq S$ is fully centralized 
in $\ca F$ then $C_S(Y)$ is a maximal $p$-subgroup of $C_{\ca L}(Y)$.  
\endproclaim 

\demo {Proof} Set $R=N_S(X)$, let $A$ be a $p$-subgroup of $N_{\ca L}(X)$ containing $R$, and set 
$B=N_A(R)$. By 3.8 there exists $g\in\ca L$ with $B^g\leq S$, and then $B^g\leq N_S(X^g)$. As $X$ is 
fully normalized in $\ca F$ we have $dim(R)\geq dim(B^g)$, so $dim(R^g)\geq dim(B^g)$, and then  
$(R^g)^\star=(B^g)^\star$ since $R^g\leq B^g$. As $(R^g)^\star=(R^\star)^g$, conjugation by $g\i$ yields 
$B\leq R^\star$. Thus $X\norm B\leq S$, and so $B=R$. Then $A=R$ by 3.2, and $R$ is a maximal $p$-subgroup 
of $N_{\ca L}(V)$. The proof that $C_S(Y)$ is a maximal $p$-subgroup of $C_{\ca L}(Y)$ is similar, and may 
be omitted. 
\qed 
\enddemo 

\proclaim {Lemma 3.10} Let $P\in\D$. Then $N_S(P)\in Syl_p(N_{\ca L}(P))$ if and only if $P$ is 
fully normalized in $\ca F$. 
\endproclaim 

\demo {Proof} Set $R=N_S(P)$ and $H=N_{\ca L}(P)$. Suppose that $P$ is fully normalized in $\ca F$. 
Then $R$ is a maximal $p$-subgroup of $H$ by 3.9, and then $(H,\G,R)$ is a locality, where $\G$ is the 
set of all subgroups of $R$. Then $R\in Syl_p(H)$ by 3.8. 

For the converse; suppose that $R\in Syl_p(H)$ 
and let $P'$ be an $\ca F$-conjugate of $P$. Then $P'=P^g$ for some $g\in\ca L$ by 2.3(c), and 
conjugation by $g$ is an isomorphism $H\to N_{\ca L}(P')$ by 2.3(b). Set $T=N_S(P')^{g\i}$. Then $T$ is a 
$p$-subgroup of $H$, and so there exists $h\in H$ with $(T^g)^h\leq R$. As $(g\i,h)\in\bold D$ via $P'$, 
we have $T^{gh}\leq R$. Thus $dim(T)\leq dim(R)$, and so $P$ is fully normalized in $\ca F$. 
\qed 
\enddemo 

\proclaim {Lemma 3.11} Let $\ca H$ be a subset of $\ca L$ such that $\ca H$ is a partial group, and such 
that the inclusion map $\ca H\to\ca L$ is a homomorphism of partial groups. Set $R=S\cap\ca H$, and set 
$$ 
\Omega_R(\ca H)=\{R\cap S_w\mid w\in\bold W(\ca H)\}. 
$$ 
Assume: 
\roster 

\item "{(1)}" $R$ is a subgroup of $\ca H$, and 

\item "{(2)}" if $x\in R$, $h\in\ca H$, and $(h\i,x,h)\in\bold D(\ca L)$, then 
$(h\i,x,h)\in\bold D(\ca H)$. 

\endroster 
Then $R$ is a $p$-subgroup of $\ca H$, and the poset $\Omega_R(\ca H)$ (partially ordered by 
inclusion) is finite-dimensional. 
\endproclaim 

\demo {Proof} The inclusion map of $R$ into $S$ is a homomorphism, so (1) imples that $R$ is a 
$p$-subgroup of $\ca H$. Let $X=\in\Omega_R(\ca H)$. Evidently $X\leq X^*$ and $X\leq R$. But also,  
(2) implies  
$$ 
X=\bigcap\{S_w\cap R\mid w\in\bold W(\ca H),\ X\leq S_w\},  
$$ 
and hence $X^*\cap R\leq X$. Thus $X^*\cap R=X$, so the mapping $X\maps X^*$ of $\Omega_R(\ca H)$ 
into $\Omega_S(\ca L)$ is injective, and so $\Omega_R(\ca H)$ is finite-dimensional. 
\qed 
\enddemo

\proclaim {Lemma 3.12} Let $H$ be a subgroup of $\ca L$. Then $Syl_p(H)\neq\nset$, and 
$Syl_p(H)$ is the set of maximal $p$-subgroups of $H$. Moreover, for any $R\in Syl_p(H)$ there exists a 
framing $\{H_n\}$ of $H$ by finite subgroups such that $R\cap H_n\in Syl_p(H_n)$ for all $n$. 
\endproclaim 

\demo {Proof} Let $R$ be a maximal $p$-subgroup of $H$ and let $\G$ be the set of all subgroups of $R$. 
Then $(H,\G,R)$ is a locality, and so $R\in Syl_p(H)$. By the remark following 3.1 one may choose 
$R$ so that $R\cap H_n\in Syl_p(H_n)$ for all $n$. One observes that for any $h\in H$, $\{(H_n)^h\}$ is 
again a framing of $H$. Since $Syl_p(H)$ forms a single $H$-conjugacy class we then have the desired 
``compatibility" with a suitable framing, for any Sylow $p$-subgroup of $H$. 
\qed 
\enddemo 

The proofs of the three parts of the following result are essentialy the same as the corresponding 
(elementary) proofs for finite groups, and they are therefore omitted. 

\proclaim {Lemma 3.13} Let $H$ be a subgroup of $\ca L$, let $K\norm H$ be a normal subgroup of $H$, 
let $R$ be a Sylow $p$-subgroup of $H$, and set $T=R\cap K$. Then the following hold. 
\roster 

\item "{(a)}" $T\in Syl_p(K)$ and $G=N_G(T)K$. 

\item "{(b)}" For any subgroup $Q$ of $R$ with $T\leq Q$, we have $Q\in Syl_p(KQ)$. 

\item "{(c)}" $RK/K\in Syl_p(H/K)$. 

\endroster 
\qed 
\endproclaim

The next three results concern localities formed within $\ca L$.

\proclaim {Lemma 3.14} Let $T\norm S$ be a normal subgroup of $S$, and set $\ca L_T=N_{\ca L}(T)$. Then 
$(\ca L_T,\D,S)$ is a locality. 
\endproclaim 

\demo {Proof} That $\ca L_T$ is a partial subgroup of $\ca L$ follows from 2.3(c). That $(\ca L_T,\D)$ is
objective is a consequence of 2.10. The remaining points of verification are inherited from $\ca L$ in 
a straightforward way. 
\qed 
\enddemo

Recall from definition 3.4 the notion of an $\ca F$-closed collection $\G$ of subgroups of $S$. The reader 
should have no difficulty in working through the steps which verify the following result. 

\proclaim {Lemma 3.15} Let $\G$ be an $\ca F$-closed subset of $\D$. Set 
$$ 
\ca L\mid_\G=\{g\in\ca L\mid S_g\in\G\}, \quad\text{and}\quad 
\bold D\mid_\G=\{w\in\bold W(\ca L)\mid S_w\in\G\}. 
$$ 
Then $\ca L\mid_\G$ is a partial group by restriction of the product in $\ca L$ to $\bold D\mid_\G$, and  
by restriction of the inversion in $\ca L$. Moreover, $(\ca L\mid_\G,\G,S)$ is a locality   
$($to be called {\it the restriction of $\ca L$ to $\G)$. 
\qed 
\endproclaim

\proclaim {Lemma 3.16} Let $(\ca L,\D,S)$ be a locality, and let $\w\D$ be the set of all 
subgroups $P\leq S$ such that $P^\star\in\D$. Then $(\ca L,\w\D,S)$ is a locality on $\ca F$.   
\endproclaim 

\demo {Proof} It follows from 2.13 that $\w\D$ is $\ca F$-invariant, and then from 2.12(b) that $\w\D$ is 
$\ca F$-closed. Thus $(\ca L,\w\D)$ satisfies condition (O2) in definition 2.1.  
Also, 2.13 shows that condition (O1) is satisfied, and so $(\ca L,\w\D)$ is an objective partial group. 
The conditions that $(\ca L,\w\D,S)$ must then fulfill in order to be a locality are given by the 
partial group structure of $\ca L$. The $\ca F$-homomorphisms are compositions of restrictions of 
conjugation maps $c_g:S_g\to S$, where $S_g$ is determined by the structure of $\ca L$ as a partial 
group, so $(\ca L,\w\D,S)$ is a locality on $\ca F$. 
\qed 
\enddemo

We end this section by investigating the situation in which $\ca L$ is a locality in more than 
one way.

\proclaim {Lemma 3.17} Let $(\ca L,\D,S)$ be a locality, let $(x_1,\cdots,x_n)\in\bold D$, and let 
$(g_0,\cdots,g_n)\in\bold W(N_{\ca L}(S))$. Then 
$$ 
(g_0,x_1,g_1\i,g_1,x_2,\cdots,x_{n_1},g_{n-1}\i,g_{n-1},x_n,g_n)\in\bold D. \tag(*)
$$ 
In particular, conjugation by $g\in N_{\ca L}(S)$ is an automorphism of the partial group $\ca L$. 
\endproclaim 

\demo {Proof} Set $w=(x_1,\cdots,x_n)$ and set $P=(S_w)^{g_0\i}$. Then the word displayed in (*) is in 
$\bold D$ via $P$. 
\qed 
\enddemo

\definition {Definition 3.18} Let $\ca L$ be a partial group. An automorphism $\a$ of $\ca L$ is {\it inner}  
if there exists $g\in\ca L$ such that $\a$ is given by conjugation by $g$. That is: 
\roster 

\item "{(1)}" $x\in\ca L$ $\implies$ $x^g$ is defined and is equal to $x\a$, and  

\item "{(2)}" $(x_1,\cdots,x_n)\in\bold D$ $\implies$ $(x_1^g,\cdots,x_n^g)\in\bold D$ 
and $x_1^g\cdots x_n^g=(x_1\cdots x_n)^g$. 

\endroster 
Write $Inn(\ca L)$ for the set of all inner automorphisms of $\ca L$. 
\enddefinition

\proclaim {Proposition 3.19} Let $(\ca L,\D,S)$ be a locality and let $K$ be 
the set of all $g\in\ca L$ such that conjugation by $g$ is an inner automorphism of $\ca L$. Let 
$\ca S$ be the set of all subgroups $S'$ of $\ca L$ having the property that for some set $\D'$ of 
subgroups of $S'$, the conditions (L1) and (L2) in definition 3.1 hold with $(\D',S')$ in place of 
$(\D,S)$. 
\roster 

\item "{(a)}" $K$ is a subgroup of $\ca L$.

\item "{(b)}" For each $w=(x_1,\cdots,x_n)\in\bold D$ and each $(g_0,\cdots,g_n)\in\bold W(K)$, 
we have 
$$ 
(g_0,x_1,g_1\i,g_1,x_2,\cdots,g_{n-1}\i,g_{n-1},x_n,g_n)\in\bold D. \tag*
$$

\item "{(c)}" Set $Q=O_p(K)$. Then $Q\in\D$ and $K=N_{\ca L}(Q)$. Moreover, for each $w\in\bold D$ 
there exists $P\leq Q\cap S_w$ such that $P\in\D$ and $P^w\leq Q$. 

\item "{(d)}" $\ca S=Syl_p(K)$. In particular, $K$ acts transitively on $\ca S$ by conjugation. 

\item "{(e)}" $Inn(\ca L)\norm Aut(\ca L)$, and every automorphism of $\ca L$ can be factored as an 
inner automorphism followed by an automorphism which leaves $S$ invariant. 

\endroster 
\endproclaim 

\demo {Proof} For $S'\in\ca S$ note that $N_{\ca L}(S')$ is a 
subgroup of $\ca L$. The condition (*) in 3.17 then holds with $N_{\ca L}(S')$ in the role of $K$. 
Let $\ca U$ be the union (taken over all $S'\in\ca S$) of the groups $N_{\ca L}(S')$, and let $H$ be the 
partial subgroup of $\ca L$ generated by $\ca U$. The above observation concerning 3.17 together with a 
straightforward argument by induction on word-length then yields $\bold W(\ca U)\sub\bold D$, 
and thus $H$ is a subgroup of $\ca L$. Moreover, $\bold W(\ca U)$ has the following property:  
For any $(g_1,\cdots,g_n)\in\bold D(\ca L)$, and any $(u_0,\cdots,u_n)$ with each $u_i\in\bold W(\ca U)$,  
the word  
$$ 
u_0\circ g_1\circ u_1\i\circ u_1\circ g_2\circ\cdots\circ u_{n-1}\i\circ u_{n-1}\circ g_n\circ u_n 
$$ 
is in $\bold D$, by induction on the sum of the lengths of the words $u_i$. Condition (*) in (b), above, 
then holds for all $(g_0,\cdots,g_n)\in\bold W(H)$, by $\bold D$-associativity. In particular, 
conjugation by $g\in H$ is an automorphism of $\ca L$, and thus $H$ is a subset of $K$. 

Now let $g\in K$. Then $(\ca L,\D^g,S^g)$ is a locality, so $S^g\in\ca S$. Then $S^g\in Syl_p(H)$, and 
there exists $h\in H$ with $S^g=S^h$. The product $gh\i$ is defined, and 
then $gh\i\in N_{\ca L}(S)$. Thus $gh\i\in H$, so $g\in H$, and we conclude that $H=K$. This completes 
the proof of (a), (b), and (d). Point (e) is immediate from (d), so it remains only to prove (c).   

By 2.17 there is a largest $Q\in\D$ such that $K$ normalizes $Q$. Then, since $S\in Syl_p(K)$, we obtain 
$Q:=O_p(K)$. Also by 2.17 there exists a word $u\in\bold W(K)$ such that 
$Q=S_u$. Then also $Q=S_{u\i}$. Let $w\in\bold D$, and let $f\in N_{\ca L}(Q)$. 
The word $v:=u\circ w\circ u$ is in $\bold D$ by (*), and  
clearly $S_v\leq Q$ and $S_{v\i}\leq Q$. Set $P=(S_v)^{\Pi(u)}$. Then $P\leq S_w\cap Q$ and 
$P^{\Pi(w)}\leq Q$. As $S_v\in\D$ we have $P\in\D$. Write $w=(x_1,\cdots,x_n)$. Then 
$$ 
(f\i,x_1,f,f\i,\cdots,f,f\i,x_n,f)\in\bold D\ \ \text{via $P$}, 
$$ 
and this shows that conjugation by $f$ is an inner automorphism of $\ca L$. Thus $K=N_{\ca L}(Q)$, 
and the proof is complete. 
\qed 
\enddemo 

The proof of the following corollary is left to the reader. 

\proclaim {Corollary 3.20} Let $(\ca L,\D,S)$ be a locality on $\ca F$, and define $\ca S'$ and $K$ 
as in 3.19. Let $S'\in\ca S$, let $g\in K$ with $S'=S^g$, and set $\ca F'=\ca F_{S'}(\ca L)$. 
Let $(\Omega,\star)$ and $(\Omega',\star')$ be the stratifications on $\ca F$ and on $\ca F'$, 
respectively, induced from $\ca L$.  Then conjugation by $g$ induces isomorphisms $\ca F\to\ca F'$ and 
$\Omega\to\Omega'$. Moreover, for each subgroup $X\leq S$ we have $(X^\star)^g=(X^g)^{\star'}$. 
\qed 
\endproclaim

\proclaim {Lemma 3.21} Let $\ca L$ be a locality and let $S$ be a Sylow $p$-subgroup of $\ca L$. 
Then there is a unique smallest set $\D=\D_0$ and a unique largest set $\D=\D_1$ of subgroups of $S$ 
such that the conditions (L1) and (L2) of 3.1 are satisfied by $(\ca L,\D,S)$.  
\endproclaim 

\demo {Proof} Set Take $\D_0$ to be the overgroup-closure in $S$ of the set of all $S_w$ for 
$w\in\bold D$. Take $\D_1$ to be the union of all the the sets $\G$ of subgroups of $S$ which fulfill 
(L1) and (L2). 
\qed 
\enddemo

\vskip .2in 
\noindent 
{\bf Section 4: Partial normal subgroups} 
\vskip .1in 

Throughout this section we fix a locality $(\ca L,\D,S)$ and a partial normal subgroup $\ca N\norm\ca L$. 
Recall that this means that $\ca N$ is a partial subgroup of $\ca L$ and that $\Pi(g\i,x,g)\in\ca N$ 
for all $x\in\ca N$ and all $g\in\ca L$ for which $(g\i,x,g)\in\bold D$. Set $T=S\cap\ca N$.  

\vskip .1in 
Recall from 3.4 that there is a fusion system $\ca F=\ca F_S(\ca L)$ associated with $\ca L$. A subgroup 
$R$ of $S$ is {\it strongly closed} in $\ca F$ if $x^w\in R$ whenever $w\in\bold W(\ca L)$ with 
$x\in R\cap S_w$. 

\proclaim {Lemma 4.1} 
\roster

\item "{(a)}" $T$ is strongly closed in $\ca F$. 

\item "{(b)}" Let $x\in\ca N$ and let $P$ be a subgroup of $S_x$. Then $PT=P^x T$.

\item "{(c)}" $T$ is maximal in the poset of $p$-subgroups of $\ca N$. 

\endroster
\endproclaim

\demo {Proof} (a) Let $g\in\ca L$ and let $t\in S_g\cap T$. Then $t^g\in S$, and $t^g\in\ca N$ as 
$\ca N\norm\ca L$. Thus $t^g\in T$. Iteration of such conjugation maps shows that $t^w\in T$ for all 
$w\in\bold W(\ca L)$ with $t\in S_w$. 
\vskip .1in 
\noindent 
(b) Let $a\in P$. Then $(P^x)^a\leq S$ and $P^a=P$. Setting $w=(a\i,x\i,a,x)$ we then have $w\in\bold D$ 
via $P^{xa}$. Now $\Pi(w)=a\i a^x\in S$, while also $\Pi(w)=(x\i)^ax\in\ca N$, and so $\Pi(w)\in T$. 
Then $a^x\in aT$, and we have thus shown that $P^x\leq PT$. Then $P^x T\leq PT$. The equality $P^xT=PT$  
follows via symmetry, with $x\i$ and $P^x$ in place of $x$ and $P$. 
\vskip .1in 
\noindent 
(c) Let $X$ be a $p$-subgroup of $\ca N$ containing $T$. By  3.8 there exists 
$g\in\ca L$ with $X^g\leq S$, and then $X^g\leq S\cap\ca N=T$. Set $P=S_g$ and $Q=P^g$. Then $X^g\leq Q$ 
and conjugation by $g\i$ yields $X\leq P$. Thus $X\leq P\cap\ca N=T$, and so $X=T$. 
\qed
\enddemo

\proclaim {Lemma 4.2} Let $x,y\in\ca N$ and let $f\in N_{\ca L}(T)$. 
\roster 

\item "{(a)}" If $(x,f)\in\bold D$ then $(f,f\i,x,f)\in\bold D$, $xf=fx^f$, 
and $S_{(x,f)}=S_{(f,x^f)}=S_x\cap S_f$. 

\item "{(b)}" If $(f,y)\in\bold D$ then $(f,y,f\i,f)\in\bold D$, 
$fy=y^{f\i}f$, and $S_{(f,y)}=S_{(y^{f\i},f)}=S_{y^{f\i}}\cap S_f$. 

\endroster 
\endproclaim 

\demo {Proof} (a): Set $Q=S_{(x,f)}$. Then $T\leq S_f$ by hypothesis. Since $Q^xT=QT$ by 4.1(b), 
we then have $Q\leq S_f$. Thus $Q\leq P:=S_x\cap S_f$. But also $P^xT=PX$, so $P=Q$. Moreover,  
$(f,f\i,x,f)\in\bold D$ via $Q$, and then $\Pi(f,f\i,x,f)=xf=fx^f$. As 
$Q=S_f\cap S_x\leq S_{(f,f\i,x,f)}$ we obtain $S_{(x,f)}=S_{(f,x^f)}$. Thus, (a) holds. 

For point (b): Set $R=S_{(f,y)}$. Then $R^{fy}T=R^fT\leq S_{f\i}$, so 
$(f,y,f\i,f)\in\bold D$ via $R$, and $fy=y^{f\i}f$. The remainder of (b) now 
follows as an application of (a) to $(y^{f\i},f)$. 
\qed 
\enddemo

\proclaim {Lemma 4.3} Let $w\in\bold W(N_{\ca L}(T))\cap\bold D$, set $g=\Pi(w)$, and let $x,y\in\ca N$. 
\roster 

\item "{(a)}" Suppose $(x)\circ w\in\bold D$ and set $P=S_{(x)\circ w}$. Then 
$u:=w\i\circ(x)\circ w\in\bold D$, and $S_u=P^g$. 

\item "{(b)}" Suppose that $w\circ(y)\in\bold D$ and set $Q=S_{w\circ(y)}$. Then 
$v:=w\circ(y)\circ w\i\in\bold D$, and $S_v=Q$. 

\endroster 
\endproclaim 

\demo {Proof} We prove only (a), leaving it to the reader to supply a similar argument for (b). 
As $P^xT=PT$ by 4.1(b), and since both $P^x$ and $T$ are contained in $S_w$, we obtain $P\leq S_w$. 
Then $P^g\leq S$, and $P^g\leq S_u$. In particular $u\in\bold D$ via $P^g$. As  
$S_u\leq S_{w\i}$, and $(S_u)^{g\i}\leq P$, we obtain $S_u\leq P^g$. Thus $S_u=P^g$.  
\qed 
\enddemo

\proclaim {Lemma 4.4 (``Frattini Calculus")} Let $w=(f_1,g_1,\cdots,f_n,g_n)\in\bold D$ via $P$, with 
$f_i\in N_{\ca L}(T)$ and with $g_i\in\ca N$ for all $i$. Set $u_n=g_n$, and for all $i$ with 
$1\leq i<n$ set 
$$ 
u_i=(f_n\i,\cdots,f_{i+1}\i,g_i,f_{i+1},\cdots,f_n). 
$$ 
Then $u_i\in\bold D$ via $P$ and, upon setting $\bar g_i=\Pi(u_i)$ and 
$w'=(f_1,\cdots,f_n,\bar g_1,\cdots,\bar g_n)$, we have $w'\in\bold D$, 
$S_w=S_{w'}$, and $\Pi(w)=\Pi(w')$. Thus:  
\roster 

\item "{(a)}" $\Pi(w)=(f_1\cdots f_n)(\bar g_1\cdots\bar g_n)$, where $\bar g_i=(g_i)^{f_{i+1}\dots f_n}$. 

\endroster  
Moreover, we similarly have: 
\roster 

\item "{(b)}" $\Pi(w)=(\w g_1\cdots\w g_n)(f_1\cdots f_n)$ where $\w g_i=(g_i)^{(f_1\cdots f_i)\i}$.

\endroster 
\endproclaim 

\demo {Proof} For the proof of (a), suppose that we are given an index $k$ such that 
$$ 
w_k=(f_1,g_1,\cdots,f_k,g_k)\circ(f_{k+1},\cdots,f_n)\circ(\bar g_{k+1},\cdots,\bar g_n)\in\bold D,  
$$ 
with $P:=S_w=S_{w_k}$, $\Pi(w)=\Pi(w_k)$, $u_i\in\bold D$ for all $i$ with 
$i>k$, and where $\bar g_i=\Pi(u_i)$. For example, $k=n$ is such an index. If  
$k=0$ then (a) holds. So assume $k>0$. Set 
$$ 
v=(f_1,g_1,\cdots,f_{k-1},g_{k-1},f_k). 
$$ 
Then 4.3(a) as applied to $(g_k)\circ(f_{k+1},\cdots,f_n)$ yields 
$u_k\in\bold D$ via $Q$, where $Q$ is the image of $P$ under conjugation by 
$\Pi(v)$. Then also $w_{k-1}\in\bold D$ via $P$, $S_w=S_{w_{k-1}}$, and 
1.1(3) yields $\Pi(w_{k-1})=\Pi(w_k)$. Iteration of this procedure yields 
(a). A similar procedure involving 4.3(b) yields (b).  
\qed 
\enddemo 

\definition {Definition 4.5} Let $\ca L\circ\D$ be the set of all pairs $(f,P)\in\ca L\times\D$ such that 
$P\leq S_f$. Define a relation $\uparrow$  on $\ca L\circ\D$ by $(f,P)\uparrow(g,Q)$ if there exist elements 
$x\in N_{\ca N}(P,Q)$ and $y\in N_{\ca N}(P^f,Q^g)$ such that $xg=fy$. 
\enddefinition 

This relation may be indicated by means of a commutative square: 
$$ 
\CD 
Q@>g>>Q^g \\ 
@AxAA  @AAyA  \\ 
P@>f>>P^f 
\endCD \tag* 
$$
of conjugation maps, labeled by the conjugating elements, and in which the 
horizontal arrows are isomorphisms and the vertical arrows are injective 
homomorphisms. The relation $(f,P)\up(g,Q)$ may also be expressed by: 
$$
\text{$w:=(x,g,y\i,f\i)\in\bold D$ via $P$, and $\Pi(w)=\1$.} 
$$ 
\vskip .1in
It is easy to see that $\up$ is reflexive and transitive. We say that $(f,P)$ is {\it maximal} in 
$\ca L\circ\D$ if $(f,P)\up(g,Q)$ implies $dim(P)=dim(Q)$. As $\ca L$ is finite-dimensional 
there exist maximal elements in $\ca L\circ\D$. Since $(f,P)\up(f,S_f)$ for $(f,P)\in\ca L\circ\D$, 
it follows that $P=S_f$ for every 
maximal $(f,P)$. For this reason, we will say that $f$ is {\it $\up$-maximal}
in $\ca L$ (with respect to $\ca N$) if $(f,S_f)$ is maximal in $\ca L\circ\D$.

\proclaim {Lemma 4.6} Let $f\in\ca L$. 
\roster 

\item "{(a)}" If $f\in N_{\ca L}(S)$ then $f$ is $\up$-maximal. 

\item "{(b)}" If $f$ is $\up$-maximal then so is $f\i$. 

\item "{(c)}" If $f$ is $\up$-maximal and $(f,S_f)\up(g,Q)$, then $g$ is $\up$-maximal and 
$Q=S_g$. 

\endroster  
\endproclaim 

\demo {Proof} Point (a) is immediate from definition 4.5. Now suppose that $f$ is $\up$-maximal, 
and let $g\in\ca L$ with $(f\i,S_{f\i})\up(g\i,S_{g\i})$. Since $S_{f\i}=(S_f)^f$ and 
$S_{g\i}=(S_g)^g$, there is a diagram 
$$ 
\CD 
(S_g)^g@>g\i>>S_g \\ 
@AxAA      @AAyA \\ 
(S_f)^f@>f\i>>S_f 
\endCD 
$$ 
as in definition 4.5, from which it is easy to read off the relation 
$(f,S_f)\up(g,S_g)$. Then $dim(S_f)=dim(S_g)$ as $f$ is $\up$-maximal. As $S_{x\i}=(S_x)^x$ for any 
$x\in\ca L$ we obtain also 
$dim(S_{f\i})=dim(S_{g\i})$. Thus $f\i$ is $\up$-maximal, and (b) holds. Point (c) is 
immediate from the transitivity of $\up$. 
\qed 
\enddemo 

\proclaim {Lemma 4.7} Let $(g,Q),(h,R)\in\ca L\circ\D$ with $(g,Q)\up(h,R)$, and suppose that 
$T\leq R$. Then there exists a unique $y\in\ca N$ with $g=yh$. Moreover: 
\roster 

\item "{(a)}" $Q^y\leq R$, and $Q\leq S_{(y,h)}$. 

\item "{(b)}" If $N_T(Q^g)\in Syl_p(N_{\ca N}(Q^g))$, then $N_T(Q^y)\in Syl_p(N_{\ca N}(Q^y))$. 

\endroster
\endproclaim

\demo {Proof} By definition of the relation $\up$, there exist elements $u\in N_{\ca N}(Q,R)$ and 
$v\in N_{\ca N}(Q^g,R^h)$ such that $(u,h,v\i,g\i)\in\bold D$ via $Q$, and such that $\Pi(w)=\1$. 
$$ 
\CD 
R@>h>>R^h \\ 
@AuAA  @AAvA \\  
Q@>>g>Q^g 
\endCD 
$$
In particular, $uh=gv$. Since $T\leq R$, points (a) and (b) of 4.1 yield 
$$
\text{$T=T^h$, $Q^uT=QT\leq R$, and $Q^gT=Q^{gv}T\leq R^h$}. 
$$
Then
$$
w:=(u,h,v\i,h\i)\in\bold D\quad{\text via}\quad (Q,Q^u,Q^{uh},Q^{uhv\i}=Q^g,Q^{gh\i}). 
$$
Set $y=\Pi(w)$. Then $y=u(v\i)^{h\i}\in N_{\ca N}(Q,R)$. By 1.1 both 
$(u,h,v\i,h\i,h)$ and $(g,v,v\i)$ are in $\bold D$, and hence  
$yh=uhv\i=g$. This yields (a). The uniqueness of $y$ is given by right cancellation. 

Suppose now that $N_T(Q^g)\in Syl_p(N_{\ca N}(Q^g))$. As $N_T(Q^y)^h=N_T(Q^g)$, it follows from 
2.3(b) that $N_T(Q^y)\in Syl_p(N_{\ca N}(Q^y))$. 
\qed 
\enddemo

\proclaim {Proposition 4.8} Let $g\in\ca L$ and suppose that $g$ is $\up$-maximal with respect to 
$\ca N$. Then $T\leq S_g$. 
\endproclaim 

\demo {Proof} Set $P=S_g$ and $Q=P^g$. We first show: 
\roster 

\item "{(1)}" Let $y\in N_{\ca N}(P,S)$. Then $(g,P)\up(y\i g,P^y)$, and $y\i g$ is $\up$-maximal.

\endroster 
Indeed, we have a commutative diagram 
$$ 
\CD 
P^y@>y^{-1}g>> P^g \\ 
@AyAA          @AA\1A \\ 
P@>>g> P^g 
\endCD 
$$ 
as in 4.5, and then (1) is given by 4.6(c). 

Suppose next that $N_T(P)\in Syl_p(N_{\ca N}(P))$. Then $N_T(P)^g\in Syl_p(N_{\ca N}(Q))$ by 2.3(b), 
and so there exists $x\in N_{\ca N}(Q)$ such that $N_T(Q)\leq (N_T(P)^g)^x$. Here $(g,x)\in\bold D$ via $P$, 
so $(N_T(P)^g)^x=N_T(P)^{gx}$, and then $(g,P)\up(gx,N_T(P)P)$. As $g$ is $\up$-maximal, we conclude that 
$N_T(P)\leq P$. Then (L2), as applied to the inclusion $P\leq TP$, yields $T\leq P$. Thus: 
\roster 

\item "{(2)}" If $T\nleq P$ then $N_T(P)\notin Syl_p(N_{\ca N}(P))$. 

\endroster 

We next show: 
\roster 

\item "{(3)}" Suppose that there exists $y\in\ca N$ such that $P\leq S_y$ and such that 
$N_T(P^y)\in Syl_p(N_{\ca N}(P^y))$. Then $T\leq P$. 

\endroster 
Indeed, under the hypothesis of (3) we have $(g,P)\up(y\i g,P^y)$ by (1). Then 
$(y\i g,P^y)$ is $\up$-maximal and $P^y=S_{y\i g}$. If $T\nleq P$ then $T\nleq P^y$, and then (2) 
applies to $(y\i g,P^y)$ in the role of $(g,P)$ and yields a contradiction. So, (3) holds. 

With these preliminaries in place, we now assume that $g$ is a counter-example (i.e. $T\nleq P$) 
with $dim(P)$ is as large as possible. By 3.5 and 3.10 there exists 
$f\in\ca L$ such that $Q^f\leq S$ and such that $N_S(Q^f)\in Syl_p(N_{\ca L}(Q^f))$. Set $h=gf$ 
(where the product is defined via $P$) and set $R=P^h$. Let $(h',P')$ be maximal in $\ca L\circ\D$ 
with $(h,P)\up(h',P')$, and set $R'=R^{h'}$. If $T\leq P'$ then 4.7 yields $h=yh'$ for some $y\in\ca N$ 
such that $P\leq S_y$ and such that $N_T(P^y)\in Syl_p(N_{\ca N}(P^y))$. The existence of 
such an element $y$ contradicts (3), so we conclude that $T\nleq P'$. Then $(h',P')$ is a counter-example 
to the proposition, and the maximality of $dim(P)$ yields $dim(P)=dim(P')$. Then $(h,P)$ is maximal in 
$\ca L\circ\D$, as is $(h\i,R)$ by 4.6(b), and so $T\nleq R$. But $N_T(R)\in Syl_p(N_{\ca N}(R))$ since 
$R=Q^f$, and so (2) applies with $(h\i,R)$ in the role of $(g,P)$. Thus $T\leq R$, so $T\leq P$, 
and the proof is complete. 
\qed
\enddemo

\proclaim {Lemma 4.9} Suppose that $S=C_S(T)T$ and that $N_{\ca N}(T)\leq N_{\ca L}(S)$. 
Then every element of $N_{\ca L}(T)$ is $\up$-maximal with respect to $\ca N$. 
\endproclaim 

\demo {Proof} Let $f\in N_{\ca L}(T)$ and set $P=S_f$. Then $T\leq P$ and $T\leq P^f$. Let 
$(g,Q)\in\ca L\circ\D$ with $(f,P)\up(g,Q)$, and let $x,y\in\ca N$ be chosen as in definition 4.5. Then 
$P\leq Q$ by 4.8, and we have $P=P^x$ and $(P^f)^y=P^f\leq Q^g$ by 4.1(b). In order to show that $f$ is 
$\up$-maximal it suffices now to show that $P=Q$, and hence it suffices to show that $N_Q(P)\leq P$. 

Set $D=N_Q(P)$. Then $D=C_D(T)T$ by the hypothesis on $S$, and then $[C_D(T),x]$ may be computed in the 
group $N_{\ca L}(P)$. By assumption, $x\in N_{\ca L}(S)$, so $[C_D(T),x]\leq C_S(T)\cap\ca N=Z(T)$, and 
hence $x\in N_{\ca L}(D)$. Similarly $y\in N_{\ca L}(D^g)$, and so $(x,g,y\i)\in\bold D$ via $D$. As $xg=fy$ 
by the setup of definition 4.5, we have $xgy\i=f$ by 1.3(g), and thus $D\leq S_f$. That is,  
$N_Q(P)\leq P$ as required. 
\qed 
\enddemo

If $X$ and $Y$ are subsets of $\ca L$ then one has the notion of the product $XY$, introduced in section 1, 
as the set of all $\Pi(x,y)$ with $(x,y)\in\bold D\cap(X\times Y)$.

\proclaim {Corollary 4.10 (Frattini Lemma)} Let $(\ca L,\D,S)$ be a locality, let $\ca N\norm\ca L$ be a 
partial normal subgroup, and let $\L$ be the set of $\up$-maximal elements of $\ca L$ with respect to 
$\ca N$. Then $\ca L=\ca N\L=\L\ca N$. In particular, we have $\ca L=N_{\ca L}(T)\ca N=\ca N N_{\ca L}(T)$. 
\endproclaim 

\demo {Proof} Let $f\in\ca L$, set $P=S_f$, and choose $(g,Q)\in\ca L\circ\D$ so that 
$(f,P)\up(g,Q)$ and so that $g$ is $\up$-maximal. Then 4.6 yields $f=xg$ for some $x\in\ca N$, 
and then 4.2 shows that $f=gy$ where $y=x^g$. 
\qed 
\enddemo

The following result is fundamental to the theory being developed here. It was 
discovered and proved by Bernd Stellmacher, in his reading of an early draft of [Ch1]. The proof given 
here is essentially his.

\proclaim {Lemma 4.11 (Splitting Lemma)} Let $(x,f)\in\bold D$ with $x\in\ca N$ and with 
$f$ $\up$-maximal. Then $S_{(x,f)}=S_{xf}=S_{(f,x^f)}$. 
\endproclaim 

\demo {Proof} Appealing to 4.2: Set $y=x^f$ and $g=xf$ (so that also $g=fy$), and set $Q=S_{(x,f))}$ (so 
that also $Q=S_{(f,y)}$). Thus $Q\leq S_f\cap S_g$. Also, 4.2(a) yields $Q=S_f\cap S_x$.  Set 
$$ 
P_0=N_{S_f}(Q),\ \ P_1=N_{S_g}(Q),\ \ P=\<P_0,P_1\>, 
$$ 
and set $R=P_0\cap P_1$. Then $Q\leq R$. In fact, 4.3(b) shows that $y=f\i g$ and that 
$(R^f)^y=R^g$, so $R\leq Q$, and thus $P_0\cap P_1=Q$. Assume now that $(x,f)$ is a 
counter-example to the lemma. That is, assume $Q<S_g$ (proper inclusion). Then 
$Q<P_1$ by (L2), and so $P_1< P_0$. Thus: 
\roster 

\item "{(1)}" $P_1<S_f$. 

\endroster 
Among all counter-examples, take $(x,f)$ so that $dim(Q)$ is as large as possible. We consider two 
cases, as follows. 
\vskip .1in 
\noindent 
CASE 1: $x\in N_{\ca L}(T)$. 
\vskip .1in 
As $f\in N_{\ca L}(T)$ by 4.8, we have $T\leq Q$, and then $x\in N_{\ca L}(Q)$ by 4.1(b). 
Thus $Q^g=Q^{xf}=Q^f$. Set $Q'=Q^g$. Then 2.2(b) yields an isomorphism 
$c_f:N_{\ca L}(Q)\to N_{\ca L}(Q')$. Here $f=x\i g$ so $c_f=c_{x\i}\circ c_g$ by 2.2(c). 
As $x\in N_{\ca N}(Q)\norm N_{\ca L}(Q)$, we obtain $(P_1)^{x\i}\leq N_{\ca N}(Q)P_1$, and then 
$$ 
(P_1)^f=((P_1)^{x\i})^g\leq (N_{\ca N}(Q)P_1)^g\leq N_{\ca N}(Q')N_S(Q'). 
$$ 
Also $(P_0)^f\leq N_S(Q')$, so 
\roster 

\item "{(2)}" $P^f\leq N_{\ca N}(Q')N_S(Q')$. 

\endroster  
Since $T\leq Q'$, $T$ is a Sylow $p$-subgroup of $N_{\ca N}(Q')$ by 4.1(c), and thus $N_S(Q')$ is a Sylow 
$p$-subgroup of $N_{\ca N}(Q')N_S(Q')$. By (2) there is then an element $v\in N_{\ca N}(Q')$ such that 
$P^{fv}\leq N_S(Q')$. In particular, we have: 
\roster 

\item "{(3)}" $P_0\leq S_{(f,v)}$ and $P\leq S_{fv}$. 

\endroster 
Set $u=v^{f\i}$. Then $(u,f)\in\bold D$, and we then have $uf=fv$ and $S_{(u,f)}=S_{(f,v)}$ by 4.2. If  
$S_{(f,v)}=S_{fv}$ then $(f,v)\in\bold D$ via $P$, so that $P\leq S_f$, contrary to (1). Thus 
$S_{(f,v)}\neq S_{fv}$, and so $(u,f)$ is a counter-example to the lemma. Then (3) and the maximality of 
$dim(Q)$ in the choice of $(x,f)$ yields $Q=P_0=N_{S_f}(Q)$, and so $Q=S_f$. As 
$(f,Q)\up(g,P)$ via $(x,\1)$, we have a contradiction to the $\up$-maximality of $f$. 
\vskip .1in 
\noindent 
CASE 2: The case $x\notin N_{\ca L}(T)$. 
\vskip .1in 
Let $h$ be $\up$-maximal, with $(g,S_g)\up(h,S_h)$. Then $T\leq S_h$ by 4.9, and there exists  
$r\in\ca N$ with $g=rh$ by 4.8. 

Set $w=(f\i,x\i,r,h)$, observe that $w\in\bold D$ via $Q^g$, and find 
$$ 
\Pi(w)=(f\i x\i)(rh)=g\i g=\1. 
$$ 
Then 2.3 yields $h=r\i xf$. Since both $f$ and $h$ are in 
$N_{\ca L}(T)$, 2.2(c) yields $r\i x\in N_{\ca L}(T)$, and so $r\i x\in N_{\ca N}(T)$. Then  
Case 1 applies to $(r\i x,f)$, and thus $S_h=S_{(r\i x,f)}\leq S_f$ (using 4.2). By definition of 
$\up$ there exist $a,b\in\ca N$ such that one has the usual sort of commutative diagram: 
$$ 
\CD 
S_h@>h>>S_{h\i} \\ 
@AaAA   @AAbA  \\ 
S_g@>>g>S_{g\i} 
\endCD. 
$$ 
As $T\leq S_h$, 4.1(b) yields 
$$ 
S_g\leq S_gT=(S_g)^aT\leq S_h, 
$$ 
and so $S_g\leq S_f$. This again contradicts (1), and completes the proof. 
\qed
\enddemo

The splitting lemma yields a useful criterion for partial normality, as follows.

\proclaim {Corollary 4.12} Let $\ca L$ be a locality, let $\ca N\norm\ca L$, and let 
$\ca K\norm\ca N$ be a partial normal subgroup of 
$\ca N$. Suppose that $\ca K$ is $N_{\ca L}(T)$-invariant. $($I.e. 
$x^h\in\ca K$ for all $(h\i,x,h)\in\bold D$ such that $x\in\ca K$ and 
$h\in N_{\ca L}(T).)$ Then $\ca K\norm\ca L$. 
\endproclaim 

\demo {Proof} Let $x\in\ca K$ and let $f\in\ca L$ such that $x^f$ is defined. 
By the Frattini Lemma, we may write $f=yg$ with $y\in\ca N$ and with $g$ $\up$-maximal, 
and then the splitting lemma yields $S_f=S_{(y,g)}$. Set $u=(f\i,x,f)$ and 
$v=(g\i,y\i,x,y,g)$. Then $S_u=S_v\in\D$, and $x^f=\Pi(u)=\Pi(v)=(x^y)^g$. Thus 
$x^f\in\ca K$, and $\ca K\norm\ca L$. 
\qed 
\enddemo

A subset $X$ of $\ca L$ of the form $\ca N f$, $f\in\ca L$, will be called a {\it coset} of $\ca N$. 
A coset $\ca Nf$ is {\it maximal} if it is not a proper subset of any coset of $\ca N$.

\proclaim {Proposition 4.13} Let $\ca N$ be a partial normal subgroup of the locality $\ca L$. 
\roster 

\item "{(a)}" $\ca Nf=f\ca N$ for all $f\in N_{\ca L}(T)$, and if $f$ is $\up$-maximal with respect to 
$\ca N$ then $\ca Nf=\ca N f\ca N=f\ca N$. 

\item "{(b)}" Let $f,g\in\ca L$. Then 
$$ 
(g,S_g)\up(f,S_f)\ \iff\ \ca Ng\sub\ca Nf\ \iff\ g\in\ca Nf. 
$$ 

\item "{(c)}" $g\in\ca L$ is $\up$-maximal relative to $\ca N$ if and only if $\ca Ng$ is a maximal coset 
of $\ca N$. 

\item "{(d)}" $\ca L$ is partitioned by the set $\ca L/\ca N$ of maximal cosets of $\ca N$. 

\item "{(e)}" Let $u:=(g_1,\cdots,g_n)\in\bold D$ and let $v:=(f_1,\cdots,f_n)$ be a sequence of 
$\up$-maximal elements of $\ca L$ such that $g_i\in\ca Nf_i$ for all $i$. Then $v\in\bold D$, 
$TS_u\leq S_v$, and $\ca N\Pi(u)\sub\ca N\Pi(v)$.  

\endroster 
\endproclaim 

\demo {Proof} (a): That $\ca Nf=f\ca N$ for $f\in N_{\ca L}(T)$ is given by 4.2. Now let $f$ be 
$\up$-maximal relative to $\ca N$. Then $f\in N_{\ca L}(T)$ by 4.8. Let $x,y\in\ca N$ such that 
$(x,f,y)\in\bold D$. Then $(x,fy)\in\bold D$, and $(x,fy)=(x,y'f)$ where $y'=fyf\i$. The splitting 
lemma (4.11) then yields $(x,y',f)\in\bold D$, and thus $\ca Nf\ca N\sub\ca Nf$. The reverse inclusion 
is obvious, and yields (a). 

\vskip .1in 
\noindent 
(b): Let $f,g\in\ca L$ and suppose that $(g,S_g)\up(f,S_f)$. Then $g=xf$ for some $x\in\ca N$ by 4.10, 
and then 4.11 yields $S_g=S_{(x,f)}\leq S_f$. Let $y\in\ca N$ such that $(y,g)\in\bold D$. 
Since $(S_{(y,g)})^y\leq S_g=S_{(x,f)}$ we get $(y,x,f)\in\bold D$, and $yg=(yx)f\in\ca Nf$. Thus 
$(g,S_g)\up(f,S_f)\implies \ca Ng\sub\ca Nf$. Clearly $\ca Ng\sub\ca Nf\implies g\in\ca Nf$. The 
required circle of implications is then completed by 4.11. 
\vskip .1in 
\noindent 
(c): Immediate from (b). 
\vskip .1in 
\noindent 
(d): Let $f$ and $g$ be $\up$-maximal, and let $h\in\ca Nf\cap\ca Ng$. Thus there exist $x,y\in\ca N$ 
with $(x,f)\in\bold D$, $(y,g)\in\bold D$, and with $h=xf=yg$. Then also $(x\i,x,f)\in\bold D$, so 
$(x\i,yg)=(x\i,xf)\in\bold D$. The splitting lemma then yields $(x\i,y,g)\in\bold D$, and 
we thereby obtain $f=x\i yg\in\ca Ng$. Now (b) implies that $\ca Nf\sub\ca Ng$, and symmetry gives 
the reverse inclusion. Thus $\ca Nf=\ca Ng$ if $\ca Nf\cap\ca Ng\neq\nset$. 
\vskip .1in
\noindent
(e): Let $x_i\in\ca N$ with $g_i=x_if_i$. Set $Q_0=S_u$, and set 
$Q_i=(Q_{i-1})^{g_i}$ for $1\leq i\leq n$. Then 
$Q_{i-1}\leq S_{g_i}=S_{(x_i,f_i)}\leq S_{f_i}$ by 4.11, and then also 
$Q_{i-1}T=(Q_{i-1})^{x_i}T\leq S_{f_i}$, by 4.1(b) and 4.8. Thus: 
$$
(Q_{i-1}T)^{f_i}=(Q_{i-1})^{x_,f_i}T=Q_iT\leq S_{f_{i+1}}, 
$$
and $v\in\bold D$ via $Q_0T$. Now $\Pi(u)\in\ca N\Pi(v)$ by the Frattini calculus (4.4), and 
then (b) implies that $\ca NPi(u)\sub\ca N\Pi(v)$. 
\qed
\enddemo

Let $\equiv$ be the equivalence relation on $\ca L$ defined by the partition 
in 4.13(d). In view of points (a) and (b) of 4.13 we refer 
to the $\equiv$-classes as the {\it maximal cosets} of $\ca N$ in $\ca L$.

\proclaim {Lemma 4.14} Let $\ca H$ be a partial subgroup of the locality $\ca L$, containing the partial 
normal subgroup $\ca N\norm\ca L$. Then $\ca H$ is the disjoint union of the maximal cosets of $\ca N$ 
contained in $\ca H$. 
\endproclaim 

\demo {Proof} Let $f\in\ca H$. Apply the Frattini lemma (4.10) to obtain 
$f=xh$ for some $x\in\ca N$ and some $h\in N_{\ca L}(T)$ such that $h$ is 
$\up$-maximal with respect to $\ca N$. Then $h=x\i f$ by 1.3(d), and thus 
$h\in\ca H$ as $\ca N\leq\ca H$. Then also $\ca Nh\sub\ca H$, where 
$\ca Nh$ is a maximal coset of $\ca N$ by 4.13(b). 
\qed 
\enddemo

The set $\ca L/\ca N$ of maximal cosets of $\ca N$ may also be denoted $\bar{\ca L}$. Let 
$\r:\ca L\to\bar{\ca L}$ be the mapping which sends $g\in\ca L$ to the unique maximal coset of $\ca N$ 
containing $g$. Set $\bold W:=\bold W(\ca L)$ and $\bar{\bold W}=\bold W(\bar{\ca L})$, and let 
$\r^*:\bold W\to\bar{\bold W}$ be the induced mapping of free monoids. For any subset or element $X$ of 
$\bold W$, write $\bar X$ for the image of $X$ under $\r^*$, and similarly if $Y$ is a subset or element 
of $\ca L$ write $\bar Y$ for the image of $Y$ under $\r$. In particular, $\bar{\bold D}$ is the image of 
$\bold D$ under $\r^*$. 

For $w\in\bold W$, we shall say that $w$ is $\up$-maximal if every entry of $w$ is $\up$-maximal.

\proclaim {Lemma 4.15} There is a unique mapping $\bar\Pi:\bar{\bold D}\to\bar{\ca L}$, a unique 
involutory bijection $\bar g\maps{\bar g}\i$ on $\bar{\ca L}$, and a unique element $\bar{\1}$ 
of $\bar{\ca L}$ such that $\bar{\ca L}$, with these structures, is a partial group, and such that 
$\r$ is a homomorphism of partial groups. Moreover, we have $Ker(\r)=\ca N$, and the homomorphism of free 
monoids $\r^*:\bold W(\ca L)\to\bold W(\bar{\ca L})$ maps $\bold D(\ca L)$ onto $\bold D(\bar{\ca L})$. 
\endproclaim 

\demo {Proof} Let $u=(g_1,\cdots,g_n)$ and $v=(h_1,\cdots,h_n)$ be members of $\bold D$ such that 
$\bar u=\bar v$. By 4.13(d) there exists, for each $i$, an $\up$-maximal $f_i\in\ca L$ with 
$g_i,h_i\in\ca Nf_i$. Set $w=(f_1,\cdots,f_n)$. Then $w\in\bold D$ by 
4.13(e), and then 4.3(a) shows that $\Pi(u)$ and $\Pi(v)$ are elements of 
$\ca N\Pi(w)$. Thus $\bar{\Pi(u)}=\bar{\Pi(w)}=\bar{\Pi(v)}$, and there 
is a well-defined mapping $\bar\Pi:\bar{\bold D}\to\bar{\ca L}$ given by 
$$ 
\bar\Pi(w)=\bar{\Pi(w)}.\tag* 
$$ 
For any subset $X$ of $\ca L$ write $X\i$ for the set of inverses of elements 
of $X$. For any $f\in\ca L$ we then have $(\ca Nf)\i=f\i\ca N\i$ by 
1.1(4). Here $\ca N\i=\ca N$ as $\ca N$ is a partial group, and then 
$(\ca Nf)\i=\ca Nf\i$ by 4.13(a). The inversion map $\ca Nf\maps\ca Nf\i$ is 
then well-defined, and is an involutory bijection on $\bar{\ca L}$. Set $\bar{\1}=\ca N$. 

We now check that the axioms in 1.1, for a partial group, are satisfied by 
the above structures. Since $\bar{\bold D}$ is the image of $\bold D$ under 
$\r^*$, we get $\bar{\ca L}\sub\bar{\bold D}$. Now let 
$\bar w=\bar u\circ\bar v\in\bar{\bold D}$, let $u,v$ be $\up$-maximal 
pre-images in $\bold W$ of $\bar u$, and $\bar v$, and set $w=u\circ v$. Then 
$w$ is $\up$-maximal, and so $w\in\bold D$ by 4.13(e). Then $u$ and $v$ are in $\bold D$, 
and so $\bar u$ and $\bar v$ are in $\bar{\bold D}$. Thus 
$\bar{\bold D}$ satisfies 1.1(1). Clearly, (*) implies that $\bar\Pi$ 
restricts to the identity on $\bar{\ca L}$, so $\bar\Pi$ satisfies 1.1(2). 

Next, let $\bar u\circ\bar v\circ\bar w\in\bar{\bold D}$, and choose corresponding 
$\up$-maximal pre-images $u,v,w$. Set $g=\Pi(v)$. Then 
$\bar g=\bar\Pi(\bar v)$ by (*). By 1.1(3) we have both 
$u\circ v\circ w$ and $u\circ(g)\circ w$ in $\bold D$, and these two words 
have the same image under $\Pi$. Applying $\r^*$ we obtain words in 
$\bar{\bold D}$ having the same image under $\bar\Pi$, and thus $\bar\Pi$ 
satisfies 1.1(3). By definition, $\bar\Pi(\nset)=\bar{\1}$, and then the 
condition 1.1(4) is readily verified. Thus, $\bar{\ca L}$ is a partial 
group. 

By definition, $\bar{\bold D}$ is the image of $\bold D$ under $\r^*$. So, 
in order to check that $\r$ is a homomorphism of partial groups it suffices 
to show that if $w\in\bold D$ then $\bar\Pi(w\r^*)=\Pi(w)\r$. But this is simply the statement (*). 
Moreover, it is this observation which establishes that the given partial 
group structure on $\bar{\ca L}$ is the unique one for which $\r$ is a 
homomorphism of partial groups. We have $f\in Ker(\r)$ if and only if 
$f\g=\bar{\1}=\ca N$. Since $\ca Nf\sub\ca N$ implies $f\in\ca N$, 
and since $\ca N$ is the maximal coset of $\ca L$ containing $\1$, we 
obtain $Ker(\r)=\ca N$. 
\qed 
\enddemo

\vskip .2in 
\noindent 
{\bf Section 5: Quotient localities} 
\vskip .1in 

We continue the setup of the preceding section, where $(\ca L,\D,S)$ is a fixed locality and 
$\ca N\norm\ca L$ is a partial normal subgroup. We have seen in 4.15 that the set $\ca L/\ca N$ of 
maximal cosets of $\ca N$ inherits from $\ca L$ a partial group structure via the projection map 
$\r:\ca L\to\ca L/\ca N$. The aim now is to go further, and to show that $\ca L/\ca N$ is a locality.

\proclaim {Lemma 5.1} Let $(\ca L,\D,S)$ be a locality and let $\ca N\norm\ca L$. Then 
$(\ca NS,\D,S)$ is a locality. 
\endproclaim 

\demo {Proof} It follows from 2.17 that $\ca NS$ is a partial subgroup of $\ca L$. One observes that 
$\bold D(\ca NS)$ is the subset $\bold D_\D$ of $\bold W(\ca NS)$, as defined in 2.1, and this suffices 
to show that $(\ca NS,\D)$ is objective. Thus $(\ca NS,\D,S)$ is a pre-locality, and is finite-dimensional 
by 2.16(b). The conditions (L1) and (L2) in definition 3.1 are inherited by $\ca NS$ from $\ca L$ in 
an obvious way, so $(\ca NS,\D,S)$ is a locality.  
\qed 
\enddemo 

\proclaim {Lemma 5.2} Let $P\leq S$ be a subgroup of $S$, and let $X$ be a subset of $S$ containing $P$ 
and having the property that $P^x\sub X$ for all $x\in X$. Then either $X=P$ or $P$ is a proper 
subset of $X\cap N_S(P^\star)$. 
\endproclaim 

\demo {Proof} Let $\{S_n\}$ be a framing of $S$ by finite subgroups, and set $P_n=P\cap S_n$ and 
$X_n=X\cap S_n$. Then $(P_n)^x\sub X_n$ for all $x\in X_n$. By finite dimensionality, 
$P^\star=(P_n)^\star$ for all sufficiently large $n$, so we may start the indexing so that 
$P^\star=(P_0)^\star$. 
Notice that $X_n\sub X_{n+1}$ for all $n$ and that $X=\bigcup\{X_n\}$. Assuming that $P\neq X$, there is 
then an index $n$ such that $P_n\neq X_n$. Again, we may assume that $n=0$. 

Let $\ca C$ be the set of subgroups $Q$ of $S_0$ such that $X_0\cap Q=P_0$. Then $P_0\in\ca C$, 
and since $S_0$ is finite there exists a maximal $Q\in\ca C$ with respect to inclusion. As $P_0\neq X_0$  
we have $Q\neq S_0$, and so $Q<N_{S_0}(Q)$. Thus $N_{S_0}(Q)\notin\ca C$, and so there exists 
$x\in X_0\cap N_{S_0}(Q)$ with $x\notin P_0$. Then $(P_0)^x\sub X_0\cap Q=P_0$, and so 
$(P_0)^x=P_0$. Then $x\in N_S(P^\star)$ by 2.13, and $x\notin P$ since $P\cap X_0\sub P\cap S_0=P_0$. 
\qed 
\enddemo

\proclaim {Theorem 5.3} Let $(\ca L,\D,S)$ be a locality, let $\bar{\ca L}$ be a partial group, and 
let $\b:\ca L\to\bar{\ca L}$ be a homomorphism of partial groups such that the induced map 
$\b^*:\bold W(\ca L)\to\bold W(\bar{\ca L})$ of free monoids restricts to a surjection 
$\bold D(\ca L)\to\bold D(\bar{\ca L})$. 
Set $\ca N=Ker(\b)$ and $T=S\cap\ca N$. Further, set $\bold D=\bold D(\ca L)$, 
$\bar{\bold D}=\bold D(\bar{\ca L})$, $\bar S=S\b$, and $\bar\D=\{P\b\mid P\in\D\}$. Then 
$(\bar{\ca L},\bar\D,\bar S)$ is a locality. Moreover: 
\roster 

\item "{(a)}" The fibers of $\b$ are the maximal cosets of $\ca N$. 

\item "{(b)}" For each $\bar w\in\bold W(\bar{\ca L})$ there exists $w\in\bold W(\ca L)$ such that 
$\bar w=w\b^*$ and such that each entry of $w$ is $\up$-maximal relative to the partial normal subgroup 
$\ca N$ of $\ca L$. For any such $w$ we then have $\bar{S_w}=\bar S_{\bar w}$, and 
$w\in\bold D$ if and only if $\bar w\in\bar{\bold D}$.  

\item "{(c)}" Let $P,Q\in\D$ with $T\leq P\cap Q$. Then $\b$ restricts to a surjection 
$N_{\ca L}(P,Q)\to N_{\bar{\ca L}}(P\b,Q\b)$; and to a surjective homomorphism if $P=Q$.  

\item "{(d)}" $\b$ is an isomorphism if and only if $\ca N=1$.  

\endroster 
\endproclaim 

\demo {Proof} The hypothesis that $\bold D\b^*=\bar{\bold D}$ implies that $\b^*$ maps the set of words of 
length 1 in $\ca L$ onto the set of words of length 1 in $\bar{\ca L}$. Thus $\b$ is surjective.  

Let $M$ be a subgroup of $\ca L$. The restriction of $\b$ to $M$ is then a homomorphism of partial groups, 
and hence a homomorphism $M\to M\b$ of groups by 1.13. In particular, $\bar S$ is a $p$-group, and $\bar\D$ 
is a set of subgroups of $\bar S$. 

We have $\ca N\norm\ca L$ by 1.14. Let $\L$ be the set of elements $g\in\ca L$ such that $g$ is 
$\up$-maximal relative to $\ca N$. For any $g\in\L$, $\b$ is constant on the maximal coset $\ca Ng$
(see 4.14) of $\ca N$, so $\b$ restricts to a surjection of $\L$ onto $\bar{\ca L}$. This shows that $\b^*$ 
restricts to a surjection of $\bold W(\L)$ onto $\bold W(\bar{\ca L})$. If $\bar w\in\bar{\bold D}$ then 
there exists $w\in\bold D$ with $w\b^*=w'$, and then 4.13(e) shows that such a $w$ may be chosen to be in 
$\bold W(\L)$. Set $\bold D(\L)=\bold D\cap\bold W(\L)$. Thus: 
\roster 

\item "{(1)}" $\b^*$ maps $\bold D(\L)$ onto $\bar{\bold D}$.

\endroster 
Let $\bar w\in\bar{\bold D}$, let $w\in\bold D(\L)$ with $w\b^*=\bar w$, let $\bar a,\bar b\in\bar S$, 
and let $a,b\in S$ with $a\b=\bar a$ and $b\b=\bar b$. Then $a,b\in\L$ by 4.6(a), and 
$(a)\circ w\circ(b)\in\bold D(\L)$ by 2.17. Then $(\bar a)\circ \bar w\circ(\bar b)\in\bar{\bold D}$. 
This shows: 
\roster 

\item "{(2)}" $\bar{\bold D}$ is a $\bar S$-biset (cf. 2.7). 

\endroster 
Let $\bar a\in\bar S$, let $a\in S$ be a preimage of $\bar a$, and let $h\in\ca L$ be any preimage of 
$\bar a$. Then $(a,h\i)\in\bold D$, and $ah\i\in\ca N$, so: 
\roster 

\item "{(3)}" The $\b$-preimage of an element $\bar a\in\bar S$ is a maximal coset $\ca Na$, where $a\in S$. 

\endroster 
Fix $\bar g\in\bar{\ca L}$, let $g\in\L$ with $g\b=\bar g$, set $P=S_g$, and set 
$$ 
\bar S_{\bar g}=\{\bar x\in\bar S\mid \bar x^{\bar S}\in \bar S\}.   
$$ 
Let $\bar a\in\bar S_{\bar g}$ and set $\bar b=\bar a^{\bar g}$. As in the proof of 2.7, we show 
that $\bar P^{\bar a}\sub\bar S_{\bar g}$. Namely, from $(\bar g\i,\bar a,\bar g)\in\bar{\bold D}$ 
and $\bar Pi(\bar g\i,\bar a,\bar g)=\bar b$ we obtain (from two applications of (2)):
$$ 
(\bar P^{\bar a})^{\bar g}=\bar P^{\bar b}\leq\bar S. \tag*
$$ 
Thus, the set $\bar\G$ of all $\bar S_{\bar g}$-conjugates of $\bar P$ is a set of subgroups of  
$\bar S_{\bar g}$. Let $\G$ be the set of all preimages in $S$ of members of $\bar\G$, and set 
$X=\bigcup\G$. Then $X=X\i$ is a subset of $S$, containing $P$, and having the property that 
$P^x\sub X$ for all $x\in X$. As $P=S_g\in\Omega$, 5.2 implies 
that either $X=P$ or $P$ is properly contained in $N_S(P)\cap S$. 

Assuming now that $\bar P\neq\bar S_{\bar g}$, we may choose an element $x\in N_X(P)-P$. 
Set $Q=P\<x\>$. Then $Q^g$ is defined (and is a subgroup of $N_{\ca L}(P^g)$) by 2.3(b). As 
$\bar Q^{\bar g}\leq\bar S$ we obtain $Q^g\leq\ca NS$. As $\ca NS$ is a locality by 5.1, $S$ is a 
Sylow $p$-subgroup of $\ca NS$ by 3.8. Thus there exists $f\in\ca N$ 
with $(Q^g)^f\leq S$. Here $(g,f)\in\bold D$ via $P$, so $Q\leq S_{gf}$. This contradicts the 
$\up$-maximality of $g$, so we conclude: 
\roster 

\item "{(4)}" $(S_g)\b=\bar S_{g\b}$ for each $g\in\L$. 

\endroster 

For $\bar w\in\bold W(\bar{\ca L})$ define $\bar S_{\bar w}$ to be the set of all $\bar x\in\bar S$ such that 
$\bar x$ is conjugated successively into $\bar S$ by the entries of $\bar w$. As an immediate consequence of 
(4): 
\roster 

\item "{(5)}" $(S_w)\b=\bar S_{w\b^*}$ for all $w\in\bold W(\L)$. 

\endroster 
Notice that (5) implies  point (b). 

We now verify that $(\bar{\ca L},\bar\D)$ is objective. Thus, let $\bar w\in\bold W(\bar{\ca L})$ with 
$\bar S_{\bar w}\in\bar\D$. Let $w\in\bold W(\L)$ with $w\b^*=w'$. Then (5) yields $S_w\in\D$, so 
$w\in\bold D$, and hence $\bar w\in\bar{\bold D}$. Thus $(\bar{\ca L},\bar\D)$ satisfies the condition (O1) 
in definition 2.1 of objectivity. Now let $\bar P\in\bar D$, let $\bar f\in N_{\bar{\ca L}}(\bar P,\bar S)$, 
and set $\bar Q=(\bar P)^{\bar f}$. Then (4) yields $\bar Q=Q\b$ for some $Q\in\D$, and thus 
$\bar Q\in\bar\D$. Any overgroup of $\bar Q$ in $\bar S$ is the image of an overgroup 
of $Q$ in $S$ as $\b$ maps $S$ onto $\bar S$, so $(\bar{\ca L},\bar\D)$ satisfies (O2) in 2.1. Thus 
$(\bar{\ca L},\bar\D)$ is objective. As $\bar{\D}$ is a set of subgroups of $\bar S$, 
$(\bar{\ca L},\bar\D,\bar S)$ is then a pre-locality, and  
finite-dimensional since, by (5), the poset $\Omega_{\bar S}(\bar{\ca L})$ is a 
homomorphic image of $\Omega_S(\ca L)$. 

Let $P,Q\in\D$ with $T\leq P\cap Q$, set $\bar P=P\b$ and $\bar Q=Q\b$, and let $\bar g\in\bar{\ca L}$ 
such that $\bar P^{\bar g}$ is defined and is a subset of $\bar Q$. Let $g$ be a preimage of $\bar g$ 
in $\L$. Then $P\leq S_g$ by (4). As $\b$ is a homomorphism, $(P^g)\b$ is a subgroup of $\bar Q$. 
Then $P^g\leq Q$ since $\b$ restricts to an epimorphism $S\to\bar S$ with kernel $T$. As $\b$ maps 
subgroups of $\ca L$ homomorphically to subgroups of $\bar{\ca L}$ (by 1.13) we obtain (c). 
In the case that $P=Q=S$ we obtain in this way an epimorphism from $N_{\ca L}(S)$ to 
$N_{\bar{\ca L}}(\bar S)$. As $S$ is a Sylow $p$-subgroup of $N_{\ca L}(S)$ it follows 
that $\bar S$ is a maximal $p$-subgroup of $\bar{\ca L}$. Thus 
$\bar{\ca L}$ satisfies the condition (L1) in definition 3.1. The condition (L2) is inherited by 
$\bar{\ca L}$ in an obvious way, and hence $\bar{\ca L}$ is a locality.  

Let $g,h\in\L$ with $\g\b=h\b$. Then $S_g=S_h$ by (4), so $(g,h\i)\in\bold D$, and $(gh\i)\b=\1$. Thus 
$g\in\ca Nh$, and then $\ca Ng=\ca Nh$ by 4.14(b). This yields (a), and it remains only to prove (d). 

If $\b$ is an isomorphism then $\b$ is injective, and $\ca N=\1$. On the other hand, suppose that 
$\ca N=\1$. Then (a) shows that $\b$ is injective, and so $\b$ is a bijection. That
$\b\i$ is then a homomorphism of partial groups is given by (5). Thus (d) holds, and the proof is 
complete.  
\qed 
\enddemo

\definition {Definition 5.4} Let $(\ca L,\D,S)$ and $(\ca L',\D',S')$ be localities, and let 
$\b:\ca L\to\ca L'$ be a homomorphism of partial groups. Then $\b$ is a {\it projection} if: 
\roster 

\item "{(1)}" $\bold D\b^*=\bold D'$, and  

\item "{(2)}" $\D'=\{P\b\mid P\in\D\}$. 

\endroster 
\enddefinition 

\proclaim {Corollary 5.5} Let $(\ca L,\D,S)$ be a locality, let $\ca N\norm\ca L$ be a partial normal 
subgroup and let $\r:\ca L\to\ca L/\ca N$ be the mapping which sends $g\in\ca L$ to the unique maximal 
coset of $\ca N$ containing $g$. Set $\bar{\ca L}=\ca L/\ca N$, set $\bar S=S\r$, and let $\bar\D$ be 
the set of images under $\r$ of the members of $\D$. Regard $\bar{\ca L}$ as a partial group in the 
unique way (given by 4.15) which makes $\r$ into a homomorphism of partial groups. Then 
$(\bar{\ca L},\bar\D,\bar S)$ is a locality, and $\r$ is a projection. 
\endproclaim 

\demo {Proof} Immediate from 4.15 and 5.3. 
\qed 
\enddemo 

\proclaim {Theorem 5.6 (``First Isomorphism Theorem")} Let $(\ca L,\D,S)$ and $(\ca L',\D',S')$ be 
localities, let $\b:\ca L\to\ca L'$ be a projection, and let $\ca N\norm\ca L$ be a partial normal 
subgroup of $\ca L$ contained in $Ker(\b)$. Let $\r:\ca L\to\ca L/\ca N$ be the 
projection given by 5.3. Then there exists a unique homomorphism 
$$ 
\g:\ca L/\ca N\to\ca L' 
$$ 
such that $\r\circ\g=\b$; and $\g$ is then a projection. Moreover, $\g$ is an isomorphism if and only if 
$\ca N=Ker(\b)$. 
\endproclaim  

\demo {Proof} Set $\ca M=Ker(\b)$, and let $h\in\ca L$ be $\up$-maximal relative to $\ca M$. Then 
$\ca Mh$ is a maximal coset of $\ca M$ in $\ca L$ by 4.13(b), and $\ca Mh=\ca Mh\ca M$ by 4.13(a). Let 
$g\in\ca Mh$. As $\ca N\leq\ca M$ by hypothesis, the splitting lemma (4.11, as applied to $\ca M$ and $h$) 
yields 
$$ 
\ca Ng\ca N\sub\ca N(\ca Mh)\ca N\sub\ca M(\ca Mh)\ca M=\ca Mh\ca M. 
$$  
The definition 4.6 of the relation $\up$ on $\ca L\circ\D$ shows that $\ca Ng\ca N$ contains an element 
$f$ which is $\up$-maximal with respect to $\ca N$. Then $f\in\ca Mh\ca M$, so $f\in\ca Mh$, and another 
application of the splitting lemma yields $\ca Nf\sub\ca Mh$. Here $\ca Nf$ is the maximal coset of $\ca N$ 
containing $f$. We have thus shown: 
\roster 

\item "{(*)}" The partition $\ca L/\ca N$ of $\ca L$ is a refinement of the partition $\ca L/\ca M$. 

\endroster 
By 5.3(a) $\b$ induces a bijection $\ca L/\ca M\to\ca L'$.  Set $\bar{\ca L}=\ca L/\ca N$. Then (*) 
implies that there is a mapping $\g:\bar{\ca L}\to\ca L'$ which sends the maximal coset $\ca Nf\ca N$ 
to $f\b$. Clearly, $\g$ is the unique mapping $\bar{\ca L}\to\ca L'$ such that $\r\circ\g=\b$. 

Let $\bar w\in\bold D(\bar{\ca L})$. Then 5.3(b) yields a word $w\in\bold D$ such that $w\r^*=\bar w$ and 
such that the entries of $w$ are $\up$-maximal relative to $\ca N$. We have $\bar w\g^*=w\b^*$, so 
$\g^*$ maps $\bold D(\bar{\ca L})$ into $\bold D(\ca L')$. Let $\Pi'$ and $\bar\Pi$ be the products in 
$\ca L'$ and $\bar{\ca L}$, respectively. As $\b$ and $\r$ are homomorphisms we get 
$$ 
\Pi'(\bar w\g^*)=\Pi'(w\b^*)=(\Pi(w))\b=(\bar{\Pi(w)})\g=(\bar\Pi(\bar w))\g,  
$$ 
and thus $\g$ is a homomorphism. As $\b=\r\circ\g$ is a projection, and $\r$ is a projection, one 
verifies that $\g^*$ maps $\bold D(\bar{\ca L})$ onto $\bold D(\ca L')$ and that $\g$ maps 
$\bar\D$ onto $\D'$. Thus $\g$ is a projection. 

We have $\ca M=\ca N$ if and only if $Ker(\g)=\1$. Then 5.3(d) shows that $\g$ is an isomorphism if 
and only if $\ca M=\ca N$; completing the proof. 
\qed 
\enddemo

\proclaim {Proposition 5.7 (Partial Subgroup Correspondence)} Let $(\ca L,\D,S)$ and 
$(\bar{\ca L},\bar\D,\bar S)$ be localities, and let $\b:\ca L\to\bar{\ca L}$ be a projection. Set 
$\ca N=Ker(\b)$ and set $T=S\cap\ca N$. Then $\b$ induces a bijection $\s$ from the set $\frak H$ of 
partial subgroups $\ca H$ of $\ca L$ containing $\ca N$ to the set $\bar{\frak H}$ of partial subgroups 
$\bar{\ca H}$ of $\bar{\ca L}$. Moreover, for any $\ca H\in\frak H$, we have $\ca H\b\norm\ca L'$ if 
and only if $\ca H\norm\ca L$. 
\endproclaim 

\demo {Proof} Any partial subgroup of $\ca L$ containing $\ca N$ is a union of maximal cosets of $\ca N$ by 
4.14. Then 5.3(a) enables the same argument that one has for groups, for proving that $\r$ induces a 
bijection $\frak H\to\bar{\frak H}$. Since each maximal coset of $\ca N$ contains an element which is 
$\up$-maximal with respect to $\ca N$, one may apply 5.3(b) in order to show that a partial subgroup 
$\ca H\in\frak H$ is normal in $\ca L$ if and only if its image is normal in $\bar{\ca L}$. 
\qed 
\enddemo

\definition {5.8 Remark}
A comprehensive ``second isomorphism theorem" appears to be out of reach, for two reasons. First, given a 
partial subgroup $\ca H\leq\ca L$ and a partial normal subgroup $\ca N\norm\ca L$, there appears to be no 
reason why the image of $\ca H$ under the projection $\r:\ca L\to\ca L/\ca N$ should be a partial subgroup of 
$\ca L/\ca N$, other than in special cases. Second, there seems to be no way, in general, to define the 
quotient of $\ca H$ over the partial normal subgroup $\ca H\cap\ca N$ of $\ca H$. On the other hand, a
``third isomorphism theorem" may easily be deduced from 4.4 
and from the observation that a composition of projections is again a projection. 
\enddefinition

\proclaim {Lemma 5.9} Let $\ca N\norm\ca L$ and let $\r:\ca L\to\ca L/\ca N$ be the canonical projection. 
Further, let $\ca H$ be a partial subgroup of $\ca L$ containing $\ca N$ and let $X$ be an 
arbitrary subset of $\ca L$. Then $(X\cap\ca H)\r=X\r\cap\ca H\r$. 
\endproclaim 

\demo {Proof} By 4.14, $\ca H$ is a union of maximal cosets of $\ca N$, and then $\ca H\r$ is the set of 
those maximal cosets. On the other hand $X\r$ is the set of all maximal cosets $\ca Ng$ of $\ca N$ such 
that $X\cap\ca Ng\neq\nset$. Thus $X\r\cap\ca H\r\sub(X\cap\ca H)\r$. The reverse inclusion is obvious. 
\qed 
\enddemo 

\proclaim {Corollary 5.10} Let $\ca N\norm\ca L$, and let $\ca M$ be a 
partial normal subgroup of $\ca L$ containing $\ca N$. Let $\r:\ca L\to\ca L/\ca N$ be the canonical 
projection. Then $(S\cap\ca M)\r$ is a maximal $p$-subgroup of $\ca M\r$. 
\endproclaim 

\demo {Proof} Write $(\bar{\ca L},\bar\D,\bar S)$ for the quotient locality given by 5.5, and set 
$\bar{\ca M}=\ca M\r$. Applying 5.9 with $S$ in the role of $X$, we obtain 
$(S\cap\ca M)\r=\bar S\cap\bar{\ca M}$. Since $\bar{\ca M}\norm\bar{\ca L}$, it follows from 4.1(c) that 
$\bar S\cap\bar{\ca M}$ is maximal in the poset of $p$-subgroups of $\bar{\ca M}$, completing the proof. 
\qed 
\enddemo

\proclaim {Proposition 5.11} Let $\ca N\norm\ca L$, set $T=S\cap\ca N$, and set $\ca L_T=N_{\ca L}(T)$. 
Set $\bar{\ca L}=\ca L/\ca N$ and let $\r:\ca L\to\bar{\ca L}$ be the canonical projection. 
Then the partial subgroup $\ca L_T$ of $\ca L$ is a locality $(\ca L_T,\D,S)$, and the  
restriction of $\r$ to $\ca L_T$ is a projection $\ca L_T\to\bar{\ca L}$. 
\endproclaim 

\demo {Proof} That $\ca L_T$ is a partial subgroup of $\ca L$ having the structure of a locality 
$(\ca L_T,\D,S)$ is given by 3.14.
Let $\r_T$ be the restriction of $\r$ to $\ca L_T$. Then $\r_T$ is 
a homomorphism of partial groups, and 4.14(e) shows that $\r_T$ maps $\bold D(\ca L_T)$ onto 
$\bold D(\bar{\ca L})$. That is, $\r_T$ is a projection. 
\qed 
\enddemo

\vskip .2in 
\noindent 
{\bf Section 6: An application} 
\vskip .1in 

The aim of this brief section is to show that, under suitable conditions on $\D$ and on $S$, the locality 
$(\ca L,\D,S)$ has a partial normal subgroup $\Theta$ which is the set-theoretic union of 
groups $\Theta(P)$ for $P\in\D$, and where $\Theta(P)$ is the largest normal $p'$-subgroup of $N_{\ca L}(P)$.  
We begin by establishing the existence of such a normal $p'$-subgroup, in the following lemma. 
But first: let us say that a group $G$ is {\it lim-finite} if $G$ is countable and locally finite.   

\proclaim {Lemma 6.1} Let $G$ be a lim-finite group. Then there is a largest normal subgroup $K\norm G$ 
such that $K$ contains no elements of order $p$. Moreover, if there exists a maximal $p$-subgroup 
$Z$ of $G$ with $Z\leq Z(G)$, then $G$ is the direct product $Z\times K$. 
\endproclaim 

\demo {Proof} Let $\{G_n\}_{n=1}^\infty$ be a framing of $G$ by finite subgroups. Then each $G_n$ has 
a largers normal subgroup $K_n$ of order prime to $p$, and we have $K_n\leq K_{n+1}$ for all $n$. Take  
$K$ to be the union of the groups $K_n$. Then $K\norm G$ and $K$ has no elements of order $p$. If 
$K'$ is any other such normal subgroup of $G$ then $K'\cap G_n\leq K$ for all $n$, and thus $K'\leq K$. 

Now let $Z$ be a maximal $p$-subgroup of $G$, and suppose that $Z\leq Z(G)$. For each $n$ let 
$P_n$ be a Sylow $p$-subgroup of $G_n$. Then $P_nZ$ is a $p$-group, so $P_n\leq Z$. Thus 
$P_n=Z\cap G_n$, whence $G_n=P_n\times K_n$, and then $G=Z\times K$. 
\qed 
\enddemo

\definition {Definition 6.2} Let $P$ be a $p$-group. Then $P$ has the {\it normalizer-increasing property} 
if for each pair of subgroups $U,V$ of $P$ with $U<V$, we have $U<N_V(U)$. 
\enddefinition 

\definition {Remark} Not every lim-finite $p$-group has the normalizer-increasing property. For example, 
Take $P_1$ to be a group of order $p$, and recursively define $P_{n+1}$ to be the wreath product 
$P_n\wr P_1$. There is an obvious inclusion of $P_n$ in $P_{n+1}$, and thus $\{P_n\}_{n=1}^\infty$ 
is a framing of a $p$-group $P$ by finite subgroups. Observe that $P$ contains a normal 
elementary abelian subgroup $E$ of infinite rank, and a subgroup $Q$ such that $P\cong Q$ 
and such that $P$ is the semi-direct product $E\rtimes Q$. Then $Q<P$ but $Q=N_P(Q)$. 
\enddefinition

\proclaim {Proposition 6.3} Let $(\ca L,\D,S)$ be a locality, and assume: 
\roster 

\item "{(1)}" $S$ has the normalizer-increasing property, and 

\item "{(2)}" $C_S(P)\leq P$ for all $P\in\D$. 

\endroster  
For $P\in\D$ let $\Theta(P)$ be the largest normal $p'$-subgroup of $N_{\ca L}(P)$, as given by 6.1. Set  
$$ 
\Theta=\bigcup\{\Theta(P)\}_{P\in\D}. 
$$ 
Then $\Theta\norm\ca L$, $S\cap\Theta=1$, and the canonical projection $\r:\ca L\to\ca L/\Theta$ resricts 
to an isomorphism $S\to S\r$. Moreover, upon identifying $S$ with $S\r$: 
\roster 

\item "{(a)}" $(\ca L/\Theta,\D,S)$ is a locality.  

\item "{(b)}" $\ca F_S(\ca L/\Theta)=\ca F_S(\ca L)$.  

\item "{(c)}" For each $P\in\D$, the restriction 
$$ 
\r_P:N_{\ca L}(P)\to N_{\ca L/\Theta}(P)
$$ 
of $\r$ induces an isomorphism 
$$ 
N_{\ca L/\Theta}(P)\cong N_{\ca L}(P)/\Theta(P), 
$$ 
and $N_{\ca L/\Theta}(P)$ is of characteristic $p$. 

\endroster 
\endproclaim 

\demo {Proof} Let $x\in\Theta$, set $P=C_{S_x}(x)$, and set $Q=N_{S_x}(P)$. There exists 
$P_0\in\D$ with $x\in\Theta(P_0)$, whence $P_0\leq P$, and so $P\in\D$. The condition (2) implies that 
$Z(P)$ is a maximal $p$-subgroup of $C_{\ca L}(P)$, so 6.1 yields $C_{\ca L}(P)=Z(P)\times\Theta(P)$. 
As $x$ is a $p'$-element of $C_{\ca L}(P)$ we conclude that $x\in\Theta(P)$. As 
$\Theta(P)\norm N_{\ca L}(P)$ we then have 
$$ 
[Q,x]\leq QQ^x\cap\Theta(P)\leq S\cap\Theta(P)=\1.  
$$ 
Thus $Q=P$, and then $P=S_x$ by (1). By an argument similar 
to the preceding one we obtain $x\in\Theta(R)$ for all $R\in\D$ with $E\leq S_x$. Thus: 
\roster 

\item "{(*)}" Let $x\in\Theta$, and let $P\in\D$ with $P\leq S_x$. Then $x\in\Theta(P)$. 

\endroster 

Clearly, $\1\in\Theta$, and $\Theta$ is closed under inversion. Let 
$$ 
w=(x_1,\cdots,x_n)\in\bold W(\Theta)\cap\bold D, 
$$ 
and set $P=S_w$. By (*), and by induction on $n$, we 
obtain $x_i\in\Theta(P)$ for all $i$, and hence $\Pi(w)\in\Theta(P)$. Thus $\Theta$ is a partial subgroup 
of $\ca L$. Now let $x\in\Theta$ and let $g\in\ca L$ be given such that $(g\i,x,g)\in\bold D$ via 
some $Q\in\D$. Then $Q^{g\i}\leq S_x$, so (*) yields $x\in\Theta(Q^{g\i})$, and then 
$x^g\in\Theta(Q)$ by 2.3(b). This completes the proof that $\Theta\norm\ca L$. 

Set $\bar{\ca L}=\ca L/\Theta$ and adopt the usual ``bar"-convention for images of elements, 
subgroups, and collections of subgroups under the quotient map $\r:\ca L\to\bar{\ca L}$. 
Since $\Theta$ is a set of $p'$-elements of $\ca L$ we have $S\cap\Theta=\1$, and we may 
therefore identify $S$ with $\bar S$, and $\D$ with $\bar\D$. Point (a) is then given by 5.5. 

For each $P\in\D$ let $\r_P$ be the restriction of $\r$ to $N_{\ca L}(P)$. 
Then $\r_P$ is an epimorphism $N_{\ca L}(P)\to N_{\bar{\ca L}}(P)$ by 5.3(c), with kernel 
$\Theta(P)$. This yields point (c). 

By 5.3(c) the conjugation maps $c_g:P\to Q$ in $\ca F$, with $P,Q\in\D$ and with $g\in\ca L$, are the same 
as the conjugation maps $c_{\bar g}:P\to Q$ with $\bar g\in\ca L/\Theta$. Since $\ca F_S(\ca L)$ is 
generated by such conjugation maps $c_g$, we obtain $\ca F_S(\ca L)=\ca F_S(\ca L/\Theta)$. That is, (b) 
holds, and the proof is complete. 
\qed 
\enddemo

\vskip .1in 
\noindent 
{\bf Section 7: Products of partial normal subgroups} 
\vskip .1in 

The main result of this section (Theorem 7.7) is that the product of any collection of partial 
normal subgroups in a locality is again a partial normal subgroup. The proof, aside from some 
relatively minor details, is the same as that given by Ellen Henke [He] for finite localities. 
As in the finite case, the argument is based on the study of products of pairs of partial normal 
subgroups. For that reason, it will be convenient to establish the following notation. 

\definition {Hypothesis 7.1} $\ca M$ and $\ca N$ are partial normal subgroups of the locality 
$(\ca L,\D,S)$. Set 
$U=S\cap\ca N$ and $V=S\cap\ca N$. Also, set $\ca K=\ca M\cap\ca N$, and set $T=S\cap\ca K$.  
Let $\bar{\ca L}$ be the 
quotient locality $\ca L/\ca K$, and let $\r:\ca L\to\bar{\ca L}$ be the canonical projection. 
Write $\bar X$ for the image under $\r$ of a subset or element $X$ of $\ca L$, and write $\bar{\bold D}$ 
for the domain of the product in $\bar{\ca L}$. 
\enddefinition 

\proclaim {Lemma 7.2} Assume the setup of 7.1, and suppose that $\ca M\cap\ca N\leq S$. Then 
$\ca M\leq N_{\ca L}(V)$ and $\ca N\leq N_{\ca L}(U)$. 
\endproclaim 

\demo {Proof} Let $g\in\ca M$, set $P=S_g$, and let $x\in N_V(P)$. Then $(x\i,g\i,x,g)\in\bold D$ via 
$P^{gx}$, and then 
$$ 
[x,g]\in\ca M\ca N\leq S\cap\ca N=V. 
$$ 
Thus $N_V(P)\leq P$, and then $V\leq P$ and $[V,g]\leq V$. That is, we have $\ca M\leq N_{\ca L}(V)$, 
and the lemma follows. 
\qed 
\enddemo

\proclaim {Lemma 7.3} Assume the setup of 7.1, and suppose that $\ca M\cap\ca N=\1$. 
\roster 

\item "{(a)}" For each $g\in\ca M\ca N$ there exists $x\in\ca M$ and $y\in\ca N$ such that 
$(x,y)\in\bold D$, $g=xy$, and $S_g=S_{(x,y)}$. 

\item "{(b)}" $\ca M\ca N=\ca N\ca M$ is a partial normal subgroup of $\ca L$, and 
$S\cap\ca M\ca N=UV$. 

\endroster 
\endproclaim 

\demo {Proof} Let $\ca T$ be the set of all triples $(g,x,y)\in\ca M\ca N\times\ca M\times\ca N$ 
such that $g=xy$ and such that $g$ is a counter-example to (a). Among all such triples, let 
$(g,x,y)$ be chosen so that $dim(S_{(x,y)})$ is as large as possible. Set $Q=S_{(x,y)}$ and set 
$P=N_{S_g}(Q)$. It suffices to show that $P=Q$ in order to obtain obtain a contradiction and to thereby 
establish (a). 

By 4.2 we have $(y,y\i,x,y)\in\bold D$ and $g=yx^y$. Suppose that $P\leq S_y$. Then $P^y\leq S$, 
and since $P^g\leq S$ we conclude that $P\leq S_{(y,x^y)}$, and hence $P=Q$, as desired. 
Thus we may assume: 
\roster 

\item "{(1)}" $P\nleq S_y$. 

\endroster 
Let $h$ be $\up$-maximal in the maximal coset of $\ca M$ containing $g$. Then the Frattini Lemma (4.10) 
yields an element $r\in\ca M$ such that $g=rh$, and the Splitting Lemma (4.11) yields $S_g=S_{(r,h)}$. 
Then $Q\leq S_{(r,h)}$, so 
$(y\i,x\i,r,h)\in\bold D$ via $Q^g$ and $\Pi(y\i,x\i,r,h)=\Pi(g\i,g)=\1$. Thus: 
$$ 
y=x\i rh\quad\text{and}\quad h=r\i xy. \tag* 
$$ 
Since $y,h\in N_{\ca L}(U)$, it follows that $r\i x\in N_{\ca M}(U)$, and 
then that $h=(r\i x)y\in\ca M\ca N$. 

Suppose that $h$ does not provide a counter-example to (a). That is, suppose that there exists 
$x'\in\ca M$ and $y'\in\ca N$ such that  $(x',y')\in\bold D$, $x'y'=h$, and $S_{(x',y'})=S_h$. As 
$r\i xy=h=x'y'$ we get $xy=rx'y'$, and $(rx',y')\in\bold D$ with $rx'y'=rh=g$. The idea now is to 
replace $(x,y)$ with $(rx',y')$ and to contradict the assumption that $S_g\neq Q$. In order to achieve this, 
observe first of all that $S_g\leq S_r$ since $S_{(r,h)}=S_{rh}=S_g$. Then observe that 
$(S_g)^r\leq S_h$, and that $S_h=S_{(x',y')}\leq S_{x'}$. Thus $(S_g)^r\leq S_{x'}$, and so 
$S_g\leq S_{rx'}$. As $rx'y'=g$ we conclude that $S_g\leq S_{(rx',y')}$, which yields the desired 
contradiction. Thus $h$ is itself a counter-example to (a). 

Since $h=r\i xy$ by (*), and since $h$ and $y$ are in $N_{\ca L}(U)$, we have $r\i x\in N_{\ca M}(U)$, 
and then $U\leq S_{(r\i x,y)}$ since $h\in N_{\ca L}(U)$. Note furthermore that 
$Q=S_{(x,y)}\leq S_g=S_{(r,h)}\leq S_r$, and thus $Q^rU\leq S_{(r\i x,y)}$. The maximality of 
$dim(Q)$ in our initial choice of $(g,x,y)$ then yields $Q^r=Q^rU=S_{(r\i x,y)}$. Thus $U^r\leq Q^r$, 
and conjugation by $r\i$ yields $U\leq Q$. A symmetric argument yields $V\leq Q$. 
Setting $H=N_{\ca L}(Q)$, it follows from 4.1(b) that $x,y\in H$. Then $Q=O_p(H)$. 

Set $X=H\cap\ca M$ and $Y=H\cap\ca N$. Then $X,Y$, and $UV$ are normal 
subgroups of $H$, and $XY/UV$ is a $p'$-group. Set $\bar H=H/(X\cap Y)UV$. 
Here $P\leq H$ and $[P,g]\leq S$. Since $g\in XY$ we obtain 
$$
[\bar P,\bar g]=[\bar P,\bar x][\bar P,\bar y]\leq \bar X\bar Y,  
$$  
and since $\bar X\bar Y$ is a $p'$-group we get $[\bar P,\bar g]=1$. As $\bar X\cap\bar Y=1$ we have  
$C_{\bar X\bar Y}(\bar P)=C_{\bar X}(\bar P)\times C_{\bar Y}(\bar P)$. 
As $\bar g=\bar x\bar y$ it follows that $\bar x$ and $\bar y$ centralize 
$\bar P$. Thus $P^x\leq(X\cap Y)P$ and $P\in Syl_p((X\cap Y)P)$. Thus 
there exists $z\in X\cap Y$ with $P^x=P^z$. Replacing $(x,y)$ with 
$(xz\i,zy)$ we get $g=(xz\i)(zy)$ and $P\leq S_{(xz\i,zy)}$. This 
contradicts the maximality of $dim(Q)$ and (at long last) yields a contradiction which proves (a). 

In order to prove (b): let $w=(g_1,\cdots,g_n)\in\bold W(\ca M\ca N)\cap\bold D$, and set $Q=S_w$.  
By 4.11 we 
may write $g_i=x_iy_i$ with $x_i\in\ca M$, $y_i\in\ca N$, and with $S_{g_i}=S_{(x_i,y_i)}$. Set 
$w'=(x_1,y_1,\cdots,x_n,y_n)$. Then $w'\in\bold D$ via $Q$ and $\Pi(w)=\Pi(w')$. Since each $y_i$ normalizes 
$U$, it follows from 4.4 that $\Pi(w')=\Pi(w'')$ for some $w''$ such that $w''=u\circ v\in\bold D$, where 
$u\in\bold W(\ca M)$, and where $v\in\bold W(\ca N)$. Thus $\ca M\ca N$ is closed under $\Pi$. 
In order to show that $\ca M\ca N=(\ca M\ca N)\i$ we note that if $(x,y)\in\bold D\cap(\ca M\times\ca N)$ 
then $(y\i,x\i)\in\bold D$ and that $y\i x\i\in\ca M\ca N$ by 4.2. Thus $\ca M\ca N$ is a partial 
subgroup of $\ca L$. Moreover, we have shown that $\ca M\ca N=\ca N\ca M$.

Let $g\in\ca M\ca N$ and let $f\in\ca L$ with $(f\i,g,f)\in\bold D$. As usual we may write $f=hr$ with 
$r\in\ca N$, $h\in N_{\ca L}(V)$, and $S_f=S_{(h,r)}$. Write $g=xy$ as in (a). By assumption we have 
$(f\i,g,f)\in\bold D$ via some $P\in\D$. Setting $v=(r\i,h\i,x,y,h,r)$ it follows that 
$v\in\bold D$ via $P$ and that $g^f=\Pi(v)$. Here $(h\i,h,y,h)\in\bold D$ via $S_{(y,h)}$ by 4.2, so 
$v':=(r\i,h\i,x,h,h\i,y,h,r)\in\bold D$ via $P$. Then  
$$
g^f=\Pi(v)=\Pi(v')=(x^h y^h)^r\in(\ca M\ca N)^r. 
$$ 
Since $r\in\ca N$, and $\ca M\ca N$ is a partial group, we conclude that 
$g^f\in\ca M\ca N$. Thus $\ca M\ca N\norm\ca L$. 

Set $M=N_{\ca L}(S)$, $N=N_{\ca N}(S)$, and let $s\in S\cap\ca M\ca N$. Then (a) yields $s=fg$ with 
$f\in\ca M$, $g\in\ca N$, and with $S=S_{(f,g)}$. Thus $f\in M$ and $g\in N$, where $M$ and $N$ are 
normal subgroups of the group $N_{\ca L}(S)$. Then $UV$ is a normal Sylow $p$-subgroup of $MN$, and since 
$s=fg\in MN$ we obtain $s\in UV$. Thus $S\cap\ca M\ca N=UV$, and the proof of (b) is complete. 
\qed 
\enddemo 

Recall from 7.1 the notation concerning $\ca K$, $T$, and $\bar{\ca L}$. Recall also the notion from 4.5, 
of an $\up$-maximal element of $\ca L$ with respect to $\ca K$. 

\proclaim {Lemma 7.4} Assume the setup of 7.1, and let $\bar g\in\bar{\ca M}\bar{\ca N}$. Then there 
exist elements $x\in\ca M$ of and $y\in\ca N$ satisfying the following conditions. 
\roster 

\item "{(a)}" $(x,y)\in\bold D$ and $\bar g=\bar x\bar y$. 

\item "{(b)}" $x$, $y$, and $xy$ are $\up$-maximal with respect to $\ca K$. 

\item "{(c)}" $S_{xy}=S_{(x,y)}$. 

\endroster 
\endproclaim 

\demo {Proof} By 5.7, $\bar{\ca M}$ and $\bar{\ca N}$ are partial normal subgroups of the locality 
$\bar{\ca L}$, and $\bar{\ca M}\cap\bar{\ca N}=\bar{\1}$. Then 7.3(a) yields a pair of elements 
$\bar x\in\bar{\ca M}$ and $\bar y)\in\bar{\ca N}$ such that $(\bar x,\bar y)\in\bar{\bold D}$, 
$\bar g=\bar x\bar y$, and $\bar S_{\bar g}=\bar S_{(\bar x,\bar y)}$. Note that any preimage $x$ 
of $\bar x$ in $\ca L$ lies in $\ca M$ by 5.3(a) (and similarly for a preimage $y$ of $\bar y$). 
A coset $\ca Kf$ of $\ca K$ is maximal if and only if $f$ is $\up$-maximal, by 4.13(c), so 
5.3(a) implies that we may choose such preimages $x$ and $y$ so that $x$ and $y$ are $\up$-maximal with 
respect to $\ca K$. Here $x$ and $y$ are in $N_{\ca L}(S\cap\ca K)$ by 4.8, and $(x,y)\in\bold D$ by 
the definition of $\bar{\bold D}$ in 4.15. Thus (a) holds. 

We now require: 
\roster 

\item "{(1)}" $S_f\leq S_{(x,y)}$ for each $f\in\ca L$ such that $\bar f=\bar g$; and if $f\in N_{\ca L}(T)$ 
then $S_f=S_{(x,y)}$. 

\endroster 
In order to prove (1) note first of all that since $x$ and $y$ are $\up$-maximal, we have 
$\bar S_{\bar x}=\bar{S_x}$ and $\bar S_{\bar y}=\bar{S_y}$ by 5.3(c). Hence: 
$$ 
\bar S_{\bar g}=\bar S_{(\bar x,\bar y)}=\{\bar a\in\bar S_{\bar x}\mid \bar a^{\bar x}\in\bar S_{\bar y}\}=
\{\bar a\in\bar S\mid a\in S_x\ \text{and}\ \bar a^{\bar x}\in\bar S_{\bar y}\}. 
$$ 
Here $a\in S_x$ for each $a\in S$ such that $\bar a\in\bar S_{\bar x}$, so 
$$ 
\bar S_{\bar g}=\{\bar a\in\bar S\mid a\in S_x\ \text{and}\ a^x\in S_y\}=S_{(\bar x,\bar y)}.\tag 2 
$$
Now let $f\in\ca L$ such that $\bar f=\bar g$. As $T:=S\cap\ca K\leq S_{(x,y)}$, (2) yields 
$TS_f=S_{(x,y)}$. In the case that $f$ $T\leq S_f$ we get $S_f=S_{(x,y)}$, and this completes the 
proof of (1). 

We now apply (1) to the element $f=xy$, and obtain $S_{xy}=S_{(x,y)}$. This completes the proof of (c),  
and it only remains to show that $xy$ is $\up$-maximal with respect to $\ca K$. Suppose that 
$xy$ is not $\up$-maximal. Then there exists $z\in\ca K$ such that $(z,xy)\in\bold D$ and such 
that $S_{xy}$ is a proper subgroup of $S_{z(xy)}$. Then $T\leq S_{z(xy)}$. An application of (1) 
to $f=z(xy)$ then yields a contradiction, and completes the proof of (b). 
\qed 
\enddemo

\proclaim {Lemma 7.5} Assume the setup of 7.1. Let $g\in\ca L$ with $\bar g\in\bar{\ca M}\bar{\ca N}$. 
Then $g\in\ca M\ca N$, and there exist elements $x\in\ca M$ and $y\in\ca N$ with $(x,y)\in\bold D$, 
$xy=g$, and $S_{(x,y)}=S_g$. 
\endproclaim 

\demo {Proof} That $g$ is in $\ca M\ca N$ is a consequence of 5.7. 
By 7.4 we may choose $x\in\ca M$ and $y\in\ca N$ with $(x,y)\in\bold D$, $xy$ $\up$-maximal 
with respect to $\ca K$, $\bar g=\bar x\bar y$, and with $S_{xy}=S_{(x,y)}$. Then $\ca Kxy$ is a maximal 
coset of $\ca K$ by 4.13(c), and $g\in\ca Kxy$ by 5.3(a). Thus there exists $z\in\ca K$ with 
$(z,xy)\in\bold D$ and such that $z(xy)=g$. As $S_{xy}=S_{(x,y)}$ we have also 
$S_{(z,xy)}=S_{(z,x,y)}$, so that $(z,x,y)\in\bold D$ and $g=\Pi(z,x,y)$. Then $g=(zx)y$ by 
$\bold D$-associativity. As $\ca K\leq\ca M$ we get $zx\in\ca M$, and it now suffices to show that 
$S_{(zx,y)}=S_g$. As $xy$ is $\up$-maximal, the Spltting Lemma (4.12) shows that 
$S_g=S_{(z,xy)}$, and so: 
$$ 
S_g=S_{(z,x,y)}\leq S_{(zx,y)}. 
$$
The reverse inclusion $S_{(zx,y)}\leq S_g$ is immediate from $g=(zx)y$. Thus, the lemma holds with 
$zx$ in the role of $x$. 
\qed 
\enddemo

Given $g\in\ca M\ca N$ and a pair of elements $(x,y)\in\ca M\times\ca N$ satisfying the conclusion of 
lemma 7.5, we shall say that $(x,y)$ is an {\it $(\ca M,\ca N)$-decomposition} for $g$.

\proclaim {Theorem 7.7} Let $\ca M$ and $\ca N$ be partial normal subgroups of $\ca L$. Then 
the following hold. 
\roster 

\item "{(a)}" $\ca M\ca N=\ca N\ca M\norm\ca L$.

\item "{(b)}" $S\cap\ca M\ca N=(S\cap\ca M)(S\cap\ca N)$. 

\item "{(c)}" There exists an $(\ca M,\ca N)$-decomposition for each $g\in\ca M\ca N$. 

\endroster 
\endproclaim 

\demo {Proof} By 7.3(b), $\bar{\ca M}\bar{\ca N}=\bar{\ca N}\bar{\ca N}$ is a partial normal 
subgroup of $\bar{\ca L}$, and then 5.7 yields (a). Point (c) is given by 7.5. Now let 
$g\in S\cap\ca M\ca N$, and let $(x,y)$ be a corresponding $(\ca M,\ca N)$-decomposition. As 
$S_g=S$ we then have $x,y\in N_{\ca L}(S)$. Set $G=N_{\ca L}(S)$, $M=\ca M\cap G$, and 
$N=\ca N\cap G$. Further, set $U=S\cap M$ and $V=S\cap M$
Then $S=O_p(H)\in Syl_p(H)$, and we have $x\in M$ and $y\in N$. Set $K=M\cap N$, set 
$\bar G=G/K$, and adopt the usual ``bar-convention" (as in 7.1) for homomorphic images. Then 
$\bar M\bar N$ is a direct product of torsion groups. By 4.1(c), $U$ and $V$ are maximal $p$-subgroups of 
$M$ and $N$. Then (L3) implies that $\bar U$ and $\bar V$ are 
Sylow $p$-subgroups in $\bar M$ and $\bar N$, and hence $\bar U\times\bar V$ is a Sylow $p$-subgroup of 
$\bar M\times\bar N$. Also $S\cap MN$ is a maximal $p$-subgroup, and hence a Sylow 
$p$-subgroup, of $MN$. Then $\bar{S\cap MN}=\bar U\bar V$. As $S\cap K\leq U\cap V$ we obtain 
$S\cap MN=UV$, completing the proof of (b). 
\qed 
\enddemo 

\proclaim {Theorem 7.8} Let $\frak N$ be a non-empty set of partial normal subgroups of the locality 
$(\ca L,\D,S)$, and let $\ca M$ be the partial subgroup of $\ca L$ generated by the union of 
$\frak N$. Then $\ca M\norm\ca L$. Moreover, for any well-ordering $\pce$ of $\frak N$,  
$\ca M$ is the set of all $\Pi(w)$ such that 
$w=(g_1,\cdots,g_n)\in\bold D$, where each $g_i$ is in some $\ca N_i\in\frak N$, 
and where $\ca N_1\pce\cdots\pce\ca N_n$. 
\endproclaim 

\demo {Proof} Fix the well-ordering $\pce$, and 
let $\frak H$ be the set of all partial normal subgroups $\ca H$ of $\ca L$ such that, 
for some non-empty subset $\frak N_0$ of $\frak N$, $\ca H$ is is the set of all elements in $\ca L$ 
of the form $\Pi(w)$, $w=(h_1,\cdots,h_n)\in\bold D$, with $h_i\in\ca N_i\in\frak N_0$, and where 
$\ca N_1\pce\cdots\pce\ca N_n$. Then $\frak N_0\sub\frak H$. Regard $\frak N_0$ as a poset via inclusion, 
let $\L$ be a totally ordered subset of this poset, and set $\ca H=\bigcup\L$. All calculations 
with finite subsets of $\ca H$ take place in some member of $\L$, so $\ca H\norm\ca L$, and thus 
$\ca H\in\frak H$. By Zorn's Lemma there then exists a maximal $\ca M\in\frak H$. Now let $\ca N\in\frak N$. 
Then $\ca M\ca N\norm\ca L$ by 7.7, so $\ca N\leq\ca M$ by maximality. Thus $\ca M=\<\bigcup\frak N\>$. 
\qed 
\enddemo

\vskip .2in 
\noindent 
{\bf Appendix A: A class of group-theoretic examples} 
\vskip .1in 

A group $G$ is a locality for the prime $p$ if and only if:  
\roster 

\item "{(1)}" $G$ is locally finite and countable, and 

\item "{(2)}" $\Omega_S(G)$ is finite-dimensional for some maximal $p$-subgroup $S$ of $G$. 

\endroster 

Indeed, definition 3.1 shows that the conditions (1) and (2) are necessary. On the other hand, given  
(1) and (2), and taking $\D$ to be the set of all subgroups of $S$, one observes that $(G,\D)$ is an 
objective partial group and then that $(G,\D,S)$ is a locality. We note also that if (1) and (2) hold 
then, by 3.8, the maximal $p$-subgroups of $G$ are the Sylow $p$-subgroups of $G$. I.e. if  
$S$ is a maximal $p$-subgroup of $G$ then $S$ contains a $G$-conjugate of every $p$-subgroup of $G$. 

Our aim in this appendix is, first of all, to provide a large set of examples of groups which are 
localities for all primes $p$. A second aim is to show that any Sylow $p$-subgroup $S$ of $G$ is either 
nilpotent or ``discrete $p$-toral", and to observe that in the latter case $S$ has two useful properties 
in common with nilpotent groups. We begin by stating the results, before giving the definition (taken 
from [BLO07]) of discrete $p$-toral group.

\proclaim {Theorem A.1} Let $F$ be the algebraic closure of a finite field, let $G^*$ be a group having a 
faithful, finite-dimensional representation over $F$, and let $G$ be a homomorphic image of $G^*$. 
Then for each prime $p$, $G$ is a discrete locality $(G,\D,S)$, where $S$ is a maximal $p$-subgroup of $G$, 
and where $\D$ is the set of all subgroups of $S$. Moreover, if $p$ is equal to the characteristic of $F$ 
then $S$ is nilpotent, and otherwise $S$ is discrete $p$-toral.  
\endproclaim 

Discrete $p$-toral groups were, as mentioned, introduced in [BLO07]. We shall briefly review them here. The 
{\it Pr\" ufer group} $\Bbb Z/{p^\infty}$ is by definition the direct limit of the set of cyclic groups 
$\Bbb Z/(p^n)$ for $n>0$, taken with respect to the natural inclusion maps. A group $T$ is 
a {\it $p$-torus} if $T$ is isomorphic to the direct product of a finite number of copies of 
$\Bbb Z/{p^\infty}$. The set of elements $x$ in a $p$-torus $T$ such that $x^p=\1$ is an elementary 
abelian $p$-group of finite order $p^k$ for some $k$, and it follows that $k$ is the number of factors in any 
decomposition of $T$ as a direct product of Pr\" ufer groups. We refer to $k$ as the {\it rank} of $T$, 
and write $rk(T)=k$. We define the $p$-torus of rank $0$ to be the identity group. A {\it discrete p-toral 
group} is a $p$-group $P$ having a subgroup $T\leq P$ of finite index, such that $T$ is a $p$-torus 
(of some rank $k$, $k\geq 0$). Thus finite $p$-groups are special cases of discrete $p$-toral groups. 

\proclaim {Lemma A.2} Let $P$ be a discrete $p$-toral group. Then: 
\roster 

\item "{(a)}" All subgroups and all homomorphic images of $P$ are discrete $p$-toral. 

\item "{(b)}" There is a unique $p$-torus $T$ having finite index in $P$. 

\item "{(c)}" Every homomorphic image of a $p$-torus is a $p$-torus. 

\endroster 
\endproclaim 

\demo {Proof} Point (a) is given by [1.3 in BLO07]. Point (b) is vacuous if $P$ is finite, so assume 
that $P$ is infinite. Thus there exists a non-identity $p$-torus $T$ of finite index in $P$. The Pr\" ufer 
group is $p$-divisible, hence so is $T$, and it follows that $T$ has no subgroups of finite index. 
This yields (b). Now let $R$ be a $p$-torus and let $\bar R$ be a homomorphic image of $R$. Then 
$\bar R$ is discrete $p$-toral by (a), and $\bar R$ is $p$-divisible, which implies (c). 
\qed 
\enddemo

\proclaim {Lemma A.3} Let $S$ be a $p$-group, and assume that $S$ is either 
nilpotent or discrete $p$-toral. Let $P$ and $Q$ be subgroups of $S$ with $P\leq Q$. 
\roster 

\item "{(a)}" If $P<Q$ (i.e. $P$ is a proper subgroup of $Q$) then $P<N_Q(P)$. 

\item "{(b)}" If $P\norm Q$ and $P\neq 1$ then $C_P(Q)\neq 1$. 

\endroster 
\endproclaim 

\demo {Proof} Both the class of nilpotent groups and the class of discrete $p$-toral groups are closed 
with respect to subgroups, so in proving (a) and (b) we may assume to begin with that $Q=S$. Suppose 
first that $S$ is nilpotent. For (a), since $P\neq S$ there is a smallest $k$ such that $Z_k(S)\nleq P$.  
One then has $Z_k(S)\leq N_S(P)$, and so (a) holds. Set $S_0=S$ and recursively define $S_k$ for 
$k>0$ by $S_k=[S_{k-1},S]$. Then $S_n=1$ for some $n$. Set $P_k=P\cap S_k$. Then $[P_{k-1},S]\leq P_k$, 
and since $P_n=1$ we then have (b) in this case. 

Assume now that $S$ is discrete $p$-toral with maximal $p$-torus $T$, and set $P_0=P\cap T$. If 
$T\leq P$ then $P<N_S(P)$ by consideration of normalizers in the finite group $S/T$. So assume $P_0<T$. 
In proving (a) there is then no loss in assuming $S=PT$. Then $P_0\norm S$, and A.2 implies that 
$T/P_0$ is the maximal $p$-torus of $S/P_0$. Let $U$ be the subgroup of $T$ containing $P_0$, such that 
$U/P_0$ is the maximal elementary abelian subgroup of $T/P_0$. The action of the finite group $P/P_0$ on the 
the finite group $U/P_0$ yields $P_0<N_U(P)$, and so we have (a). 

Finally, assume that $P\norm S$ with $P\neq 1$. If $P_0=1$ then $[P,T]=1$, and then (b) follows from the 
action of the finite group $S/C_S(P)$ on $P$. So assume $P_0\neq 1$. There is then no loss in taking 
$P=P_0$, and we again obtain $[P,T]=1$, and thus (b). 
\qed 
\enddemo 

For brevity, we shall say that a group $G$ which is locally finite and countable is {\it lim-finite}. 
For $G$ a lim-finite group and $S$ a maximal $p$-subgroup of $G$, define $\G_S(G)$ to be the poset 
(partially ordered by inclusion) of all subgroups $D$ of $S$ such that $D$ is the intersection of a 
finite number of $G$-conjugates of $S$:  
$$ 
D=S\cap S^{g_1}\cap\cdots\cap S^{g_n}. \tag*
$$ 
Let $D$ given as in (*) with $n\geq 1$. Set $h_1=g_1\i$, and for all $i$ with $2\leq i\leq n$ set 
$h_i=g_{i-1}\i g_i$. Then $D$ is the set of all $x\in S$ such that each of the $n$ conjugates 
$x^{h_1\cdots h_i}$ is an element of $S$. That is, we have $D=S_w$ where $w=(h_1,\cdots,h_n)$. 

Let $\D$ be the set of all subgroups of $S$. Then $(G,\D,S)$ is a pre-locality, and so one also has the 
notion of $\Omega_S(G)$ given by definition 2.11, which by proposition 2.14 is the set of all $S_w$ with 
$w\in\bold W(G)$. Thus definition 2.11 is equivalent to the definition given by (*), and thus 
$\G_S(G)=\Omega_S(G)$.  

We say of any poset $\Omega$ that it is {\it finite-dimensional} if there is an upper bound to the 
lengths of monotone chains in $\Omega$. 

\definition {Definition A.4} For any prime $p$ let $\Bbb G_p$ be the class of all groups $G$ such that: 
\roster 

\item "{(1)}" $G$ is locally finite and countable, and  

\item "{(2)}" for some maximal $p$-subgroup $S$ of $G$ the poset $\Omega_S(G)$ is finite-dimensional. 

\endroster 
Write $\Bbb G$ for the intersection of the classes $\Bbb G_p$, over all primes $p$. 
\enddefinition

\proclaim {Lemma A.5} Let $p$ be a prime and let $G\in\Bbb G_p$.  
\roster 

\item "{(a)}" The maximal $p$-subgroups of $G$ are the Sylow $p$-subgroups of $G$. 

\item "{(b)}" Let $S$ be a Sylow $p$-subgroup of $D$ and let $\D$ be the set of subgroups of $S$. 
Then $(G,\D,S)$ is a locality. 

\endroster 
\endproclaim 

\demo {Proof} Let $S$ be a maximal $p$-subgroup of $G$, and let $\D$ be the set of all subgroups of $S$. 
Then $(G,\D,S)$ is a pre-locality (defined in 2.6). Here $S$ is a $p$-group, every subgroup of 
$G$ is lim-finite, and $\Omega_S(G)$ is finite-diminsional, so 
the conditions in definition 3.1 for $(G,\D,S)$ to be a locality are fulfilled. Thus 
(b) holds, and then (a) is given by 3.8. 
\qed 
\enddemo

\proclaim {Lemma A.6} The classes $\Bbb G_p$ and $\Bbb G$ are closed with respect to subgroups and 
homomorphic images. 
\endproclaim 

\demo {Proof} Fix a group $G\in\Bbb G_p$ for some prime $p$, let 
$H$ be a subgroup of $G$ and let $\bar G$ be a homomorphic image of $G$. It is obvious 
that $H$ and $\bar G$ are locally finite and countable,  and then $H$ and $\bar G$ have maximal 
$p$-subgroups.  

Let $T$ be a maximal $p$-subgroup of $H$ and let $S$ be a Sylow $p$-subgroup of $G$. Then $T\leq S^g$ 
for some $g\in G$, where $S^g$ is again a Sylow $p$-subgroup of $G$. Thus we may assume $T\leq S$ in 
proving that $H\in\w{\Bbb G}$. Let $(E_0< \cdots < E_k)$ be a monotone chain in $\Omega_T(H)$. 
Thus for each $i$ there are elements $h_1,\cdots,h_{r_i}$ of $H$ such that 
$$ 
E_i=T\cap T^{h_1}\cdots\cap T^{h_{r_i}}. 
$$ 
Set $D_i=S\cap S^{h_1}\cdots\cap S^{h_{r_i}}$. Then $(D_1<\cdots<D_k)$ is a monotone chain in 
$\Omega_S(G)$, and we therefore have $dim(\Omega_T(H))\leq dim(\Omega_S(G))$. 
Thus $\Omega_T(H)$ is finite-dimensional, and so $H\in\Bbb G$. 

For any subgroup $X$ of $G$ write $\bar X$ for the image of $X$ in $\bar G$. Let 
$(G_i)_{i=0}^\infty$ be a framing of $G$ by finite subgroups, having the property that $S\cap G_i$ is 
a Sylow $p$-subgroup of $G_i$ for all $i$. Then $\bar G$ is the union of its subgroups $\bar{G_i}$, and 
$\bar S\cap\bar{G_i}$ is a Sylow $p$-subgroup of $\bar{G_i}$ for all $i$. It follows that $\bar S$ is a 
maximal $p$-subgroup of $\bar G$. Monotone chains in $\Omega_{\bar S}(\bar G)$ pull back to monotone chains 
in $\Omega_S(G)$, so $\Omega_{\bar S}(\bar G)$ is finite-dimensional, and $\bar G\in\Bbb G$. 
\qed 
\enddemo

Our aim now is to prove the following result.  

\proclaim {Theorem A.7} Let $F$ be the algebraic closure of a finite field and let $n$ be a positive 
integer. Then $GL_n(F)\in\Bbb G$. Moreover if $F$ is of characteristic $p$ then a Sylow $p$-subgroup  
of $G$ is nilpotent, and otherwise $S$ is discrete $p$-toral.  
\endproclaim 

Theorem A.1 is immediate from A.6 and A.7, and from the closure of the class of discrete $p$-toral groups 
(and of nilpotent groups) with respect to subgroups and homomorphic images.

\vskip .1in 
We begin the proof of Theorem A.7 by considering the case where $p$ is unequal to the characteristic of $F$. 
By a {\it torus} we shall mean a direct product of finitely many copies of the multiplicative group of $F$.

\proclaim {Lemma A.8} Let $F$ be the algebraic closure of a finite field, let $G$ be the group 
$GL_n(F)$, and let $p$ be a prime such that $p\neq char(F)$. Let $H$ be the group of diagonal matrices 
in $G$, and let $W$ be the group of all permutation matrices in $G$. 
Let $T$ be the set of all $x\in H$ such that $|x|$ is a power of $p$, let $P$ be a 
Sylow $p$-subgroup of $W$, and set $S=PT$. Let $p^m$ be the exponent of $P$. Then: 
\roster 

\item "{(a)}" $S$ is a discrete $p$-toral group, with identity component $T$. 

\item "{(b)}" Let $A$ be a subgroup of $T$. Then $Z(C_G(A))$ is a torus, contained in $H$. 

\item "{(c)}" $S$ is a maximal $p$-subgroup of $G$. 

\item "{(d)}" Let $D\in\Omega_S(G)$. Then $D\cap T$ is a direct product $A\times B$, 
where $A$ is a $p$-torus and where $B$ is of exponent at most $p^m$. 

\endroster 
\endproclaim 

\demo {Proof} Any finite subgroup of the multiplicative group $F^\times$ is cyclic, and it is from this 
that one concludes that $T$ is a direct product of $n$ copies of the Pr\" ufer group $\Bbb Z_{p^\infty}$. 
As $P$ is finite, we obtain (a). 

Let $A\leq T$ be given, and let $a\in A$. There is then a $W$-conjugate $a'$ of $a$ such that $a'$ is 
in Jordan canonical form. Let $(\l_1,\cdots,\l_r)$ be an ordering of the eigenvalues of $a$, and for 
each $i$ with $1\leq i\leq r$ let $n_i$ be the multiplicity of $\l_i$ in $a$. Then $C_G(a')$ is a 
direct product $G_1\times\cdots\times G_r$ where $G_i\cong GL_{n_i}(F)$, and this yields (b). 

The group $\Bbb Z_{p^\infty}$ is isomorphic to each of its 
non-trivial homomorphic images, so $T$ is the unique maximal $p$-torus contained in $S$. Suppose that there 
exists a $p$-subgroup $\w S$ of $G$ properly containing $S$. As $\w S$ is countable and locally finite, 
$\w S$ has a framing $(\w S_k)_{n=0}^\infty$ by finite subgroups, and then every finite subgroup of 
$\w S$ is contained in some $s_k$. Choose $k$ so that $\w S_k\cap T$ contains the homocyclic abelian 
subgroup $A$ of $T$ of exponent $p^{m+1}$ and rank $n$. Set $S_k=\w S_k\cap S$. By taking $k$ sufficieently 
large we may further assume that $\w S_k\nleq S$, and so there exists $g\in\w S_k$ 
with $g\notin S$ and with $g\in N_G(S_k)$. Note that $D\cap T$ contains the elementary abelian 
subgroup $E$ of $T$ of order $p^n$, so $E$ is $g$-invariant. We have $H=Z(C_G(E))$, so $H$ is 
$g$=invariant, and thus $g\in HW$. As $S$ is a maximal $p$-subgroup of $HW$ we conclude that 
$g\in S$, and this contradiction proves (c). 

Every subgroup $X$ of $T$ is the direct product of a $p$-torus $T_X$ and a finite $p$-group $B$. Take 
$X$ to be $D\cap T$ for some $D\in\Omega_S(G)$. If $D=S$ then $D\cap T=T$ and we have 
(d) in that case. Thus we may assume $D\neq S$, and so 
$$ 
D=S\cap S^{g_1}\cap\cdots\cap S^{g_j}, 
$$ 
for elements $g_1,\cdots,g_j\in G$. Set $D_0=S$ if $j=1$, and otherwise set 
$D_0=S\cap S^{g_1}\cap\cdots\cap S^{g_{j-1}}$. By induction on $j$ we may assume that 
$D_0\cap T=A_0\times B_0$ where $A_0$ is a $p$-torus and where $B_0$ has exponent at most $p^m$. 

As $D\leq D_0$ we have $A\leq A_0$. 
Assume now that there exists an element $b\in B$ of order $p^{m+1}$, and set $e=b^{p^m}$. 
For all $i$ with $1\leq i\leq j$ set $b_i=b^{(g_i)\i}$ and $e_i=e^{(g_i)\i}$. As $b\in D$ we have $b_i\in S$ 
for all $i$, and then $e_i\in T$. Point (b) now shows that $Z(C_G(e))$ and $Z(C_G(e_i))$ are tori, 
contained in $H$, and then conjugation by $(g_i)\i$ sends the $p$-torus $T\cap Z(C_G(e))$ to 
$T\cap Z(C_G(e_i))$. Thus $T\cap Z(C_G(e))\leq D$, and then $T\cap Z(C_G(e))\leq A$ since $A$ is the 
identity component of $D$. But $b\in B$, so $e\in B$, while $A\cap B=1$. This contradiction completes 
the proof of (d). 
\qed 
\enddemo

\proclaim {Corollary A.9} Let $F$, $G$, and $S$ be as in the preceding lemma. 
Then $\Omega_S(G)$ is finite-dimensional.  
\endproclaim 

\demo {Proof} Let $\s=(D_0>\cdots>D_k)$ be a monotone chain in $\Omega_S(G)$. Then there is a refinement 
of $\s$ to a chain $\s'=(E_0>\cdots>E_{k'})$ with $E_0=S$, and where there exist elements $g_i\in G$ 
with $E_i=E_{i-1}\cap S^{g_i}$ for all $i$ with $1\leq i\leq k'$. Since we are seeking an upper bound 
for $k$, and since $k\leq k'$, we may take $\s=\s'$. 

Set $R_i=D_i\cap T$, where $T$ is the identity component of $S$, let $\tau$ be the chain 
$(T_0\geq\cdots\geq T_k)$, and set $\ell=log_p|S/T|$. Then $R_i>R_{i+1}$ for all but at most $\ell$ 
indices $i$. By A.8(d) we may write $R_i=A_i\times B_i$ where $A_i$ is a $p$-torus and where 
$B_i$ has exponent at most $p^\ell$. Let $i$ be an index with $R_i>R_{i+1}$. Then either the rank 
of $A_{i+1}$ is smaller than that of $A_i$, or else $A_i=A_{i+1}$ and $|B_i|>|B_{i+1}|$. Since the 
rank of $A_0$ is $n$, and since $|B_1|\leq p^{\ell n}$, we conclude that the length of the chain 
$\tau$ is at most $\ell+n+\ell n$. This is then an upper bound for the length of $\s$.  
\qed 
\enddemo

With A.9 we have reduced the proof of Theorem A.7 to the case where the characteristic of $F$ is $p$. 
Thus, Theorem A.7 will follow from the following result. 

\proclaim {Lemma A.10} Let $F$ be an algebraically closed field, set $G=GL_n(F)$, and let $S$ be a 
maximal unipotent subgroup of $G$. Then $\Omega_S(G)$ is finite-dimensional. 
\endproclaim 

\demo {Proof} The maximal unipotent subgroups of $G$ are all conjugate, so we may take $S$ to be 
the group of upper triangular unipotent matrices. Viewed in this way, it is then obvious that $U$ is 
isomorphic to an affine variety over $F$, and so $U$ is a hyperplane of a projective variety. 
Set $\Omega=\Omega_S(\ca L)$, and for $X\in\Omega$ write $d(X)$ for the dimension of $X$ as a 
variety. Let $X,Y\in\Omega$ with $X>Y$. By [Exercise 2.11 in Ha], $X$ and $Y$ are linear varieties 
and $d(X)>d(Y)$. As $d(S)$ is finite, the lemma follows. 
\qed 
\enddemo

\enddocument

Here the authors were unable to find a direct proof of the finite-dimensionality of $\Omega_S(G)$, based 
for example on the action of $G$ on its building. Instead, we shall rely on model-theoretic results 
concerning groups of finite Morley rank. We thank Tuna Altinel, Alexandre Borovik, and Katrin Tent for 
helping us through the arguments. 

We view groups as ``structures" for the ``language" $L$ consisting of the set $\{e,\cdot,^{-1}\}$ of symbols, 
where $e$ is a constant symbol, and where $\cdot$ and $^{-1}$ are mapping symbols (respectively binary and 
unary). We rely on [BN] for basic notions, especially 
concerning definability and uniform definability.  

\proclaim {Lemma A.10} Let $G$ be a group and let $X$ and $Y$ be definable subsets of $G$ (i.e. definable 
in terms of the language $L$). Then:  
\roster 

\item "{(a)}" $XY$ is definable.  

\item "{(b)}" The set $X^G$ of all $G$-conjugates of $X$ is uniformly definable. 

\endroster 
\endproclaim 

\demo {Proof} Let $\a$ be a formula which defines $X$ and let $\b$ be a formula which defines $Y$. Then 
the formula 
$$ 
\phi(z) = (\exists x\a(x))\wedge(\exists y\b(y))\wedge(z = x\cdot y) 
$$ 
defines $XY$, and similarly 
$$ 
\a_g(z) = (\exists x\a(x))\wedge (z=x^g) 
$$ 
is a formula which defines $X^g$. As $G$ is definable (for example via the formula $y\cdot y^{-1}= e$), 
the set of formulas $\a_g$ for $g\in G$ defines $X^G$ uniformly. 
\qed 
\enddemo

\proclaim {Lemma A.11} Let $K$ be an algebraically closed field, let $n$ be a positive integer, let  
$G$ be the group $GL_n(K)$. Then $G$ is of finite Morley rank. 
\endproclaim 

\demo {Proof} Algebraically closed fields, without additional structure beyond their arithmetic operations, 
are strongly minimal (and hence have Morley rank 1) by [5.7 in TZ]. Any structure definable within a strongly 
minimal structure has again finite Morley rank by [6.4.3 in TZ]. 
Here $G$ is definable in $K$ via a finite set of polynomial equations, and thus $G$ has finite Morley rank. 
\qed 
\enddemo 

The following result is [Lemma 2.8 in ABC]. 

\proclaim {Lemma A.12} Let $G$ be a group of finite Morley rank, and let $\ca H=\{\ca H_b\}_{b\in B}$ 
be a uniformly definable family of subgroups of $G$. Then: 
\roster 

\item "{(a)}" There is an absolute bound $n=n_{\ca H}$ on the lengths of chains of subgroups in $\ca H$. 

\item "{(b)}" The set $\bigcap\ca H$ of arbitrary intersections of subgroups from $\ca H$ is uniformly 
definable. 

\endroster 
\qed 
\endproclaim 

\proclaim {Corollary A.13} Let $G$ be a group of finite Morley rank and let $X$ be a non-empty subset of $G$. 
Then $C_G(X)$ is definable. 
\endproclaim 

\demo {Proof} Let $\ca H$ be the set of subgroups $\ca H_x=C_G(x)$ for $x\in G$. Each $\ca H_x$ is 
definable by the formula $\phi_x(z)$ given by $xz=zx$, and $\ca H$ is then uniformly definable. Now 
A.12(b) implies that each $C_G(X)$ is definable. 
\qed 
\enddemo  

We may now complete the proof of Theorem A.7. Let $K$ be an algebraically closed field, let $G$ be the 
group $GL_n(K)$, and let $S\leq G$ be the group of upper triangular unipotent matrices. The root-group 
$Z(S)$ is equal to $C_G(C_G(Z(S))$, and hence $Z(S)$ is a definable subgroup of $G$ by A.13. Every 
$G$-conjugate of $S$ is then definable (for example by A.10(b)). Since $S$ is the product of a finite 
number of root-groups, and since all root-subgroups of $G$ are conjugate to $Z(S)$, it follows from 
A.9(a) that $S$ is definable. Then $S^G$ is uniformly definable by A.10(b), and A.12 then yields the 
finite-dimensionality of $\Omega_S(G)$. The special case where $K$ is the algebraic closure of a 
finite field then yields Theorem A.7, and thereby completes the proof of A.1.

\vskip .2in 
\noindent 
{\bf Appendix B: Limits, and some colimits in the category of partial groups} 
\vskip .1in 

This appendix was inspired by some remarks of Edoardo Salati, who has shown [Sa] that the category of 
partial groups is complete (has all limits) and co-complete (has all colimits). The discussion here 
will establish somewhat less. 

By a {\it pointed set} we mean a set with a distinguished base-point, and there is then a category 
$Set^*$ of pointed sets with base-point-preserving maps. Let $Part$ be the category of partial groups. 
There is then a forgetful functor $Part\to Set^*$, given by regarding a partial group as a pointed set 
having the identity element as its base-point. 

In order to discuss limits and colimits in $Part$ (and their relation with limits and colimits in $Set^*$) 
we begin by reviewing some definitions. 

\definition {Definition B.1} Let $J$ be a small category and let $\ca C$ be a category. By a {\it $J$-shaped 
diagram in $\ca C$} we mean a covariant functor $F:J\to\ca C$. 
\enddefinition 

As always, composition of mappings will be written from left to right. 

\definition {Definition B.2} Let $F:J\to\ca C$ be a $J$-shaped diagram in $\ca C$. A {\it cone} to $F$ consists 
of an object $M$ of $\ca C$ together with a family $\phi=(\phi_X:M\to F(X))_{X\in Ob(J)}$ of $\ca C$-morphisms, 
such that for each $J$-morphism $f:X\to Y$ we have $\phi_Y=\phi_X\circ F(f)$. The cone 
$(M,\phi)$ is a {\it limit} of $F$ if for every cone $(N,\psi)$ to $F$ there exists a 
unique $\ca C$-morphism $u:N\to M$ such that $\psi_X=u\circ\phi_X$ for all $X\in Ob(J)$. 
\enddefinition

Consider now the case in which $\ca C$ is the category of sets (and mappings of sets), and let $F:J\to\ca C$ 
be a $J$-shaped diagram. If the only $J$-morphisms are identity morphisms then the direct product 
$\wh M$ of the sets $F(X)$ for $X\in Ob(J)$, together with the set $\wh\phi$ of associated projection maps 
is a limit of $F$. More generally,  let $M$ be the subset of $\wh M$ consisting of all $Ob(J)$-tuples 
$(a_X)_{X\in Ob(J)}$ such that, for 
each $J$-morphism $f:X\to Y$, we have $a_Y=(a_X)F(f)$ ($a_Y$ is equal to the image of $a_X$ under $F(f)$). 
Then $M$, together with the set $\phi$ of maps $\phi_X:M\to F(X)$ where $\phi_X$ is the restriction to $M$  
of the projection $\wh\phi_X:\wh M\to F(X)$, is a limit of $F$. If instead $\ca C$ is taken to be the 
category $Set^*$ of pointed sets, then $\wh M$ and $M$ are pointed sets (via the $Ob(J)$-tuple of base-points 
$*_X\in F(X)$) and one observes that $(M,\phi)$ is again a limit of $F$.

\proclaim {Theorem B.3} Let $Part$ be the category of partial groups, let $Set^*$ be the category of pointed 
sets, let $J$ be a small category, and let $F:J\to Part$ be a $J$-shaped diagram.  
Let $F_0:J\to Sets^*$ be the composition of $F$ with the 
forgetful functor $Part\to Sets^*$. Then there exists a limit $(\ca M,\phi)$ of $F$, and the forgetful 
functor $Part\to Set^*$ sends $(\ca M,\phi)$ to a limit of $F_0$. 
\endproclaim 

\demo {Proof} We shall only outline the steps to the proof, leaving most details to the reader. 
Let $\wh{\ca M}$ be the pointed set obtained as the direct product of the partial 
groups $F(X)$ for $X\in Ob(J)$. Let $\wh{\bold D}$ be the direct product of the pointed sets 
$\bold D(F(X))$. Thus the members of $\wh{\bold D}$ are $Ob(J)$-tuples $(w_X)_{X\in Ob(J)}$, 
with $w_X\in\bold D(F(X))$. Let $\Pi_X:\bold D(F(X))\to F(X)$ be the product. There is then a mapping 
$$ 
\wh\Pi:\wh{\bold D}\to\wh{\ca M} 
$$ 
which sends $(w_X)_{X\in Ob(J)}$ to $(\Pi_X(w_X))_{X\in Ob(J)}$. It is now straightforward to check that 
$\wh{\ca M}$ is a partial group via the product $\wh\Pi$ and via the inversion map which sends an 
element $(g_X)$ of $\wh{\ca M}$ to the $Ob(J)$-tuple $(g_X\i)$ of inverses. 

Let $\ca M$ be the subset of $\wh{\ca M}$ consisting of all $Ob(J)$-tuples 
$(g_X)_{X\in Ob(J)}$ such that, for each $J$-morphism $f:X\to Y$, we have $g_Y=(g_X)F(f)$. Let 
$\phi$ be the $Ob(J)$-tuple $(\phi_X)$ of maps $\phi_X:\ca M\to F(X)$ obtained by restriction to $\ca M$ 
of the projection $\wh\phi_X:\wh{\ca M}\to F(X)$. One observes that each $\phi_X$ is a homomorphism 
of partial groups, and that $(\wh{\ca M},\phi)$ is a cone of $F$. 

Now let $(\ca N,\psi)$ be any cone of $F$. For $g\in\ca N$ define $g\mu$ to be the $Ob(J)$-tuple 
$(g\psi_X)$. One checks that each such $g\mu$ is an element of $\ca M$, and then that the mapping 
$\mu:\ca N\to\ca M$ is a homomorphism of partial groups. Finally, one observes that 
$\psi_X=u\circ\phi_X$ for all $X\in Ob(J)$ and that $u$ is necessarily the unique homomorphism 
$\ca N\to\ca N$ having this property. Thus $(\ca M,\phi)$ is a limit of $F$. But also
$(\ca M,\phi)$ is a limit of $F_0$, which completes the proof. 
\qed 
\enddemo

The situation for colimits of partial groups is not as straightforward as that of limits.  

\definition {Definition B.4} Let $F:J\to\ca C$ be a $J$-shaped diagram in $\ca C$. A {\it co-cone} to $F$ 
consists of an object $M$ of $\ca C$ together with a family $\phi=(\phi_X:F(X)\to M)_{X\in Ob(J)}$ of 
$\ca C$-morphisms, such that for each $J$-morphism $f:X\to Y$ we have $\phi_X=F(f)\circ\phi_Y$. The co-cone 
$(M,\phi)$ is a {\it colimit} of the diagram $F$ if for every co-cone $(N,\psi)$ to $F$ there exists a 
unique $\ca C$-morphism $u:M\to N$ such that $\psi_X=\phi_X\circ u$ for all $X\in Ob(J)$. 
\enddefinition 

Again, it will be fruitful to review the case where $\ca C$ is the category of sets. 
Thus, let $F:J\to\ca C$ (with $\ca C=Sets$) be a $J$-shaped diagram. Let $\wh M$ be the disjoint union of 
the sets $F(X)$ for $X\in Ob(J)$. Let $\sim$ be 
the relation on $\wh M$ given by $a\sim b$ if there exists a $J$-morphism $f:X\to Y$ such that 
$a\in F(X)$, $b\in F(Y)$, and $b=(a)F(f)$ is the image of $a$ under $F(f)$. Let $\approx$ be the 
symmetrization of $\sim$ (so that $a\approx b$ if either $a\sim b$ or $b\sim a$). 
As $\sim$ is transitive, $\approx$ is then an equivalence relation.  
 Let $M$ be the set $\wh M/\approx$ of equivalence classes and 
let $\phi$ be the set of all $\phi_X:F(X)\to M$, where $\phi_X$ 
is the mapping which sends $a\in F(X)$ to the $\approx$-equivalence class of $a$ in $\wh M$.  
Then $(M,\phi)$ is a colimit of $F$. 

Next, take $\ca C$ to be the category $Set^*$ and let $F:J\to\ca C$ be a $J$-shaped diagram in $\ca C$. 
Here we take $\wh M$ to be the pointed set obtained from the disjoint union of the pointed sets $F(X)$ 
(over all $X\in Ob(J)$) by identifying base-points. For any $J$-morphism $f:X\to Y$, the 
morphism $F(f)$ of pointed sets sends the base-point of $F(X)$ to the base-point of $F(Y)$, and we may 
therefore define the equivalence relation $\approx$ on $\wh M$ as in the preceding paragraph. Take 
$M=\wh M/\approx$. Again, 
for each $X\in Ob(J)$ one has the mapping $\phi_X:X\to M$ which sends $a\in X$ to the 
$\approx$-equivalence class of $a$ in $\wh M$, and $(M,\phi)$ is a colimit of $F$.

In passing now to the case where $\ca C$ is the category of partial groups, we face the problem that, 
in general, there will be no partial normal subgroup, and no ``quotient" partial group corresponding to the 
equivalence relation $\approx$. For this reason, restrictions must be placed on the sort of diagrams  
$F:J\to Part$ that may be considered.

\proclaim {Theorem B.5} Let $Part$ be the category of partial groups, let $Set^*$ be the category of pointed 
sets, and let $J$ be a small category, and let $F:J\to Part$ be a $J$-shaped diagram. Assume: 
\roster 

\item "{(1)}" For each  ordered pair $(X,Y)$ of objects of $J$, there exists at most one $J$-morphism 
$X\to Y$.  

\item "{(2)}" For each $J$-morphism $f:X\to Y$ the kernel of the homomorphism $F(f):F(X)\to F(Y)$ of 
partial groups is trivial. 

\endroster 
Let $F_0:J\to Sets^*$ be the composition of $F$ with the 
forgetful functor. Then there exists a colimit $(\ca M,\phi)$ of $F$, and the forgetful 
functor sends $(\ca M,\phi)$ to a limit of $F_0$. 
\endproclaim 

\demo {Proof} Let $\wh{\ca M}$ be the pointed set obtained as the disjoint union of all of the partial 
groups $F(X)$, for $X\in Ob(J)$, and with base-points identified. Define $\wh{\bold D}$ to be the 
disjoint union of the domains $\bold D(F(X))$. There is then a mapping $\wh\Pi:\wh{\bold D}\to\wh{\ca M}$ 
whose restriction to 
$\bold D(F(X))$ is the product $\Pi_X$ on $F(X)$. The union of the inversion maps on the partial groups 
$F(X)$ is an involutory bijection on $\wh{\ca M}$, and one may check that $\wh{\ca M}$ is a partial 
group via these structures. 

Let $\sim$ be the equivalence relation on $\wh{\ca M}$ given by $a\sim b$ if there exists a 
$J$-morphism $f:X\to Y$ such that $a\in F(X)$, $b\in F(Y)$, and $F(f):a\maps b$. Let $\approx$ be the 
equivalence relation given by symmetrizing $\sim$, and extend $\approx$ to an equivalence relation on 
$\bold W(\wh{\ca M})$ in the component-wise way. That is, if $u=(a_1,\cdots,a_m)$ and 
$v=(b_1,\cdots,b_n)$ are words in the alphabet $\wh{\ca M}$ then $u\approx v$ if and only if $m=n$ and 
$a_i\approx b_i$ for all $i$. For $a\in\wh{\ca M}$ we write $[a]$ for the $\approx$-class of $a$. Then 
the $\approx$-class of a word $u=(a_1,\cdots,a_n)$ is the word $([a_1],\cdots,[a_n])$.
Let $\ca M$ be the pointed set $\wh{\ca M}/\approx$ (whose base-point is the 
equivalence class of the base-point of $\wh{\ca M}$, and let $\bold D$ be the set of all words 
$([a_1],\cdots,[a_n])$ having a representative $(a_1,\cdots,a_n)\in\wh{\bold D}$.]  

Let $u=(a_1,\cdots,a_n)\in\wh{\bold D}$, and assume that there exists at least one index $k$ such that 
$a_k$ is not the identity element of $\wh{\ca M}$. Then there is a unique object $X$ of $J$ such that 
$u\in\bold D(F(X))$. Let also $v=(b_1,\cdots,b_n)\in\wh{\bold D}$, and assume $u\approx v$. Then (2) 
implies that $b_k$ is not the identity element of $\wh{\ca M}$, and there is a unique object 
$Y$ of $J$ with $v\in\bold D(F(Y))$. Let $i$ be any index from $1$ to $n$ such that not 
both $a_i$ and $b_i$ are identity elements. Then neither $a_i$ nor $b_i$ is an identity element, 
and (1) implies that either there is a unique $J$-morphism $f:X\to Y$ and $F(f):a_i\maps b_i$, or there 
is a unique $J$-morphism $g:Y\to X$ and $F(g):b_i\maps a_i$. If there exist both a $J$-morphism 
$f:X\to Y$ and a $J$-morphism $g:Y\to X$ then (1) implies that $f$ and $g$ are isomorphisms, and are 
inverse to each other, whence $F(f):a_i\maps b_i$ if and only if $F(g):b_i\maps a_i$. We may 
therefore assume without loss of generality that there exists a $J$-morphism $f:X\to Y$ and that 
$F(f)$ maps $u$ to $v$ component-wise. As $F(f)$ is a homomorphism $F(X)\to F(Y)$ of partial groups  
we then have $\wh\Pi(u)=\wh\Pi(v)$. We have thus shown that $\wh\Pi$ induces a mapping 
$\Pi:\bold D\to\ca M$, and the reader may check that $\Pi$ is a product, as defined in 1.1. 
If $\l:\ca H\to\ca K$ is a homomorphism of partial groups and $a\in\ca H$, then $(a\i)\l =(a\l)\i$, 
so there is a well-defined inversion mapping $\ca M\to\ca M$ given by $[a]\i=[a\i]$. Again it is left to 
the reader to check that with these structures, $\ca M$ is a partial group. 

For each $X\in Ob(J)$ define $\phi_X:F(X)\to\ca M$ by $a\phi_X=[a]$. Then $\phi_X$ is a homomorphism. If 
$a'\in F(X)$ with $[a]=[a']$ then $a=a'$, since the only $J$-morphism $X\to X$ is the identity morphism. 
Thus $\phi_X$ is injective. For any $J$-morphism $f:X\to Y$ and any $a\in F(X)$ we have $[(a)F(f)]=[a]$, 
so $(\ca M,\phi)$ is a co-cone of $F$. 

Let $(\ca N,\psi)$ be an arbitrary co-cone of $F$. Thus $\psi_X=F(f)\circ\psi_Y$ whenever 
$f:X\to Y$ is a $J$-morphism. That is, we have $a\psi_X=b\psi_Y$ if $F(f):a\maps b$, and thus 
there is a well-defined mapping $\s:\ca M\to\ca N$ given by $[a]\maps[a\psi_X]$ for $a\in F(X)$. 
Moreover, we have $\phi_X\circ\s=\psi_X$, and $\s$ is the unique such mapping $\ca M\to\ca N$. One 
checks that $\s$ is a homomorphism of partial groups, in order to complete the proof that 
$(\ca M,\phi)$ is a colimit of $F$. 
\qed 
\enddemo

\vskip .2in 
\noindent 
{\bf Section 8: Compact localities and $p$-local compact groups} 
\vskip .1in 

Let $(\ca L,\D,S)$ be a pre-locality. There is then a category $\bold T(\ca L)$ whose set of objects is $\D$, 
and where $Hom_{\bold T(\ca L)}(P,Q)$ is the set of all triples $(P,g,Q)$ such that 
$P\in\D$, $g\in\ca L$ with $P\leq S_g$, and $P^g\leq Q$. Composition of morphisms is then given by: 
$$ 
(P,g,Q)\circ(Q,h,R)=(P,gh,R). 
$$   
The morphisms $(P,\1,Q)$ ($P,Q\in\D$ with $P\leq Q$) are called {\it inclusion morphisms}. 
The category $\bold T(\ca L)$ is called the {\it transporter category} of $(\ca L,\D,S)$, and it 
has the following properties. 
\roster 

\item "{(1)}" All morphisms in $\bold T(\ca L)$ are monic and epic (i.e. are subject to both left and right 
cancellation). 

\item "{(2)}" Each $\bold T(\ca L)$-morphism can be factored in a unique way as a 
$\bold T(\ca L)$-isomorphism followed by an inclusion morphism. 

\endroster 
If $\phi$ and $\psi$ are $\bold T(\ca L)$-isomorphisms, then we say that $\psi$ is an {\it extension of 
$\phi$} and write $\phi\up\psi$, if there is a commutative diagram: 
$$ 
\CD 
\w P@>\psi>>\w Q \\ 
@AAA        @AAA \\ 
P   @>>\phi> Q 
\endCD 
$$ 
in which the vertical arrows are inclusion morphisms. We say also that $\phi$ is a {\it restriction} of 
$\psi$, and write $\psi\down\phi$. 

Let $\approx$ be the equivalence relation generated by $\up$ and $\down$ on the set of 
$\bold T(\ca L)$-isomorphisms. Thus $\phi\approx\psi$ if there is a finite sequence $(\phi_0,\cdots,\phi_n)$ 
of $\bold T(\ca L)$-isomorphisms such that $\phi=\phi_0$, $\phi_n=\psi$, and such that $\phi_k$ is either 
an extension or a restriction of $\phi_{k-1}$ for all $k$ from $1$ to $n$. This leads to a third property 
of $\bold T(\ca L)$:  
\roster 

\item "{(3)}" Each $\approx$-class $C$ contains a unique isomorphism $\psi$ such that $\phi\up\psi$ 
for all $\phi\in C$. (In brief: $C$ contains a unique {\it maximal} representative.) 

\endroster 
The proof of (3) consists of the observation that if $\phi=(P,g,Q)$ then every member of the 
$\approx$-class of $\phi$ is a restriction of $\psi=(S_g,g,S_{g\i})$. 
\vskip .1in 
One further point concerning the category $\bold T_\D(\ca L)$ is that it contains $\bold T_\D(S)$ as a 
full sub-category, where $\bold T_\D(S)$ is formed in the preceding way from the pre-locality $(S,\D,S)$.  

\vskip .1in
Now let $\ca T$ be a category such that $Ob(\ca T)=\D$, and assume that there is given a faithful functor 
$$ 
\d:\bold T_\D(S)\to\ca T, 
$$ 
such that the image of $\d$ is a full subcategory of $\ca T$. The image of an 
inclusion morphism $(P,\1,Q)$ will then be called an {\it inclusion 
morphism in $\ca T$}, and we may then define relations $\up$, $\down$, and $\approx$ on the set of 
$\ca T$-isomorphisms, exactly as above. 

\definition {Definition 8.1} Let $\ca T$ be a small category, and let $\d:\bold T_\D(S)\to\ca T$ be a 
faithful functor whose image is a full subcategory. 
Then $\ca T$ is a {\it transporter system on $\D$} if the following three conditions 
hold.  
\roster 

\item "{(1)}"  All morphisms in $\bold T(\ca L)$ are monic and epic. 

\item "{(2)}" Each $\ca T$-morphism can be factored in a unique way as a $\ca T$-isomorphism followed by 
an inclusion morphism. 

\item "{(3)}" Each $\approx$-class $C$ of $\ca T$-isomorphisms contains a unique maximal representative. 

\endroster 
\enddefinition 

\proclaim {Proposition 8.2} Let $\d:\bold T(S)\to\ca T$ be a transporter system on $\D$, let 
$\ca L$ be the set of $\approx$-classes of $\ca T$-isomorphisms, and let $\bold D$ be the set of 
all words $w=(C_1,\cdots,C_n)\in\bold W(\ca L)$ (including the empty word) such that there exists a 
composable sequence $\g=(\phi_1,\cdots,\phi_n)$ of representatives $\phi_i\in C_i$. For any such $\g$ 
let $\Pi(\g)$ be the $\approx$-class of $\phi_1\circ\cdots\circ\phi_n$, and set 
$\g\i=(\phi_1^{-1},\cdots,\phi_n^{-1})$. Then the following hold. 
\roster 

\item "{(a)}" Both $\Pi(\g)$ and $\g\i$ are independent of the choice of the composable sequence $\g$ 
representing $w$. $\Pi$ may thereby be viewed as a mapping $\bold D\to\ca L$, and $(-)^{-1}$ as a 
mapping $\bold D\to\bold D$. 

\item "{(b)}" $\ca L$ is a partial group relative to $\Pi$ and $(-)^{-1}$. 

\item "{(c)}" For any $g\in S$, the image of $(S,g,S)$ under $\d$ is the unique maximal element in its 
$\approx$-class, and $S$ may thereby be viewed as a subgroup of $\ca L$. Moreover, $\D$ is then a set of 
subgroups of $\ca L$, $(\ca L,\D,S)$ is a pre-locality, and there is a commutative diagram of 
functors: 
$$ 
\CD 
\endCD
$$

natural isomorphism of categories: 
$$ 
\ca T\to\bold T(\ca L). 
$$ 

\endroster 
\endproclaim

\enddocument 

\vskip .2in 
\noindent 
{\bf Appendix: Compact localities and algebraic localities} 
\vskip .1in 

The aim of this appendix is to indicate that all of the known examples of structures related to 
``saturated" fusion systems may be viewed as localities, and to indicate as well the existence of 
further such structures. The treatment here will of necessity be expository in nature, as we have not 
as yet introduced the notion of saturation for the fusion systems of localities. Saturation, and 
related notions, will in fact be the focus of Part II of this series. 

We expect readers of this appendix to be familiar with centric linking systems over finite $p$-groups (also 
called ``$p$-local finite groups"). These were shown, in the Appendix to [Ch1], to be essentially the 
``same" as localities $(\ca L,\D,S)$ in which $S$ is a finite $p$-group, and where $\D$ is the set of 
all subgroups $P$ of $S$ such that $C_{\ca L}(Q)\leq Q$ for all $\ca F_S(\ca L)$-conjugates $Q$ of $P$.  
The expository paper [AO] by Michael Aschbacher and Bob Oliver provides a nice introduction to linking 
systems, and the book [AKO] provides a thorough treatment. Most of the terminology to be employed below 
is taken from those sources, or (in the case of infinite $p$-groups) from the papers [BLO1] and [BLO2] by 
Broto, Levi, and Oliver. Our citing of these references, which in many ways provided the inspiration 
and motivation for this work, should not disguise the fact that our aim is to supersede them. 

\definition {Definition A1} Let $(\ca L,\D,S)$ be a pre-locality. The {\it transporter system} associated 
with $\ca L$ is the category $\bold T=\bold T(\ca L)$ whose set of objects is $\D$, and in which the 
set $Mor_{\bold T}(P,Q)$ of morphisms is the set of triples $(P,g,Q)$ such that $g\in\ca L$, 
$P\leq S_g$, and $P^g\leq Q$. Composition of morphisms is then given by 
$$ 
(P,g,Q)\circ(Q,h,R)=(P,gh,R), 
$$ 
and the identity morphism in $Aut_{\bold T}(P)$ is the triple $(P,\1,P)$. 
\endproclaim 

One observation that one may immediately make is that $Aut_{\bold T}(P)$ is essentially the same 
thing as $N_{\ca L}(P)$, for $P\in\D$. 

\vskip .1in 
The general notion of transporter system over a finite $p$-group $S$ was introduced by Oliver and Ventura 
in [OV]. The definition is somewhat complicated, and will not be presented here. It will be sufficient 
for the discussion here, to know that a transporter system on a fiis a category 

whose set $\D$ of objects 

a rough idea may 
be given as follows.

The {\it Pr\" ufer group} $\Bbb Z/(p^\infty)$ is the subgroup of the circle group $\{z\in\Bbb C\mid |z|=1\}$ 
consisting of those elements $z$ whose order is finite and is a power of $p$. One observes that 
$\Bbb Z/(p^\infty)$ the nested union of subgroups isomorphic to $\Bbb Z/(p^n)$ for $n\geq 1$, and is thus 
countable and locally finite. Following [BLO1], we say that the direct product $T$ of a finite number of 
copies of $\Bbb Z/(p^\infty)$ is a {\it discrete $p$-torus}, and that a {\it discrete $p$-toral group} is 
a group $P$ having a discrete $p$-torus of finite index equal to a power of $p$. More generally, we shall 
say that a group $G$ having a discrete $p$-torus of finite index is {\it virtually $p$-toral}. Thus, 
all virtually $p$-toral groups are countable and locally finite.  

\definition {Definition 8.1} Let $(\ca L,\D,S)$ be a pre-locality. Then $\ca L$ is {\it pre-compact} if: 
\roster 

\item "{(CL1)" $S$ is a discrete $p$-toral group, and $N_{\ca L}(P)$ is virtually $p$-toral for 
all $P\in\D$. 

\endroster 
We say that $\ca L$ is a {\it compact locality} if also: 
\roster

\item "{(CL2)}" There are finitely many orbits for the action of $S$, by conjugation, on 
$\Omega_S(\ca L)$. 

\endroster 
\enddefinition 

Notice that (CL2) implies that $dim(\Omega_S(\ca L))$ is finite. Thus, compact ``localities" are 
bona fide localities in the sense of definition 3.1. 

\vskip .1in 
There is a well-developed theory of ``$p$-local compact groups" and, more generally, of ``transporter 
systems" over discrete $p$-toral groups, due to Broto, Levi, and Oliver (see [BLO1] and [BLO2]). 
The aim of this section is to establish that there is an equivalence between the notion of pre-compact 
locality and that of transporter system, and a near equivalence (in a sense which will be made precise) 
between the notion of compact locality and that of $p$-local compact group. 

A model for establishing the sort of equivalence that we have in mind is provided by the Appendix to [Ch1], 
where it is shown that localities over finite $p$-groups are the essentially the same as transporter 
systems over finite $p$-groups. We shall take the 
liberty here of giving outlines, rather than complete proofs, of statements of which the proofs are 
word-for-word identical to those in [Ch1]. We'll explain where (and why) more details are needed, and 
provide complete proofs in those instances.  

\vskip .1in 
Here then, is an outline of the correspondence that we have in mind.

\proclaim {Example 8.2} Let $G$ be a linear algebraic group defined over the algebraic closure of 
the prime field of order $p$, and let $S$ be a maximal unipotent subgroup of $G$. Then $(G,S)$ is a 
lim-finite pair. 
\endproclaim 

\demo {(Proof)} It's enough to consider the case where $G$ is reductive, and, if one prefers, 
where $G$ is simple. It is well known that that $S$ is a maximal $p$-subgroup of $G$. [REFERENCE ?] 

Fix a root system $\Phi$ for $G$ such that $S$ is generated by the root-groups $U_\a$ for $\a\in\Phi^+$, 
and let $\Pi$ be the set of fundamental roots for $\Phi^+$. A subset $X$ of 
$\Phi^+$ is {\it convex} if the set $\Bbb Z^+[X]$ (of linear combinations of $X$ with coefficients 
in $\Bbb Z^+$) intersects $\Phi^+$ in $X$. For any convex $X\sub\Phi^+$ the Chevalley commutator 
formula shows that the group 
$$ 
U_X=\<U_\a\mid\a\in X\> 
$$ 
is the product (in any order) of the root-groups $U_\a$ for $\a\in X$. Let 
$\ca U$ be the set of all such $U_X$, and set 
$$ 
\w{\ca U}=\{U^s\mid s\in S,\ U\in\ca U\}. 
$$ 

The aim is to show that $\Omega_S(G)\sub\w{\ca U}$. This will require two preliminary steps. 
\vskip .1in 
\roster 

\item "{(1)}" Let $g\in G$. Then $S\cap S^g\in\w{\ca U}$. 

\item "{(2)}" Let $X,Y$ be convex subsets of $\Phi^+$ and let $s,t\in S$. Then
$U_X^s\cap U_Y^t\in\w{\ca U}$. 

\endroster 
Assume (1) and (2), and let $\bold w=(g_1,\cdots,g_n)\in\bold W(G)$. For each $k$ with $1\leq k\leq n$ set 
$h_k=(g_1\cdots g_k)\i$. Then 
$$ 
S_{\bold w}=S\cap S^{h_1}\cap\cdots\cap S^{h_n}=(S\cap S^{h_1})\cap\cdots\cap(S\cap S^{h_n})\in\w{\ca U} \tag*
$$ 
by (1) and (2) and induction on $n$. Thus: the members of $\Omega_S(G)$ are closed and connected. One 
knows that $dim_F(U)<dim_F(V)$ whenever $U$ and $V$ are closed connected subsets of $G$ with $U$ properly 
contained in $V$. Thus $dim(\Omega_S(G))$ is bounded above by $dim_F(S)$. 

It remains to prove (1) and (2). Point (1) is a consequence of the standard decomposition $G=BNB$ of 
$G$ as a $BN$-pair (relative to $\Phi$ and a choice of a fundamental system of positive roots). Point 
(2) follows from the Chevalley commutator relations. 

In the case that $F$ is contained in the algebraic closure of a finite field one observes that 
$G$ is locally finite and countable. Then $G$ is lim-finite, and $G\in\frak G$. 
\qed 
\enddemo

\enddocument